%% file: Diagrammatic-KL-R-polynomials..tex
\documentclass[preprint,12pt]{elsarticle}

\usepackage[margin=1in,footskip=0.5in]{geometry}
\usepackage{amsbsy}
\usepackage{amssymb}
\usepackage{amsthm}
\usepackage{amsmath}
 \usepackage[usenames,dvipsnames]{pstricks}
 \usepackage[utf8]{inputenc}
\usepackage{epsfig}
\usepackage{pst-grad} 
\usepackage{pst-plot} 
\usepackage[space]{grffile} 
\usepackage{etoolbox} 
\makeatletter 
\patchcmd\Gread@eps{\@inputcheck#1 }{\@inputcheck"#1"\relax}{}{}
\makeatother

\journal{ }

\theoremstyle{plain}
\newtheorem{teo}{Theorem}[section]
\newtheorem{coro}[teo]{Corollary}

\newtheorem{defi}[teo]{Definition}
\newtheorem{conj}[teo]{Conjecture}
\newtheorem{exa}[teo]{Example}
\newtheorem{lem}[teo]{Lemma}
\newtheorem{propos}[teo]{Proposition}

\newtheorem{rem}[teo]{Remark}

\usepackage{hyperref}

\newenvironment{demo}{\noindent \textit{Proof:}  \rm}{\quad \hfill $\square$}

\input{pictures.tex}

\begin{document}

\begin{frontmatter}


\author{David Plaza\corref{cor1 dnasidnadinada dansdiadn}}
\ead{dplaza@inst-mat.utalca.cl}

\title{Diagrammatics for Kazhdan--Lusztig $\tilde{R}$-polynomials}



\address{Instituto de Matem\'atica y F\'isica\\ Universidad de Talca \\Talca, Chile}

\begin{abstract}
Let $(W,S)$ be an arbitrary Coxeter system. We  introduce a family of polynomials, $\{ \tilde{\mathcal{R}}_{u,\underline{v}}(t)\}$, indexed by pairs $(u,\underline{v})$ formed by an element $u\in W$ and a (non-necessarily reduced) word $\underline{v}$ in the alphabet $S$. The polynomial $\tilde{\mathcal{R}}_{u,\underline{v}}(t)$ is obtained by considering a certain subset of Libedinsky's light leaves associated to the pair $(u,\underline{v})$. Given a reduced expression $\underline{v}$ of an element $v\in W$; we show that $\tilde{\mathcal{R}}_{u,\underline{v}}(t)$ coincides with the Kazhdan--Lusztig $\tilde{R}$-polynomial $\tilde{R}_{u,v}(t)$. Using the diagrammatic approach, we obtain some closed formulas for $\tilde{R}$-polynomials.
\end{abstract}


\end{frontmatter}

\section{Introduction}
\label{intro}
In their landmark paper \cite{kazhdan1979representations}, Kazhdan and Lusztig introduced, for each Coxeter system $(W,S)$, a family of polynomials with integer coefficients indexed by pairs of elements of $W$. These polynomials are now known as the Kazhdan--Lusztig (KL)-polynomials 
 that are related closely with the representation theory of semisimple algebraic groups, the topology of the Schubert varieties, the Verma modules, etc. (see e.g. \cite{brenti2003kazhdan} and references therein). To demonstrate the existence of these polynomials, the authors in \cite{kazhdan1979representations} introduced another family of polynomials, the KL $R$-polynomials, which we denote by $R_{u,v}(t)\in \mathbb{Z}[t]$.  The relevance of the $R$-polynomials lies in the fact that the knowledge of them implies 
 the knowledge of the KL-polynomials.

Recently, it was proven that the KL-polynomials have non-negative coefficients \cite{elias2014hodge}. Meanwhile, the simplest examples reveal that $R$-polynomials can have negative coefficients. Fortunately, there exists another family of polynomials with positive coefficients that, after an easy change in variable, ``coincides'' with the $R$-polynomials (see Proposition \ref{propose decomposition R tilde polynomials}). These polynomials are known as the KL $\tilde{R}$-polynomials. We denote them by $\tilde{R}_{u,v}(t)$, for $u,v\in W$. 

The primary object of the study herein is a new family of polynomials, $\{\tilde{\mathcal{R}}_{u,\underline{v}}(t) \}$, indexed by pairs $(u,\underline{v})$ formed by an element $u\in W$ and a (non-necessarily reduced) word $\underline{v}$ in the alphabet $S$. We refer to these polynomials as the ``diagrammatic $\tilde{R}$-polynomials.'' Such a name is justified by the following: Given a pair $(u,\underline{v})$, we define a set, $\mathbb{L}_{\underline{v}}(u)$, which is a subset of Libedinsky's light leaves associated to the pair $(u,\underline{v})$ \cite{libedinsky2008categorie}. In \cite{elias2016soergel}, Elias and Williamson introduced a diagrammatic version of the set $\mathbb{L}_{\underline{v}}(u)$. We define the diagrammatic $\tilde{R}$-polynomial, $\tilde{\mathcal{R}}_{u,\underline{v}}(t)$,  by 
\begin{equation}
	\tilde{\mathcal{R}}_{u,\underline{v}}(t)= \sum_{l\in\mathbb{L}_{\underline{v}}(u)} t^{\deg{(l)}},
\end{equation}
where $\deg(\cdot )$ is a certain statistic defined over the set $\mathbb{L}_{\underline{v}}(u)$.
 On the other hand, we demonstrate that if $\underline{v}$ is a reduced expression of some $v\in W$, then 
\begin{equation}
	\tilde{\mathcal{R}}_{u,\underline{v}}(t)=\tilde{R}_{u,{v}}(t).
\end{equation} 
 
This last equality is the primary result herein (see Theorem \ref{teo diagrammatic R poly}), and it provides a combinatorial interpretation for $\tilde{R}$-polynomials. Obtaining a combinatorial interpretation for $\tilde{R}$-polynomials is not new; further, two of such interpretations exist in the literature: the one given by Deodhar \cite[Theorem 1.3]{deodhar1985some} in terms of distinguished subexpressions, and the one provided by Dyer \cite[Corollary 3.4]{dyer1993hecke} in terms of directed paths in the Bruhat graph of $W$. 

Our approach presents two advantages with respect to Dyer's or Deodhar's approach. On the one hand, it allows us to obtain some closed formulas for $\tilde{R}$-polynomials in a simple and combinatorial manner. In particular, we recover and generalize the closed formulas obtained previously by Marietti \cite{marietti2002closed} and Pagliacci \cite{pagliacci2001explicit}.
Here, we use the term ``combinatorial'' in its more strict meaning, i.e., we obtain the closed formulas for $\tilde{R}$-polynomials using bijective proofs. On the other hand, our approach establishes an explicit connection between $\tilde{R}$-polynomials and the theory of Soergel bimodules. 

Although not necessary, we briefly explain the aforementioned connection. As we indicated, KL-polynomials have non-negative coefficients. Elias and Williamson's proof of this fact relies on the previous work of Soergel \cite{soergel1992combinatorics,soergel2007kazhdan}. In fact, what they proved in \cite{elias2014hodge} is a deeper result known as ``Soergel's conjecture.'' One of the consequences of Soergel's conjecture is that the KL-polynomials are the graded ranks of certain Hom-spaces; in particular, they have positive coefficients. The diagrams that we use throughout this paper represent morphisms in the category of Soergel bimodules. In particular, the set $\mathbb{L}_{\underline{v}}(u)$ spans a submodule of a quotient of a certain Hom-space in the category of Soergel bimodules, such that our combinatorial interpretation of $\tilde{R}$-polynomials becomes a categorical interpretation. Namely, the diagrammatic $\tilde{R}$-polynomials are the graded rank of a submodule of a quotient of a certain Hom-space. We expect that the interpretation above together with the Hodge theoretic machinery developed in \cite{elias2014hodge} allow for a better understanding of the $\tilde{R}$-polynomials. Hereinafter, we treat the diagrams as only combinatorial objects. The reader interested in the algebraic part of the story is referred to \cite{elias2010diagrammatics,elias2016soergel,libedinsky2008categorie,libedinsky2015light,libedinsky2017gentle}.

This paper is organised as follows. In the next section, we recall the definition and some properties of the $\tilde{R}$-polynomials. In section \ref{section diagrams}, we introduce diagrammatic $\tilde{R}$-polynomials and prove that they coincide with the classical ones. Finally, in section \ref{section closed}, we provide some closed formulas for diagrammatic $\tilde{R}$-polynomials. 

\section*{Acknowledgements}
This work is supported by the Fondecyt project 11160154 and the Inserci\'on en la Academia project PAI-Conicyt 79150016. The author would like to thank Nicol\'as Libedinsky and Paolo Sentinelli for their comments and suggestions.

\section{Kazhdan--Lusztig $\tilde{R}$-polynomials}
\label{section two}

Herein, $(W,S)$ denotes an arbitrary Coxeter system. We refer the reader to \cite{bjorner2006combinatorics} and \cite{humphreys1992reflection} for the notations and definitions concerning Coxeter systems. In particular, for $u,v\in W$, we write $u\leq v$ to indicate that $u$ is smaller than or equal to $v$ in the Bruhat order of $(W,S)$. The Bruhat order can be characterised in terms of subwords as follows.

\begin{lem} \label{lema subword property}
Let $\underline{v}=s_1\ldots s_n$ be a reduced expression of some $v\in W$. Then, $u\leq v$ if and only if there exists a reduced expression $\underline{u}= s_{i_1}s_{i_2}\ldots s_{i_k}$ of $u$, where $1\leq i_1 <i_2<\ldots  < i_k\leq n$.	
\end{lem}

Let $\mathcal{A} = \mathbb{Z}[t,t^{-1}]$. The Hecke algebra $\mathcal{H}=\mathcal{H}(W,S)$ of $(W,S)$ is the associative and unital $\mathcal{A}$-algebra generated by $\{H_s \, |\, s\in S\}$ with relations 
$$\begin{array}{rl}
	H_s^2= & 1+(t^{-1}-t)H_s, \\
	\underbrace{H_sH_tH_s \cdots}_{m_{st}\mbox{-times}} = & \underbrace{H_tH_sH_t \cdots}_{m_{st}\mbox{-times}}
\end{array}
$$
where $m_{st}$ denotes the order of $st$ in $W$. 

Given a reduced expression $\underline{v}=s_{1} \cdots s_{k}$ of an element $v\in W$, we define $H_v=H_{s_1}\cdots H_{s_k}$. It is well known that $H_v$ is well defined, that is, it does not depend on the choice of a reduced expression. We also set $H_e:=1$, where $e$ denotes the identity of $W$. The set $\{H_v \, | \, v\in W\}$ is a basis of $\mathcal{H}$, which we call the standard basis of $\mathcal{H}$. It follows directly from the definition of $\mathcal{H}$ that each generator $H_s$ is invertible and we have $H_s^{-1}= H_s+(t-t^{-1})$. Hence, all elements in the standard basis of $\mathcal{H}$ are invertible as well. 

For $u,v\in W$, we denote the corresponding KL $R$-polynomials by $R_{u,v}(t)\in \mathcal{A}$. These (Laurent) polynomials are, by definition, the coefficients of the expansion of $(H_{v^{-1}})^{-1}$ in terms of the standard basis. Concretely, for $v\in W$ we have 
\begin{equation}
	(H_{v^{-1}})^{-1} = \sum_{u\in W} R_{u,v}(t) H_u
\end{equation}

\begin{rem}\rm
	Herein, we follow the normalisation of Soergel \cite{soergel1997kazhdan} rather than the typical normalisation given in \cite{kazhdan1979representations}, \cite{humphreys1992reflection} or \cite{bjorner2006combinatorics}. In particular, if $R'_{u,v}(t)$ denotes the $R$-polynomial in that sources, then it is related with our $R$-polynomial, $R_{u,v}(t)$, by the formula
	\begin{equation}
		R_{u,v}(t)=t^{l(u)-l(v)}R'_{u,v}(t^2). 
	\end{equation} 
\end{rem}

It is clear that $H_u$ appears in the expansion of $(H_{v^{-1}})^{-1} $ in terms of the standard basis if and only if $u\leq v$. In other words, we have $R_{u,v}(t)\neq 0$ if and only if $u  \leq v$. Meanwhile, one can verify that even in the simplest examples, the $R$-polynomials do not have non-negative coefficients. This drawback is overcome with the following:

\begin{propos}\cite[Propositions 5.3.1 and 5.3.2]{bjorner2006combinatorics} 
\label{propose decomposition R tilde polynomials}
	Let $u,v\in W$. Then, there exists a unique polynomial $\tilde{R}_{u,v}(t)\in \mathbb{N}[t]$ such that 
	\begin{equation}
		R_{u,v}(t)= \tilde{R}_{u,v}(t-t^{-1}).
	\end{equation}
Now, assume that $u\leq v$.  Then, $\tilde{R}_{u,v}(t)$ is a monic polynomial of degree $l(v)-l(u)$.	Furthermore, if $ s\in S$ satisfies $l(vs)<l(v)$, then 
	
	\begin{equation} \label{decomposition R tilde polynomials}
		\tilde{R}_{u,v}(t)= \left\{ \begin{array}{ll}
			\tilde{R}_{us,vs}(t),                       &   \mbox{ if }   l(us)<l(u); \\
			\tilde{R}_{us,vs}(t) + t\tilde{R}_{u,vs}(t), & \mbox{ if }   l(us)>l(u).
		\end{array} \right.
		\end{equation}
\end{propos}  

Equation (\ref{decomposition R tilde polynomials}) allows us to calculate the $\tilde{R}$-polynomials with the initial conditions $\tilde{R}_{v,v}(t)=1$ and $\tilde{R}_{u,v}(t)=0$ if $u\not \leq v$.

\section{Diagrammatic $\tilde{ \mathcal{R} }$-polynomials}  \label{section diagrams}

In this section, we associate to each (non-necessarily reduced) word $\underline{v}$ a binary tree $\mathbb{T}_{\underline{v}}$ with nodes decorated by certain diagrams. The diagrams represent morphisms in the category of Soergel bimodules. However, we only treat them as combinatorial objects herein. For the algebraic meaning of the diagrams, the reader is referred to \cite{elias2010diagrammatics,elias2016soergel,libedinsky2008categorie,libedinsky2015light,libedinsky2017gentle}. In the following, we treat $S$ as a set of colours.

\begin{defi} \rm
	An $S$-graph is a finite planar graph with its boundary embedded in a planar strip $\mathbb{R}\times [0,1]$, whose edges are coloured by elements of $S$, and all of whose vertices are of the following types\footnote{The reader familiar with diagrammatics for Soergel categories should note that the definition of $S$-graph herein does not allow trivalent vertices or polynomials to float in the graph. }:
	\begin{enumerate} 
		\item  Univalent vertices (Dots): $ \quad \univalent \quad $
        \item $2m_{st}$-valent vertices: $\quad \multivalent \quad$
	\end{enumerate}
	Here, exactly $2m_{st}$ edges originate from the vertex, and they are coloured alternately by $s$ and $t$. For instance, the picture above corresponds to $m_{\red{s} \blue{t}}=4$. 	 
\end{defi}

\begin{exa}\rm
	Let $S=\{{\red{s}},{\blue{t}},{\green{u}} \}$ with $m_{{\red{s}}{\blue{t}}}=3$,  $m_{{\red{s}}{\green{u}}}=2$ and $ m_{{\blue{t}}{\green{u}}}=3$. An example of an $S$-graph is given by  
	\begin{equation} \label{example Soergel graph}
		\exampleSgraph
	\end{equation}
\end{exa}

The boundary points of an $S$-graph on $\mathbb{R}\times \{0\}$ and on $\mathbb{R}\times \{1\}$ determine two sequences of coloured points, and hence two words in the alphabet $S$. We call these two sequences the bottom sequence and top sequence, respectively. For example, the bottom (resp. top) sequence of the $S$-graph in (\ref{example Soergel graph}) is $ststuutss$ (resp. $ust$). 
\begin{defi}\rm
	Given an $S$-graph $D$ we define its degree, $\deg (D)$, as the number of dots.
\end{defi}

For instance, the degree of the $S$-graph in (\ref{example Soergel graph}) is $4$. For the remainder of this section, we fix a word $\underline{v}=s_1s_2\ldots s_{n}$, $s_i\in S$. The construction of $\mathbb{T}_{\underline{v}}$ is by induction on the depth of the nodes and  is as follows. In depth one, we have the following tree: 
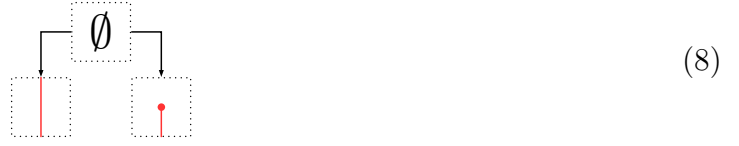
\begin{equation}
	\levelone
\end{equation}
Here, we assume that the first letter of $\underline{v}$, $s_1$, is associated with the colour red. Let $1<k\leq n$ and assume that a node $N$ of depth $k-1$ has been decorated by the diagram $D$. Further, suppose that the top boundary of $D$, say $\underline{u}$, is a reduced expression for some $u\in W$. If the letter $s_k$ is blue, then we have two possibilities.

\begin{enumerate}
	\item If $l(us_k)>l(u)$, then $N$ contains two child nodes that are decorated as follows. 
	\begin{equation}
		\treestepone
	\end{equation}
	\item If $l(us_k)<l(u)$, then $N$ contains one child node. It is well known that in this case, there exists a reduced expression $\underline{u}'$ of $u$ with a letter $s_k$ in its rightmost position. Furthermore, we can obtain $\underline{u}'$ from $\underline{u}$ by applying a sequence of braid movements. Diagrammatically, we move a blue edge to the rightmost position by placing a sequence of $2m_{st}$-valent vertices on the top of $D$. Subsequently, we connect the two blue edges as illustrated below.
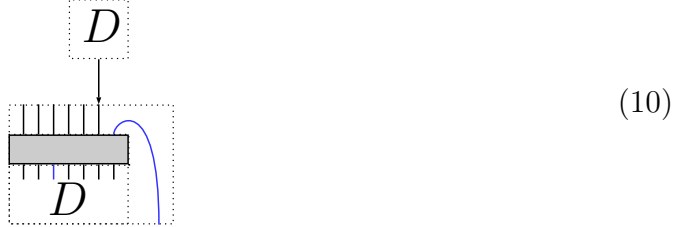
\begin{equation}
	\treesteptwo
\end{equation}
The grey region in the diagram above indicates the place where the sequence of $2m_{st}$-valent vertices is located. This completes the construction of $\mathbb{T}_{\underline{v}}$.
\end{enumerate}

The diagram-decorating leaves of $\mathbb{T}_{\underline{v}}$ are called light leaves\footnote{Herein, we use the term ``light leaves'' to indicate a subset of the typical set of light leaves, namely the ones that do not contain any trivalent vertex. See \cite{libedinsky2008categorie} for the original definition of light leaves and \cite{elias2016soergel} for its diagrammatic version. We hope that this  does not cause any confusion in the reader} of $\underline{v}$. The set of all light leaves of $\underline{v}$ is denoted by $\mathbb{L}_{\underline{v}}$. By construction, the bottom sequence of any element of $\mathbb{L}_{\underline{v}}$ is $\underline{v}$. Meanwhile, the top sequence of any element of $\mathbb{L}_{\underline{v}}$ is a reduced expression of some element of $W$, even if the word $\underline{v}$ is not reduced. Given $u\in W$, we denote by $\mathbb{L}_{\underline{v}}(u)$ the set of all light leaves of $\underline{v}$ with top sequence a reduced expression of $u$. 

\begin{rem} \rm  \label{remark ambiguity}
	Step $2$ in the construction of $\mathbb{T}_{\underline{v}}$ introduces some ambiguity in the sets $\mathbb{L}_{\underline{v}}$ and $\mathbb{L}_{\underline{v}}(u)$. The problem here is that we have multiple choices to pass from one reduced expression to another. Hence, Step 2 can be performed in many ways. In the following, we treat the sets $\mathbb{L}_{\underline{v}}$ and $\mathbb{L}_{\underline{v}}(u)$ as any admissible choice of the corresponding $S$-graphs. 
\end{rem}

\begin{exa}\rm  \label{exa simplest}
	Let $\underline{v}=\red{s} \red{s} \red{s}$. The tree $\mathbb{T}_{\underline{v}}$ associated to $\underline{v}$ is given by 
	\begin{equation}
		\simplest
	\end{equation} 
	Exactly five light leaves exist for $\underline{v}$. They are split into two classes according to their top sequence. We have
	\begin{equation}
		\mathbb{L}_{\underline{v}} = \mathbb{L}_{\underline{v}}(s) \cup \mathbb{L}_{\underline{v}}(e).
	\end{equation} 
\end{exa}

\begin{exa}\rm \label{example two different expressions}
	Let $W$ be the symmetric group on four letters with  $S=\{ {\red{s_1}}, {\blue{s_2}} , {\green{s_3}} \}$. Consider the words $\underline{v}= {\red{s_1}}{\blue{s_2}}{\red{s_1}}{\green{s_3}} {\blue{s_2}}{\red{s_1}}$ and $\underline{w}= {\blue{s_2}} {\green{s_3}}{\red{s_1}} {\blue{s_2}}{\green{s_3}}{\red{s_1}}$. Notice that both $\underline{v}$ and $\underline{w}$ are reduced words of the same element. We have
	
	\begin{equation}
	 	\mathbb{L}_{\underline{v}}(e)=    \left\{	\exampletreeA \right\}
	\end{equation} 
	\begin{equation}
	 	\mathbb{L}_{\underline{w}}(e)=    \left\{	\exampletreeB \right\}
	\end{equation} 
	
This example stresses that, although the sets $\mathbb{L}_{\underline{v}}(e)$ and $\mathbb{L}_{\underline{w}}(e)$	are different, they contain the same number of elements in each degree.
	\end{exa}

The number and degree of the elements of the set $\mathbb{L}_{\underline{v}}(u)$ are collected in a polynomial, $\tilde{\mathcal{R}}_{u,\underline{v}}(t) \in \mathbb{N}[t]$, which we call the ``diagrammatic $\tilde{R}$-polynomial'' associated to $u$ and $\underline{v}$. Concretely, we define
\begin{equation}
	\tilde{\mathcal{R}}_{u,\underline{v}}(t) = \sum_{l\in \mathbb{L}_{\underline{v}}(u) } t^{\deg (l)}. 
\end{equation}

Furthermore, if $\underline{v}$ is an empty word, we define 

\begin{equation}
	\tilde{\mathcal{R}}_{u,\underline{v}}(t) = \left\{  \begin{array}{ll}
	1, & \mbox{ if } u=e;\\
	0, & \mbox{ otherwise.}	
	\end{array} \right. . 
\end{equation}

\begin{rem}\rm
	As mentioned in Remark \ref{remark ambiguity}, the sets $\mathbb{L}_{\underline{v}}(u)$ are not uniquely defined. However, regardless of how Step $2$ is performed in the construction of the tree $\mathbb{T}_{\underline{v}}$, we always obtain the same polynomial $\tilde{\mathcal{R}}_{u,\underline{v}}(t)$. In other words, $\tilde{\mathcal{R}}_{u,\underline{v}}(t)$ is well defined, that is, it only depends on the word $\underline{v}$ and the element $u$.
\end{rem}

\begin{exa}\rm
	With the same notation as in Example \ref{example two different expressions}, we have
	\begin{equation}
		\tilde{\mathcal{R}}_{e,\underline{v}} (t)=	\tilde{\mathcal{R}}_{e,\underline{w}} (t)= t^6+3t^4+t^2.
	\end{equation}
\end{exa}

The following is the diagrammatic version of (\ref{decomposition R tilde polynomials}).

\begin{lem} \label{lema decomposition R tilde polynomials}
	Let $\underline{v}=s_1\ldots s_ns$ be a word and $u\in W $. Let $\underline{v}'$ be the word obtained from $\underline{v}$ by erasing the rightmost letter, that is, $\underline{v}'=s_1\ldots s_{n}$. Then, 
	 \begin{equation} \label{decomposition diagrammatic R tilde polynomials}
		\tilde{\mathcal{R}}_{u,\underline{v}}(t)= \left\{ \begin{array}{ll}
			\tilde{\mathcal{R}}_{us,\underline{v}'}(t),                       &   \mbox{ if }   l(us)<l(u); \\
			\tilde{\mathcal{R}}_{us,\underline{v}'}(t) + t\tilde{\mathcal{R}}_{u,\underline{v}'}(t), & \mbox{ if }   l(us)>l(u).
		\end{array} \right.
		\end{equation}
\end{lem}

\begin{demo}
Suppose we have fixed $\mathbb{T}_{\underline{v}'}$. We will construct $\mathbb{T}_{\underline{v}}$ from our particular choice of $\mathbb{T}_{\underline{v}'}$.\footnote{That is, among all possibilities for $\mathbb{T}_{\underline{v}'}$, we have fixed one and from this tree, we complete the construction of $\mathbb{T}_{\underline{v}}$ by applying Step 1 or Step 2 to each leaf node of $\mathbb{T}_{\underline{v}'}$.} In particular, we fixed the sets $\mathbb{L}_{\underline{v}'}(x)$, for all $x\in W$. Let $l\in \mathbb{L}_{\underline{v}'}(x)$. The $S$-graphs decorating the child nodes of the node decorated by $l$ in $\mathbb{T}_{\underline{v}}$ belong to $ \mathbb{L}_{\underline{v}}(x)$ or $ \mathbb{L}_{\underline{v}}(xs)$. Hence, all the elements in $ \mathbb{L}_{\underline{v}}(u)$ are obtained from the elements of $ \mathbb{L}_{\underline{v}'}(u)$ or $ \mathbb{L}_{\underline{v}'}(us)$. We split the proof into two cases.

\begin{itemize}
	\item[Case A. ] Suppose that $l(us)<l(u)$. Let $l\in \mathbb{L}_{\underline{v}'}(u)$. 
	Under this hypothesis, the node decorated by $l$ in $\mathbb{T}_{\underline{v}}$ contains one child node. Furthermore, the $S$-graph decorating this child node belongs to $\mathbb{L}_{\underline{v}}(us)$. Now, suppose that $l\in \mathbb{L}_{\underline{v}'}(us)$. In this case, the node decorated by $l$ in $\mathbb{T}_{\underline{v}}$ appears as
	
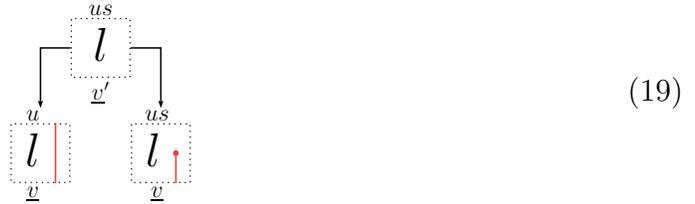
\begin{equation}
	\treeproof
\end{equation}

For the above, we conclude that all the elements of 	$ \mathbb{L}_{\underline{v}}(u)$ are obtained by adding a line in the rightmost region of each element of $ \mathbb{L}_{\underline{v}'}(us)$. Therefore, 
\begin{equation}
	\tilde{\mathcal{R}}_{u,\underline{v}}(t)= \tilde{\mathcal{R}}_{us,\underline{v}'}(t),
\end{equation}
proving the result in this case. 
\item[Case B. ] Suppose that $l(us)>l(u)$. Let $l\in \mathbb{L}_{\underline{v}'}(us)$. Under this hypothesis, the node decorated by $l$ in $\mathbb{T}_{\underline{v}}$ contains one child node. Furthermore, the $S$-graph decorating this child node contains the same degree as $l$ and belongs to $\mathbb{L}_{\underline{v}}(u)$. Therefore, each element in $\mathbb{L}_{\underline{v}'}(us)$ produces a light leaf in $\mathbb{L}_{\underline{v}}(u)$ of the same degree. 
   Now, suppose that $l\in \mathbb{L}_{\underline{v}'}(u)$. In this case, the node decorated by $l$ in $\mathbb{T}_{\underline{v}}$ appears as
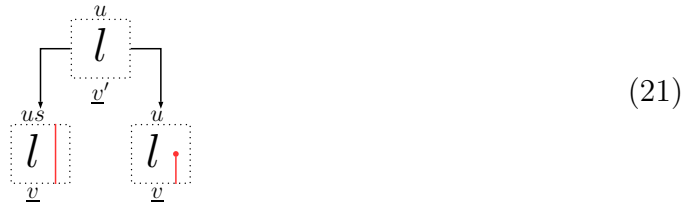
\begin{equation}
	\treeproofA
\end{equation}	

Therefore, each element $l\in \mathbb{L}_{\underline{v}'}(u)$ produces a light leaf $l'\in \mathbb{L}_{\underline{v}}(u)$ with $\deg(l')=\deg(l)+1$. Summing up, we obtain  
\begin{equation}
	\tilde{\mathcal{R}}_{u,\underline{v}}(t)= \tilde{\mathcal{R}}_{us,\underline{v}'}(t) + t\tilde{\mathcal{R}}_{u,\underline{v}'}(t).
\end{equation}

\end{itemize}
\end{demo}

\begin{teo} \label{teo diagrammatic R poly}
	Let $\underline{v}$ be a reduced expression of an element $v\in W$. Then, 
	\begin{equation}
		\tilde{\mathcal{R}}_{u,\underline{v}}(t)= \tilde{R}_{u,v}(t),
	\end{equation}
	for all $u\in W$. In particular, $\tilde{\mathcal{R}}_{u,\underline{v}}(t)$ does not depend on the particular choice of a reduced expression for $v$. 
\end{teo}

\begin{demo}
We proceed by an induction on $l(v)$. If $l(v)\leq 1$, the result is clear. Suppose that $l(v)>1$ and let $\underline{v}=s_1\ldots s_n s$ be a reduced expression of $v$. We remark that $\underline{v}'=s_1\ldots s_n$ is a reduced expression of $vs$ and that $l(vs)<l(v)$. Let $u\in W$ and suppose that $l(us)>l(u)$. Therefore, Proposition \ref{propose decomposition R tilde polynomials}, Lemma \ref{lema decomposition R tilde polynomials}, and our induction hypothesis yield
 
\begin{equation}
	\begin{array}{rl}
		\tilde{\mathcal{R}}_{u,\underline{v}}(t)= &  \tilde{\mathcal{R}}_{us,\underline{v}'}(t) + t\tilde{\mathcal{R}}_{u,\underline{v}'}(t)  \\
	= 	&  \tilde{R}_{us,vs}(t) + t\tilde{R}_{u,vs}(t)  \\
	= & \tilde{R}_{u,v}(t)
	\end{array}
\end{equation}
The case when $l(us)<l(u)$ is treated similarly.
\end{demo}

\section{ Closed formulas for $\tilde{R}$-polynomials.} \label{section closed}

In this section, we provide some examples of the computations of $\tilde{R}$-polynomials using the diagrammatic approach. In particular, we recover and generalise uniformly some closed formulas obtained by Pagliacci \cite{pagliacci2001explicit} and Marietti \cite{marietti2002closed}. 
  
\subsection{$\tilde{R}$-polynomials for permutations smaller than a transposition.} \label{section transposition}
 In this section, we restrict our attention to the symmetric group $\mathfrak{S}_n$. We express the permutations using either the cycle notation or one-line notation. For instance, $ (1,2,3)  $ and $231$ represent the same permutation. We denote by $s_i$ the simple transposition $(i,i+1)$, for $1\leq i <n$. We establish the convention that $s_i$ acts on the left on $\{ 1,2,\ldots ,n \}$. 
  
  The following was conjectured in \cite[Conjecture 7.7]{brenti1998kazhdan} and proven in \cite{marietti2002closed}. 
  
  \begin{teo} \label{Brenti conjecture}
	Let $u,v \in \mathfrak{S}_n$ be such that $u\leq v \leq (a,b)$ for some $1 \leq a< b \leq n  $. Therefore, 
	\begin{equation}
		\tilde{R}_{u,v}(t)=t^{c}(t^2+1)^{\frac{1}{2} (l(v)-l(u)-c )}
	\end{equation}
	for some $c\in \mathbb{N}$. 
\end{teo}
  
We prove a slight generalisation of Theorem \ref{Brenti conjecture} using the diagrammatic approach. We first require the following definition.

\begin{defi} \rm \label{definition UD-word}
	We say a word $\underline{v}$ is ``up and down'' (or UD-word) if 
\begin{equation} \label{almost favourite reduced expression}
\underline{v}= (s_{i_1}s_{i_2}\ldots s_{i_{r}}) s_t (s_{j_1}s_{j_2}\ldots s_{j_q}), 
\end{equation}
where 
\begin{enumerate}
	\item $1\leq i_1 <i_2 <\ldots <i_{r} <t <n$;
	\item  $1\leq j_q <\ldots <j_2 <j_1<t< n$.
\end{enumerate}

\end{defi}
  
The definition above does not exclude the possibility that one or both of the parentheses can be empty. 
It is noteworthy that in general, UD-words are not reduced expressions. 

\begin{lem} \label{lemma simple tree A}
	Let $\underline{v}$ be a UD-word. Then, the tree $\mathbb{T}_{\underline{v}}$ can be constructed 
	without using $6$-valent vertices. Furthermore, for this construction, if $u\in \mathfrak{S}_n$ and $\mathbb{L}_{\underline{v}}(u)\neq \emptyset $, then the top sequence of any element of $\mathbb{L}_{\underline{v}}(u)$ is a UD-word as well.
\end{lem}  

\begin{demo}
	We must use a $6$-valent vertex in the construction of  $\mathbb{T}_{\underline{v}}$ only if  $\underline{v}$ admits a subword of the type $s_is_{i+1}s_is_{i+1}$ or $s_{i}s_{i-1}s_is_{i-1}$. Both options are impossible for a UD-word. 
	
	Meanwhile, if $\mathbb{T}_{\underline{v}}$ is constructed without using $6$-valent vertices, then the top sequence of any element $l\in \mathbb{L}_{\underline{v}}(u)$ is obtained from $\underline{v}$ by erasing some letters. Therefore, such a sequence must be a UD-word.
\end{demo}

\begin{exa}\label{example we are forced} \rm
	Lemma \ref{lemma simple tree A} tells us that among all the possibilities to construct 
	$\mathbb{T}_{\underline{v}}$ for a UD-word, there exists at least one option in which we do not need to use $6$-valent vertices. Consider the UD-word $ \underline{v} = {\red{s_1} }{\blue{s_3} }{\green{s_4} }{\blue{s_3} }{\red{s_1}} $. The reader can easily verify that the following $S$-graph decorates a node of depth $4$ in $\mathbb{T}_{\underline{v}}$.
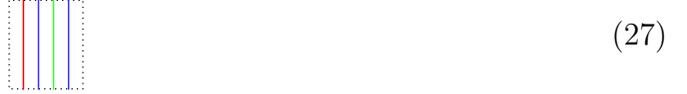
	\begin{equation}
		\exatreeUD
	\end{equation}
	This node contains only one child node that can be decorated by 
	\begin{equation} \label{two decorations}
		\exatreeUDA  \qquad \mbox{ or } \qquad \exatreeUDB\quad .
	\end{equation}
	Lemma \ref{lemma simple tree A} only claims the existence of the $S$-graph on the left of (\ref{two decorations}). The $S$-graph on the right of (\ref{two decorations}) is, despite being slightly bizarre, an admissible decoration. 
	\end{exa}

\begin{teo} \label{teo R poly UD-words}
	Let $\underline{v}$ be a $UD$-word. If $u\in \mathfrak{S}_n$ satisfies $\mathbb{L}_{\underline{v}}(u)\neq \emptyset$, then 
	
	\begin{equation} \label{equation teo diagram R dominated by a transposition}
		\tilde{\mathcal{R}}_{u,\underline{v}}(t) =t^{c}(t^2+1)^{\frac{1}{2} (l(\underline{v})-l(u)-c )}
	\end{equation}
	for some $c\in \mathbb{N}$. 	
	\end{teo}  

\begin{demo}
	By Lemma \ref{lemma simple tree A}, there is one admissible construction for $\mathbb{T}_{\underline{v}}$ without using $6$-valent vertices. This allows us to obtain the set $\mathbb{L}_{\underline{v}}(u)$ without constructing the whole tree $\mathbb{T}_{\underline{v}}$. This is performed by analysing each letter occurring in $\underline{v}$ separately. Eight cases are to be considered. In the diagrams below, we only draw the letters involved in each case. 
	\begin{itemize}
		\item[Case A1.] Exactly one $s_{i}$ appears in $\underline{v}$ and $s_i$ does not appear in $u$. In this case, the elements in $\mathbb{L}_{\underline{v}}(u)$ are forced to be of the form 
		\begin{equation}
			\CaseAone 
		\end{equation}
		\item[Case A2.] Exactly one $s_{i}$ appears in $\underline{v}$ and $s_i$ appears (necessarily once) in $u$. In this case, elements in $\mathbb{L}_{\underline{v}}(u)$ are forced to be of the form 
		\begin{equation}
			\CaseAtwo 
		\end{equation}
		\item[Case B1.] The letter $s_i$ appears twice in $\underline{v}$ and $u$. It is noteworthy that this forces a letter $s_{i+1}$ to appear in between the two $s_i$. In this case, the elements in $\mathbb{L}_{\underline{v}}(u)$ are forced to be of the form 
		\begin{equation}
			\CaseBone 
		\end{equation}
		\item[Case C1.] The letter $s_i$ appears twice in $\underline{v}$, $s_{i+1}$ appears in $u $, and $s_i $ appears once in $u$ on the left of $s_{i+1}$. In this case, the elements in $\mathbb{L}_{\underline{v}}(u)$ are forced to be of the form 
		\begin{equation}
			\CaseCone 
		\end{equation}
        \item[Case C2.] The letter $s_i$ appears twice in $\underline{v}$, $s_{i+1}$ appears in $u $, and $s_i $ appears once in $u$ on the right of $s_{i+1}$. In this case, the elements in $\mathbb{L}_{\underline{v}}(u)$ are forced to be of the form 
		\begin{equation}
			\CaseCtwo 
		\end{equation}
		\item[Case C3.] The letter $s_i$ appears twice in $\underline{v}$, and $s_{i+1}$ does not appear in $u $. In this case, the elements in $\mathbb{L}_{\underline{v}}(u)$ are forced to be of the form 
		\begin{equation}
			\CaseCthree 
		\end{equation}
		\item[Case D1.] The letter $s_i$ appears twice in $\underline{v}$, $s_{i+1}$ appears in $u $, and $s_i$ does not appear in $u$. In this case, the elements in $\mathbb{L}_{\underline{v}}(u)$ are forced to be of the form 
		\begin{equation}
			\CaseDone 
		\end{equation}
		\item[Case D2.] The letter $s_i$ appears twice in $\underline{v}$, $s_{i+1}$ does not appear in $u $, and $s_i$ does not appear in $u$. In this case, the elements in $\mathbb{L}_{\underline{v}}(u)$ have two options 
		\begin{equation} \label{last case}
			\CaseDtwo 
		\end{equation}
	\end{itemize} 
We construct the elements of $\mathbb{L}_{\underline{v}}(u)$ by considering the letters appearing in $\underline{v}$ in the decreasing order. The result shows that the process is independent at each stage. This explains the reason that the formula in (\ref{equation teo diagram R dominated by a transposition}) is a product.

The only case when more than one option exists is D2. In case D2, the factor $(t^{2}+1)$ (see the degrees of the diagrams in (\ref{last case})) appears in $\tilde{\mathcal{R}}_{u,\underline{v}}(t) $. The other cases contribute with a power of $t$. The exponent in this power can be $0$, $1$, or $2$ according to the number of dots involved in each case. The sum of these exponents yields the value of the integer $c$.
\end{demo}

\begin{coro}
	Theorem \ref{Brenti conjecture} holds.
\end{coro}  

\begin{demo}
A reduced expression for $(a,b)$ is given by 
\begin{equation}
s_as_{a+1}\ldots s_{b-2}s_{b-1}s_{b-2}\ldots s_{a+1}s_{a}
\end{equation}

It follows directly from Lemma \ref{lema subword property} that any $e< v\leq (a,b)$ has a reduced expression which is a UD-word. The result is now a direct consequence of Theorem \ref{teo R poly UD-words}.
\end{demo}

\begin{exa}\rm
Consider the UD-word  $\underline{v}=s_1s_2s_4s_5s_7s_9s_8s_7s_4s_3s_2s_1$ and the element $u=s_7s_9s_8s_3\in \mathfrak{S}_{10}$. We enlist each letter with its corresponding case as in the proof of Theorem \ref{teo R poly UD-words}. 

\begin{equation}
	\begin{tabular}{|c|c|c|c|c|c|c|c|c|}
\hline
Letter & $s_9$ & $s_8$& $s_7$& $s_5$& $s_4$& $s_3$& $s_2$& $s_1$\\
\hline
Case & A2 & A2& C1& A1  & D2& A2 & D1 &D2 \\
\hline
Exponent & $0$ &$0$ &$1$ & $1$  & & $0$& $2$&  \\
\hline
\end{tabular}
\end{equation}
As Case D2 occurs twice, the factor $(t^2+1)^2$ appears. Meanwhile, by adding the exponents in the table above, we obtain $c=4$. Therefore, 
\begin{equation}
	\tilde{\mathcal{R}}_{u,\underline{v}}(t) =t^{4}(t^2+1)^2.
\end{equation}
The four light leaves in $\mathbb{L}_{\underline{v}}(u)$ are depicted in (\ref{total leaves}), where we have drawn all edges with the same colour and have replaced the letters $s_i$ for $i$.  
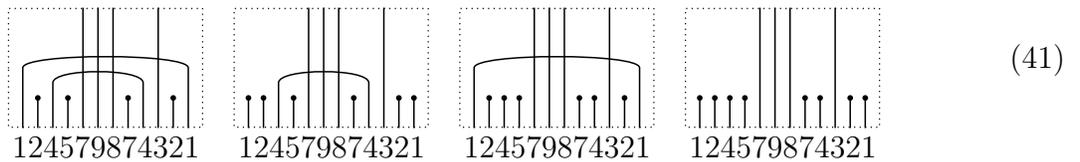
\begin{equation} \label{total leaves}
	\totalleaves
\end{equation}
\end{exa}


\vspace{.3cm}
\subsection{$\tilde{\mathcal{R}}$-polynomials for the simplest words.}

In this section, we again consider an arbitrary Coxeter system $(W,S)$. Among all the words in the alphabet $S$, the simplest ones are those formed by one letter. We compute the diagrammatic $\tilde{\mathcal{R}}$-polynomials for such words. For the remainder of this section, we set $s\in S$ and for $n\in \mathbb{N}$, we define
\begin{equation}
	\underline{n_s}= sss\cdots \qquad (n\mbox{-times}).
\end{equation}
  
A moment's though reveals that $\tilde{\mathcal{R}}_{u,\underline{n_s}}(t)=0$,    unless $u=e$ or $u=s$. To treat the cases $u=e$ and $u=s$, we introduce the Fibonacci polynomials.
  
 \begin{defi} \rm
 	For each $n\in \mathbb{N}$, we define the $n$-th Fibonacci polynomial $F_n(v)\in \mathbb{N}$ by the recurrence 
\begin{equation}
	F_n(v)=vF_{n-1}(v)+F_{n-2}(v), 
\end{equation}
with initial conditions $F_0(v) = 1$ and $F_1(v)= v$.
 \end{defi} 
  
The Fibonacci polynomials have a combinatorial description in terms of paths in the Fibonacci tree.    
  
\begin{defi}\rm
	The $n$-th Fibonacci tree $FT_n$ is the binary tree defined inductively as follows. First, we have
	\begin{equation}
		FT_1 \longrightarrow \fibotreeone 
	\end{equation}
	Suppose now that $FT_n$ has been defined. Then, $FT_{n+1}$ is obtained from $FT_n$ by adding one child in each leaf node of $FT_n$ that is a left brother in $FT_n$, and two child nodes (one to the left and one to right) for the other leaf nodes of $FT_n$. 
\end{defi}  
  
\begin{exa} \rm
$FT_6$ is depicted below.
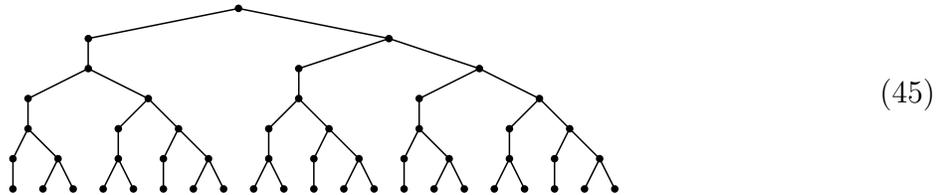
\begin{equation}
	\Fibo
\end{equation}	
\end{exa}

Let $P_n$ be the set of all paths\footnote{In this paper, a path is always going down.} from the root to the leaves of $FT_n$. To each $p\in P_n$, we define $\rho (p)$ as the number of right steps in $p$. Finally, we define $P_n'$ to be the subset of $P_n$ formed by the paths $p$ such that the last step is not a left step. It is a straightforward exercise to verify. 

\begin{equation} \label{interpretation Fibo poly}
	F_n(t)= \sum_{p\in P'_n} t^{\rho (p)}. 
\end{equation}

 \begin{teo}
 	Let $n\in \mathbb{N}$. Then, 
 	\begin{equation}
 		\tilde{\mathcal{R}}_{u,\underline{n_s}}(t)= \left\{ \begin{array}{rl}
 			F_n(t), & \mbox{ if } u=e;\\
 			F_{n-1}(t), & \mbox{ if } u=s. 
 		\end{array}  \right.
 	\end{equation}
 \end{teo}

\begin{demo}
	By disregarding the decorations in the nodes of  $\mathbb{T}_{\underline{n_{s}}}$ temporarily, we observed that $\mathbb{T}_{\underline{n_{s}}}$ coincides with $FT_n$. In this setting, we observed the light leaves of $\underline{n_s}$ as paths in $\mathbb{T}_{\underline{n_{s}}}$. With this interpretation, the right steps correspond to adding a dot in the rightmost region of a diagram; the left steps correspond to adding a straight line in the rightmost region of a diagram and central steps\footnote{That is, the ones that are not left or right.} correspond to transforming the rightmost straight line of a diagram into a loop. 
	
	As the degree of a diagram is the number of dots, we obtain that the degree of a light leaf of $\underline{n_s}$ is the number of right steps. Meanwhile, we have
	\begin{equation}
		\mathbb{L}_{\underline{n_s}}(e) = \{ l\in \mathbb{L}_{\underline{n_s}} \,|\, \mbox{last step in } l \mbox{ is not a left step } \}
	\end{equation}
	Therefore, (\ref{interpretation Fibo poly}) yields
	\begin{equation}
	\tilde{\mathcal{R}}_{e,\underline{n_s}}(t)=	\sum_{l\in \mathbb{L}_{\underline{n_s}}(e)} t^{\deg(l)} = \sum_{p\in P_n'} t^{\rho (p)} = F_n(t). 
	\end{equation}
	Finally, we found that a degree-preserving bijection exists: $ \mathbb{L}_{\underline{(n-1)_s}} (e) \longrightarrow \mathbb{L}_{\underline{n_s}} (s) $, and is produced by adding a straight line on the rightmost region of each element of $\mathbb{L}_{\underline{(n-1)_s}} (e)$. Hence, 
	\begin{equation}
		\tilde{\mathcal{R}}_{s,\underline{n_s}}(t)= \sum_{l\in \mathbb{L}_{\underline{n_s}}(s)} t^{\deg(l)}=  \sum_{l\in \mathbb{L}_{\underline{(n-1)_s}}(e)}t^{\deg(l)}= \tilde{\mathcal{R}}_{e,\underline{(n-1)_s}}(t)=F_{n-1}(t).
	\end{equation}
	   
\end{demo}

\subsection{The polynomial $\tilde{R}_{e,v}(t)$ for $ v=34\cdots n12$.}   
   In this section, we focus on the symmetric group $\mathfrak{S}_n $. We obtain a closed formula of the polynomial $\tilde{R}_{e,v}(t)$ for $ v=34 \cdots n12$. As in the previous section, this formula is given in terms of the Fibonacci polynomials. Such a formula was previously obtained by Pagliacci \cite[Theorem 4.1]{pagliacci2001explicit} using an inductive argument. We present a combinatorial proof in the stricter sense. That is, we  construct a (degree-preserving) bijection between  $\mathbb{L}_{\underline{v_n}}(e)$ and paths in the Fibonacci tree $FT_{n-3}$.  
      
Hereinafter, we set $n\in \mathbb{N}$ and $v_n=34\cdots n12$. It can be easily shown that 
\begin{equation}
	\underline{v_n} = s_2s_1s_3s_2 s_4s_3 \cdots s_{n-1}s_{n-2}
\end{equation}
   is a reduced expression for $v_n$. 
  
We recall from section 4.1 that the key to obtaining a closed formula for $\tilde{R}$-polynomials was the fact that the relevant tree was constructed without using $6$-valent vertices. However, this does not apply for $ \underline{v_n}$ due to the occurrence of subwords of the form $ s_{i}s_{i+1}s_is_{i+1}$ in $\underline{v_n}$; thus, we must use $6$-valent vertices in the construction of $\mathbb{T}_{\underline{v_n}}$. Nevertheless, some $u\leq v_n$ exist such that the sets $\mathbb{L}_{\underline{v_n}}(u) $ can be constructed without using the $6$-valent vertices. In particular, we have the following:
\begin{lem} \label{lemma pagliacci without sixvalent}
	The set $\mathbb{L}_{\underline{v_n}}(e)$ can be constructed without using $6$-valent vertices. 
\end{lem} 
\begin{demo}
Suppose we are forced to use a $6$-valent vertex at some stage of the construction of $\mathbb{T}_{\underline{v_n}}$. At this level, we must have the following diagram:
	\begin{equation} \label{six valent}
		\sixvalent
	\end{equation}

\medskip	
We recall that elements in $\mathbb{L}_{\underline{v_n}}(e)$ are precisely the light leaves with empty top sequence. Because letters $ s_{i}$ and $s_{i+1} $ no longer exist in $\underline{v_n}$ on the right of the ones shown in (\ref{six valent}), the top sequence of the light leaves emanating from this diagram contain at least two letters. In particular, they do not belong to $\mathbb{L}_{\underline{v_n}}(e)$.
\end{demo}   
   
To provide a closed formula for polynomials $\tilde{R}_{e,v_n}(t)$, we introduce slightly modified Fibonacci polynomials. For $n\geq 1$, we define
\begin{equation}
	\mathcal{F}_{n}(t)=t^{n-1}F_{n+1}(t). 
\end{equation}  

These polynomials contain the following combinatorial interpretation. Given a path $p \in P_n$, we define $\lambda (p)$ to be the number of left steps in $p$. Subsequently, we have

\begin{equation}
	\mathcal{F}_{n}(t)= \sum_{p\in P_n} t^{2(n-\lambda (p))} .
\end{equation}

\begin{teo} \label{teo pagliacci}
	For $n\geq 3$, we have
	\begin{equation} \label{equation Pagliacci}
	\tilde{R}_{e,v_n}(t)=	\tilde{\mathcal{R}}_{e,\underline{v_n}} (t) = t^2\mathcal{F}_{n-3}(t)=t^{n-2}F_{n-2}(t).
	\end{equation}
\end{teo} 
 
\begin{demo}
	We only need to show the equation in the middle of (\ref{equation Pagliacci}). The idea is to construct a bijection between the elements in $\mathbb{L}_{\underline{v_n}}(e)$ and the paths in the Fibonacci tree $FT_{n-3}$. 
	
	By Lemma \ref{lemma pagliacci without sixvalent}, we can construct the set $\mathbb{L}_{\underline{v_n}}(e)$ without using $6$-valent vertices. This allows us to treat each letter occurring in $\underline{v_n}$ separately. We recall that elements of $\mathbb{L}_{\underline{v_n}}(e)$ are the light leaves of $\underline{v_n}$ with an empty top sequence. 
	
	Let us begin by considering the letters $s_1$ and $s_{n-1}$. These letters appear only once in $ \underline{v_{n}}$; to obtain an element in $\mathbb{L}_{\underline{v_n}}(e)$, we must use dots over such letters. Now, consider the letter $s_{2}$. It appears twice in $\underline{v_{n}}$. Because we do not wish for them to occur in the top sequence, we would have to eliminate them. This can be accomplished by two options:
	\begin{itemize}
		\item[(R)] Placing a dot over each one of these letters; or
		\item[(L)] Connecting these two letters to form a loop.
	\end{itemize}
	 Option (R) increases the degree by two, and option (L) maintains the same degree. Now, consider the letter $s_3$. We remark that one of the two occurrences of the letter $s_3$ in $\underline{v_n}$ is in between the two letters $s_{2}$. If we use the first option (R) for thw letter $s_2$, then we have again options (R) and (L) for the letter $s_3$. On the contrary, if we use option (L) for the letter $s_{2}$, the letter $s_3$ in between the two $s_2$ is ``trapped'' by the loop formed by the letters $s_2$. Therefore, we must place a dot over each letter $s_{3}$. In this case, we increase the degree by two. We denote this unique option by (C). We continue with letter $s_4$. If we use option (R) or (C) for the letter $s_3$, then we can use options (R) and (L) for $s_4$. If we use option (L) for the letter $s_3$, we are forced to choose option (C) for $s_4$. The pattern is repeated continuously until we reach the letter $s_{n-2}$. 
	 
	 In summary, each element of $\mathbb{L}_{\underline{v_n}}(e)$ is determined by a word in the alphabet $\{C,L,R\}$ of length $n-3$ satisfying:
	 \begin{enumerate}
	 	\item The first letter is $R $ or $L$.
	 	\item Each letter $L$ is followed by a letter $C$.
	 	\item Each letter $C$ is preceded by a letter $L$.
	 \end{enumerate}
	 
Now, the promised bijection $\mathbb{L}_{\underline{v_n}}(e) \longrightarrow P_{n-3}$ is given by replacing letters $C$, $L$, and $R$ by the left steps, right steps, and central steps in $FT_{n-3}$, respectively. Given $l\in \mathbb{L}_{\underline{v_n}}(e)$, we denote by $p_l$ the image of $l$ under the map above. By construction, we have 
\begin{equation}
	\deg(l) =2 + 2((n-3)- \lambda (p_l)).
\end{equation}
Therefore, we obtain
\begin{equation}
	\tilde{\mathcal{R}}_{e,\underline{v_n}} (t) =  \sum_{l\in \mathbb{L}_{\underline{v_n}}(e)} t^{\deg (l)} = \sum_{p\in P_n} t^2 t^{2((n-3)- \lambda (p_l))}  = t^2  \mathcal{F}_{n-3}.
\end{equation}
\end{demo} 
   
\begin{exa}\rm In (\ref{tree Pagliacci}), we illustrate the bijection in the proof of Theorem \ref{teo pagliacci} for $n=7$. Because of space limitations, we replace letters $s_i$ by $i$. In this case, we have
$\tilde{\mathcal{R}}_{e,\underline{v_7}} (t) =t^{10}+4t^{8}+3t^{6} $.
 \begin{equation} \label{tree Pagliacci}
 	\treePagliacci
\end{equation} 	
\end{exa}

   \vspace{.3cm}
\subsection{$\tilde{R}_{e,v}(t)$ for $v$ a $321$-avoiding and $2$-repeating permutation.}   

In this section, we generalise the formula obtained in the previous section.

\begin{defi}\rm
	A permutation $v\in \mathfrak{S}_n$ is called $321$-avoiding if no reduced expression $\underline{v}$ of $v$ contains a subword of consecutive letters of the form $s_is_{i \pm 1}s_i$.   
\end{defi}

Let $v\in \mathfrak{S}_n$ be a $321$-avoiding permutation. In this case, the number of occurrences of any letter $s_i$ in $v$ is well defined. That is, it does not depend on the particular choice of a reduced expression for $v$. We define $n_v(i)$ to be the number of occurrences of the letter $s_i$ in $v$, for $1\leq i <n$. Subsequently, we say that $v$ is $2$-repeating if $n_v(i)\leq 2$, for all $1\leq i <n$. It is clear that permutations of the form $34\ldots n12$ considered in the previous section are $321$-avoiding and $2$-repeating. We now introduce a diagrammatic method of representing such permutations. Consider the set 

\begin{equation}
	\mathcal{T}= \{ (i,j)\in \mathbb{Z}^2 \, | \, j\leq 0, \, |i|\leq |j| \, ,\, i\equiv j \mbox{ mod } 2  \}. 
\end{equation}  
We call a finite subset of $\mathcal{T}$ a \emph{configuration of points}.

\begin{defi} \label{defin configuration}
We say a configuration of points $\mathcal{C}$ is admissible if it satisfies the following:
\begin{enumerate}
	\item No three points in $\mathcal{C}$ with the same $y$-coordinate exist.
	\item Let $P_1=(i_1,j_1)$ and $P_2=(i_2,j_2)$ be points in $\mathcal{C}$ with $j_1=j_2$. Then, we have $|i_1-i_2|=2$. 
	\item Let $P_1=(i_1,j_1)$ and $P_2=(i_2,j_2)$ be points in $\mathcal{C}$ with $j_1=j_2$ and $i_1<i_2$. Then, $(i_1+1,j_1+1)$ and $(i_1+1,j_1-1)$ belong to $\mathcal{C}$. 
\end{enumerate}	
\end{defi}

An example of an admissible configuration of points is illustrated in (\ref{configuration}). 
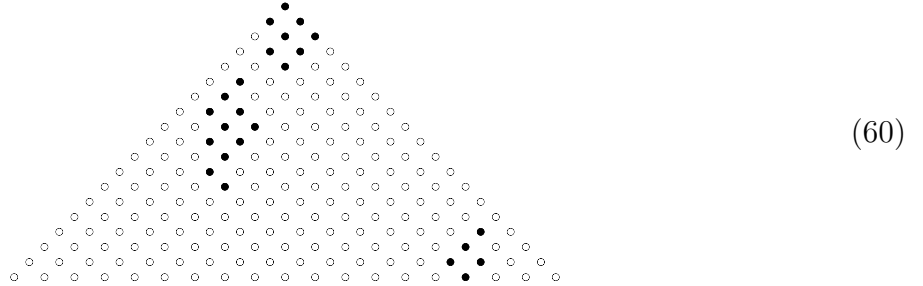
\begin{equation} \label{configuration}
	\configuration
\end{equation}	
We now associate to each admissible configuration of points, $\mathcal{C}$, a word in the alphabet $S=\{s_i\}_{i=1}^{\infty} $. First, each point $P=(i,j)\in \mathcal{C}$ is associated with the letter $s_{1-j}$. Then, the word $\underline{v_{\mathcal{C}}}$ associated to $\mathcal{C}$ is obtained by reading the points in $\mathcal{C}$ from left to right and from top to bottom. For instance, the word associated to the configuration of points in (\ref{configuration}) is given by
\begin{equation}
	s_{8}s_{10}s_{12}s_{7}s_{9}s_{11}s_{13}s_{6}s_{8}s_{10}s_{9}s_{2}s_{4}s_{1}s_{3}s_{5}s_{2}s_{4}s_{3}s_{18}s_{17}s_{19}s_{16}s_{18}.
\end{equation}

\begin{lem}
	Let $\mathcal{C}$ be an admissible configuration of points and $\underline{v_{\mathcal{C}}}$ its corresponding word. Then, 
	\begin{enumerate}
		\item The word $\underline{v_{\mathcal{C}}}$ is a reduced expression of some permutation $v_{\mathcal{C}}$. 
		\item The permutation $v_\mathcal{C}$ is $321$-avoiding.
		\item The permutation $v_\mathcal{C}$ is $2$-repeating.
	\end{enumerate}
\end{lem}

\begin{demo}
The first two claims are a consequence of \cite[Lemma 1]{billey2001kazhdan}. The last claim follows immediately by the first condition in Definition \ref{defin configuration}.
 \end{demo}

\begin{lem} \label{lema permutation to configuration}
Let $v$ be a $321$-avoiding and $2$-repeating permutation. Then, there exists an admissible configuration of points, $\mathcal{C}$, such that $\underline{v_{\mathcal{C}}}$ is a reduced expression of $v$.  
\end{lem}

\begin{demo}
After a rotation of $90$ degrees, the ``heap'' of $v$, constructed in \cite[Section 3]{billey2001kazhdan}, is a configuration satisfying the requirements of the lemma. 
\end{demo}

\begin{exa} \label{exa configuration}\rm
	As mentioned earlier, permutations $v_n=34\ldots n12$ considered in the previous section are $321$-avoiding and $2$-repeating permutations. The word associated to the configuration in (\ref{confi vn}) is a reduced word for $v_9$.
	\begin{equation}  \label{confi vn}
		\configurationA
	\end{equation} 
\end{exa}

Contrary to what Example \ref{exa configuration} might suggest, the configuration in Lemma \ref{lema permutation to configuration} is not unique. Considering the five points located in the bottom right of (\ref{configuration}), we can move all these points to two positions to the left or to the right without altering the element represented by the configuration. 

Let $v\in \mathfrak{S}_n$ be a $321$-avoiding and $2$-repeating permutation. We now define some statistics on $v$. First, we define $n_v^1$ to be the number of letters in $v$ that appear exactly once. Meanwhile, we say that a set of consecutive integers $\{a,a+1,\ldots , b\}$ is a $2$-chain of $v$ if $n_v(i)=2$, for all $i\in \{a,a+1,\ldots , b\}$, and $n_v(a-1)\neq 2$ and $n_v(b+1)\neq 2$. For instance, if $v$ is the permutation obtained from the configuration in (\ref{configuration}), then the $2$-chains of $v$ are: $\{2,3,4 \}$, $\{ 8,9,10 \}$, and $\{ 18 \}$. We define $\kappa (v) $ to be the number of $2$-chains in $v$. 
Let $c_1, c_2, \ldots, c_{\kappa(v)}$ be all the $2$-chains in $v$, indexed in a manner that the numbers in $c_i$ are less than the numbers in $c_j$, for $i<j$. Finally, we define $\lambda_{i}$ as the cardinality of $c_i$. 
With these definitions, we can enunciate a vast generalisation of Pagliacci's formula.

\begin{teo} \label{teo generalisation Pagliacci}
	Let $v\in \mathfrak{S}_n$ be a $321$-avoiding and $2$-repeating permutation. Then, 
	 \begin{equation} \label{equation final}
	 	\tilde{R}_{e,v}(t) =t^{n_v^1}\prod_{i=1}^{\kappa (v)} \mathcal{F}_{\lambda_i} (t).
	 \end{equation}	 
\end{teo}

\begin{demo}
Let $\mathcal{C}$ be any admissible configuration of points such that its associated word, $\underline{v_{\mathcal{C}}}$, is a reduced word of $v$. We must calculate the diagrammatic $\tilde{R}$-polynomial $\tilde{\mathcal{R}}_{e,\underline{v_{\mathcal{C}}}}(t)$. Hence, we have to determine the set $\mathbb{L}_{\underline{v_{\mathcal{C}}}}(e)$. The same argument used in the proof of Lemma \ref{lemma pagliacci without sixvalent} reveals that $\mathbb{L}_{\underline{v_{\mathcal{C}}}}(e)$ can be constructed without using $6$-valent vertices; therefore, we can consider each letter occurring in $ \underline{v_{\mathcal{C}}}$ separately.

	Let $N_v^1$ be the set of all letters occurring exactly once in $v$, such that $n_v^1=|N_v^1|$. Over each element of $N_v^1$, we must locate a dot. This explains the factor $t^{n_v^1}$ in (\ref{equation final}). Subsequently, we can think as if the letters in $N_v^1$ never occurred in $\underline{v_{\mathcal{C}}}$, or equivalently in $\mathcal{C}$. Let $\mathcal{C}'$ be the configuration of points obtained from $\mathcal{C}$ by eliminating the points associated to the letters in $N_v^1$. The configuration $\mathcal{C}'$ splits naturally into $ \kappa(v)$ subconfigurations of the form
	
\begin{equation}
	\configurationB
\end{equation}
	
	Let us index the subconfigurations occurring in $\mathcal{C}'$ from top to bottom as follows: $\mathcal{C}_1', \mathcal{C}_2', \ldots ,\mathcal{C}_{\kappa(v)}'$. For each $1\leq i \leq \kappa (v)$, we define $\underline{v_{\mathcal{C}_i'}}$ to be the word obtained by reading the points in $ \mathcal{C}_i'$ from left to right and from top to bottom. We remark that the configurations $\{ \mathcal{C}_i' \}$ are not admissible. However, the word $\underline{v_{\mathcal{C}_i'}}$ is still valid. By applying a similar argument to the one used in the proof of Theorem \ref{teo pagliacci}, we obtain 
	\begin{equation}
		\tilde{\mathcal{R}}_{e,\underline{v_{\mathcal{C}_i'}}}(t)=\mathcal{F}_{\lambda_i}(t),
	\end{equation}
	as $\lambda_i$ is the number of letters involved in $\underline{v_{\mathcal{C}_i'}} $. 
	
	The crucial point is that $\mathcal{C}_i'$ and $\mathcal{C}_j'$ are independent, for $i\neq j$. That is, any generator involved in $\mathcal{C}_i'$ commutes with any generator involved in $\mathcal{C}_j'$, so that we can treat each of these subconfigurations separately. We finally obtain
	\begin{equation}
	\tilde{R}_{e,v}(t) = \tilde{\mathcal{R}}_{e,\underline{v_{\mathcal{C}}}}(t)
	= t^{n_v^1} \prod_{i=1}^{\kappa (v)} \tilde{\mathcal{R}}_{e,\underline{v_{\mathcal{C}_i'}}}(t) 
	=t^{n_v^1}\prod_{i=1}^{\kappa (v)} \mathcal{F}_{\lambda_i} (t)	.
	\end{equation}
	\end{demo}

\begin{exa} \rm
	Let us illustrate Theorem \ref{teo generalisation Pagliacci} for the following configuration 
	\begin{equation}
		\configurationD
	\end{equation}

If we denote by $\mathcal{C}$  the configuration above, then $\underline{v_{\mathcal{C}}}$ is given by

\begin{equation}
	s_{3}s_{2}s_{4}s_{6}s_{1}s_{3}s_{5}s_{7}s_{2}s_{4}s_{6}s_{8}s_{7}s_{9}s_{8}.
\end{equation}

In this case, $N_v^1=\{ s_1,s_5,s_9  \}$. Thus, $n_v^1=3$. Furthermore, two sub-configurations exist: $\mathcal{C}_1'$ and $\mathcal{C}_2'$, which are depicted in (\ref{equation confifi}). In this setting, we have $\kappa_v=2$ and $\lambda_1=\lambda_2=3$. 

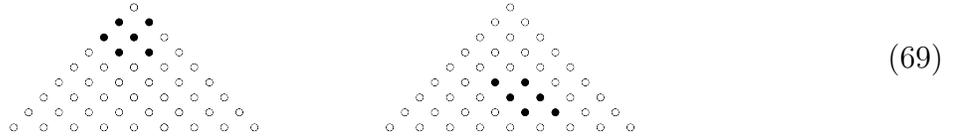
\begin{equation} \label{equation confifi}
	\confifi
\end{equation}

The words associated to the configurations $\mathcal{C}_1'$ and $\mathcal{C}_2'$ are given by 
\begin{equation}
	\underline{v_{\mathcal{C}_1'}} = s_{3}s_2s_4s_3s_2s_4 \qquad \mbox{ and } \qquad \underline{v_{\mathcal{C}_2'}}=s_6s_7s_6s_8s_7s_8.
\end{equation}

The sets $ \mathbb{L}_{\underline{v_{\mathcal{C}_1'}}}(e)$ and $ \mathbb{L}_{\underline{v_{\mathcal{C}_2'}}}(e)$ are obtained by considering the trees in (\ref{equation configuration C}), where because of space limitations, we have replaced the letter $s_i$ by $i$. 
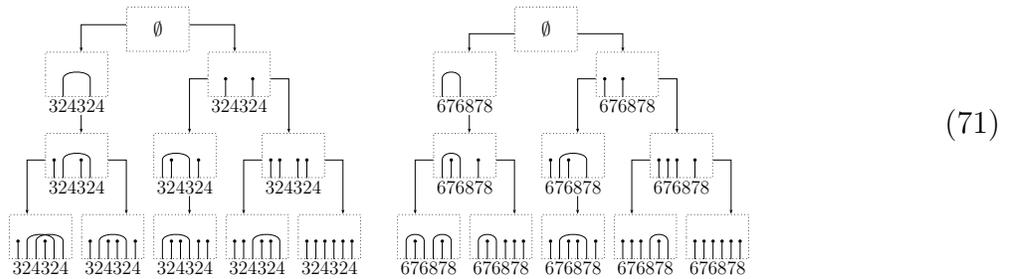
\begin{equation} \label{equation configuration C}
	\configurationC
\end{equation}  

We obtain

\begin{equation}
	\tilde{\mathcal{R}}_{e,\underline{v_{\mathcal{C}_1'}}}(t)=\tilde{\mathcal{R}}_{e,\underline{v_{\mathcal{C}_2'}}}(t) = \mathcal{F}_3(t)= v^6+3v^4+v^2.
\end{equation}

To obtain any element of $\mathbb{L}_{\underline{v_\mathcal{C}}}(e)$ we choose an element of $ \mathbb{L}_{\underline{v_{\mathcal{C}_1'}}}(e)$, and an element $ \mathbb{L}_{\underline{v_{\mathcal{C}_2'}}}(e)$, and combine them. At this point, it is crucial that the letters involved in $\mathcal{C}_1'$ and $\mathcal{C}_2'$ commute, because some lines in the diagrams might involve an intersection. Exactly $25$ elements appear in $\mathbb{L}_{\underline{v_\mathcal{C}}}(e)$. Because of space limitations, we draw only one of these $25$ elements in (\ref{equation last leaf}). This diagram is the element of the lower degree in $\mathbb{L}_{\underline{v_\mathcal{C}}}(e)$.

\begin{equation}\label{equation last leaf}
\lastleaf	
\end{equation}

Finally, we have

\begin{equation}
	\tilde{\mathcal{R}}_{e,v}(t)=\tilde{\mathcal{R}}_{e,\underline{v_{\mathcal{C}}}}(t)=t^3(v^6+3v^4+v^2)^2=t^{15} + 6t^{13} + 11t^{11} + 6t^9 + t^7.
\end{equation} 
\end{exa}

As a final remark, the common feature of the four formulas obtained in this section is that the relevant set $\mathbb{L}_{\underline{v}}(u)$ could be constructed without using $6$-valent vertices. However, this is not true in general. The occurrence of $6$-valent vertices (or more generally $2m_{st}$-valent vertices) allows for letters appearing in a word to interact between them, and gaining control on this interaction appears extremely challenging. This explains the difficulty in obtaining closed formulas for $\tilde{R}$-polynomials. Nevertheless, we believe that if a few $6$-valent vertices exist in the diagrams of $\mathbb{L}_{\underline{v}}(u)$, then the possibility of obtaining a closed formula for $\tilde{\mathcal{R}}_{u,\underline{v}}(t)$ still exists. In particular, we conjecture that the following holds.

\begin{conj}
	Let $w\in \mathfrak{S}_n$ be a $321$-avoiding and $2$-repeating permutation. If $u\leq v\leq w$, then
	\begin{equation}
		\tilde{R}_{u,v}(t)= t^{a}\prod_{i=1}^{b} \mathcal{F}_{c_i} (t),
	\end{equation}
	for some integers $a,b$ and $ c_i$.
\end{conj}

\section*{References}
\bibliographystyle{myalpha} 
\bibliography{mybibfile}

\end{document}

%% file: pictures.tex
\def\totalleaves{
\psscalebox{.5 .5} 
{
\begin{pspicture}[shift=-1.6](0,-1.6)(23.2,1.6)
\rput(2.6,-2.1){\scalebox{2}{$124579874321$}}
\rput(8.6,-2.1){\scalebox{2}{$124579874321$}}
\rput(14.6,-2.1){\scalebox{2}{$124579874321$}}
\rput(20.6,-2.1){\scalebox{2}{$124579874321$}}
\psframe[linecolor=black, linewidth=0.04, linestyle=dotted, dotsep=0.10583334cm, dimen=outer](5.2,1.6)(0.0,-1.6)
\psframe[linecolor=black, linewidth=0.04, linestyle=dotted, dotsep=0.10583334cm, dimen=outer](11.2,1.6)(6.0,-1.6)
\psframe[linecolor=black, linewidth=0.04, linestyle=dotted, dotsep=0.10583334cm, dimen=outer](17.2,1.6)(12.0,-1.6)
\psframe[linecolor=black, linewidth=0.04, linestyle=dotted, dotsep=0.10583334cm, dimen=outer](23.2,1.6)(18.0,-1.6)
\psline[linecolor=black, linewidth=0.04, dotsize=0.07055555cm 2.0]{-*}(0.8,-1.6)(0.8,-0.8)
\psline[linecolor=black, linewidth=0.04, dotsize=0.07055555cm 2.0]{-*}(1.6,-1.6)(1.6,-0.8)
\psline[linecolor=black, linewidth=0.04, dotsize=0.07055555cm 2.0]{-*}(4.4,-1.6)(4.4,-0.8)
\psline[linecolor=black, linewidth=0.04, dotsize=0.07055555cm 2.0]{-*}(3.2,-1.6)(3.2,-0.8)
\psline[linecolor=black, linewidth=0.04, dotsize=0.07055555cm 2.0]{-*}(6.8,-1.6)(6.8,-0.8)
\psline[linecolor=black, linewidth=0.04, dotsize=0.07055555cm 2.0]{-*}(7.6,-1.6)(7.6,-0.8)
\psline[linecolor=black, linewidth=0.04, dotsize=0.07055555cm 2.0]{-*}(10.4,-1.6)(10.4,-0.8)
\psline[linecolor=black, linewidth=0.04, dotsize=0.07055555cm 2.0]{-*}(9.2,-1.6)(9.2,-0.8)
\psline[linecolor=black, linewidth=0.04, dotsize=0.07055555cm 2.0]{-*}(12.8,-1.6)(12.8,-0.8)
\psline[linecolor=black, linewidth=0.04, dotsize=0.07055555cm 2.0]{-*}(13.6,-1.6)(13.6,-0.8)
\psline[linecolor=black, linewidth=0.04, dotsize=0.07055555cm 2.0]{-*}(16.4,-1.6)(16.4,-0.8)
\psline[linecolor=black, linewidth=0.04, dotsize=0.07055555cm 2.0]{-*}(15.2,-1.6)(15.2,-0.8)
\psline[linecolor=black, linewidth=0.04, dotsize=0.07055555cm 2.0]{-*}(18.8,-1.6)(18.8,-0.8)
\psline[linecolor=black, linewidth=0.04, dotsize=0.07055555cm 2.0]{-*}(19.6,-1.6)(19.6,-0.8)
\psline[linecolor=black, linewidth=0.04, dotsize=0.07055555cm 2.0]{-*}(22.4,-1.6)(22.4,-0.8)
\psline[linecolor=black, linewidth=0.04, dotsize=0.07055555cm 2.0]{-*}(21.2,-1.6)(21.2,-0.8)
\psline[linecolor=black, linewidth=0.04](2.0,-1.6)(2.0,1.6)
\psline[linecolor=black, linewidth=0.04](2.4,1.6)(2.4,-1.6)
\psline[linecolor=black, linewidth=0.04](2.8,1.6)(2.8,-1.6)
\psline[linecolor=black, linewidth=0.04](4.0,-1.6)(4.0,1.6)
\psline[linecolor=black, linewidth=0.04](8.0,-1.6)(8.0,1.6)
\psline[linecolor=black, linewidth=0.04](8.4,1.6)(8.4,-1.6)
\psline[linecolor=black, linewidth=0.04](8.8,-1.6)(8.8,1.6)
\psline[linecolor=black, linewidth=0.04](10.0,-1.6)(10.0,1.6)
\psline[linecolor=black, linewidth=0.04](14.0,-1.6)(14.0,1.6)
\psline[linecolor=black, linewidth=0.04](14.4,-1.6)(14.4,1.6)
\psline[linecolor=black, linewidth=0.04](14.8,-1.6)(14.8,1.6)
\psline[linecolor=black, linewidth=0.04](16.0,-1.6)(16.0,1.6)
\psline[linecolor=black, linewidth=0.04](20.0,-1.6)(20.0,1.6)
\psline[linecolor=black, linewidth=0.04](20.4,-1.6)(20.4,1.6)
\psline[linecolor=black, linewidth=0.04](20.8,-1.6)(20.8,1.6)
\psline[linecolor=black, linewidth=0.04](22.0,-1.6)(22.0,1.6)
\psbezier[linecolor=black, linewidth=0.04](1.2,-1.6)(1.2,-1.2)(1.2,-0.8)(1.2,-0.4)(1.2,0.0)(3.6,0.0)(3.6,-0.4)(3.6,-0.8)(3.6,-1.2)(3.6,-1.6)
\psbezier[linecolor=black, linewidth=0.04](7.2,-1.6)(7.2,-1.2)(7.2,-0.8)(7.2,-0.4)(7.2,0.0)(9.6,0.0)(9.6,-0.4)(9.6,-0.8)(9.6,-1.2)(9.6,-1.6)
\psbezier[linecolor=black, linewidth=0.04](12.4,-1.6)(12.4,-1.2)(12.4,-0.4)(12.4,0.0)(12.4,0.4)(16.8,0.4)(16.8,0.0)(16.8,-0.4)(16.8,-1.2)(16.8,-1.6)
\psbezier[linecolor=black, linewidth=0.04](0.4,-1.6)(0.4,-1.2)(0.4,-0.4)(0.4,0.0)(0.4,0.4)(4.8,0.4)(4.8,0.0)(4.8,-0.4)(4.8,-1.2)(4.8,-1.6)
\psline[linecolor=black, linewidth=0.04, dotsize=0.07055555cm 2.0]{-*}(6.4,-1.6)(6.4,-0.8)
\psline[linecolor=black, linewidth=0.04, dotsize=0.07055555cm 2.0]{-*}(10.8,-1.6)(10.8,-0.8)
\psline[linecolor=black, linewidth=0.04, dotsize=0.07055555cm 2.0]{-*}(13.2,-1.6)(13.2,-0.8)
\psline[linecolor=black, linewidth=0.04, dotsize=0.07055555cm 2.0]{-*}(15.6,-1.6)(15.6,-0.8)
\psline[linecolor=black, linewidth=0.04, dotsize=0.07055555cm 2.0]{-*}(18.4,-1.6)(18.4,-0.8)
\psline[linecolor=black, linewidth=0.04, dotsize=0.07055555cm 2.0]{-*}(19.2,-1.6)(19.2,-0.8)
\psline[linecolor=black, linewidth=0.04, dotsize=0.07055555cm 2.0]{-*}(21.6,-1.6)(21.6,-0.8)
\psline[linecolor=black, linewidth=0.04, dotsize=0.07055555cm 2.0]{-*}(22.8,-1.6)(22.8,-0.8)
\end{pspicture}
}
}

\def\lastleaf{
\psscalebox{.5 .5} 
{
\begin{pspicture}[shift=-1.8](0,-1.8)(6.4,1.8)
\rput(3.2,-2.2){\scalebox{1.45}{$3\,  2\, 4\, 6\, 1\, 3\, 5\, 7\, 2\, 4\, 6\, 8\, 7\, 9\,  8$}}
\psframe[linecolor=black, linewidth=0.04, linestyle=dotted, dotsep=0.10583334cm, dimen=outer](6.4,1.8)(0.0,-1.8)
\psline[linecolor=black, linewidth=0.04, dotsize=0.07055555cm 2.0]{-*}(2.0,-1.8)(2.0,-1.0)
\psline[linecolor=black, linewidth=0.04, dotsize=0.07055555cm 2.0]{-*}(5.6,-1.8)(5.6,-1.0)
\psline[linecolor=black, linewidth=0.04, dotsize=0.07055555cm 2.0]{-*}(0.4,-1.8)(0.4,-1.0)
\psline[linecolor=black, linewidth=0.04, dotsize=0.07055555cm 2.0]{-*}(2.4,-1.8)(2.4,-1.0)
\psbezier[linecolor=black, linewidth=0.04](0.8,-1.8)(0.8,-1.4)(0.8,-1.0)(0.8,-0.6)(0.8,-0.2)(3.6,-0.2)(3.6,-0.6)(3.6,-1.0)(3.6,-1.4)(3.6,-1.8)
\psbezier[linecolor=black, linewidth=0.04](1.2,-1.8)(1.2,-1.2666667)(1.2,-0.73333335)(1.2,-0.2)(1.2,0.33333334)(4.0,0.33333334)(4.0,-0.2)(4.0,-0.73333335)(4.0,-1.2666667)(4.0,-1.8)
\psbezier[linecolor=black, linewidth=0.04](1.6,-1.8)(1.6,-1.1333333)(1.6,-0.46666667)(1.6,0.2)(1.6,0.8666667)(4.4,0.8666667)(4.4,0.2)(4.4,-0.46666667)(4.4,-1.1333333)(4.4,-1.8)
\psbezier[linecolor=black, linewidth=0.04](4.8,-1.8)(4.8,-1.2666667)(4.8,-1.5333333)(4.8,-1.0)(4.8,-0.46666667)(6.0,-0.46666667)(6.0,-1.0)(6.0,-1.5333333)(6.0,-1.2666667)(6.0,-1.8)
\psline[linecolor=black, linewidth=0.04, dotsize=0.07055555cm 2.0]{-*}(2.8,-1.8)(2.8,-1.0)
\psline[linecolor=black, linewidth=0.04, dotsize=0.07055555cm 2.0]{-*}(3.2,-1.8)(3.2,-1.0)
\psline[linecolor=black, linewidth=0.04, dotsize=0.07055555cm 2.0]{-*}(5.2,-1.8)(5.2,-1.0)
\end{pspicture}
}
}

\def\confifi{
\psscalebox{.5 .5} 
{
\begin{pspicture}[shift=-1.7](0,-1.7006946)(16.601387,1.7006946)
\psdots[linecolor=black, dotstyle=o, dotsize=0.2, fillcolor=white](3.3006945,1.6)
\psdots[linecolor=black, dotstyle=o, dotsize=0.2, fillcolor=white](2.9006946,1.2)
\psdots[linecolor=black, dotstyle=o, dotsize=0.2, fillcolor=white](3.7006946,1.2)
\psdots[linecolor=black, dotstyle=o, dotsize=0.2, fillcolor=white](3.3006945,0.8000001)
\psdots[linecolor=black, dotstyle=o, dotsize=0.2, fillcolor=white](2.5006945,0.8000001)
\psdots[linecolor=black, dotstyle=o, dotsize=0.2, fillcolor=white](4.1006947,0.8000001)
\psdots[linecolor=black, dotstyle=o, dotsize=0.2, fillcolor=white](4.5006948,0.40000007)
\psdots[linecolor=black, dotstyle=o, dotsize=0.2, fillcolor=white](3.7006946,0.40000007)
\psdots[linecolor=black, dotstyle=o, dotsize=0.2, fillcolor=white](2.9006946,0.40000007)
\psdots[linecolor=black, dotstyle=o, dotsize=0.2, fillcolor=white](2.1006947,0.40000007)
\psdots[linecolor=black, dotstyle=o, dotsize=0.2, fillcolor=white](1.7006946,0.0)
\psdots[linecolor=black, dotstyle=o, dotsize=0.2, fillcolor=white](2.5006945,0.0)
\psdots[linecolor=black, dotstyle=o, dotsize=0.2, fillcolor=white](3.3006945,0.0)
\psdots[linecolor=black, dotstyle=o, dotsize=0.2, fillcolor=white](4.1006947,0.0)
\psdots[linecolor=black, dotstyle=o, dotsize=0.2, fillcolor=white](4.9006944,0.0)
\psdots[linecolor=black, dotstyle=o, dotsize=0.2, fillcolor=white](4.5006948,-0.39999992)
\psdots[linecolor=black, dotstyle=o, dotsize=0.2, fillcolor=white](5.3006945,-0.39999992)
\psdots[linecolor=black, dotstyle=o, dotsize=0.2, fillcolor=white](3.7006946,-0.39999992)
\psdots[linecolor=black, dotstyle=o, dotsize=0.2, fillcolor=white](2.9006946,-0.39999992)
\psdots[linecolor=black, dotstyle=o, dotsize=0.2, fillcolor=white](2.1006947,-0.39999992)
\psdots[linecolor=black, dotstyle=o, dotsize=0.2, fillcolor=white](1.3006946,-0.39999992)
\psdots[linecolor=black, dotstyle=o, dotsize=0.2, fillcolor=white](0.9006946,-0.79999995)
\psdots[linecolor=black, dotstyle=o, dotsize=0.2, fillcolor=white](0.5006946,-1.1999999)
\psdots[linecolor=black, dotstyle=o, dotsize=0.2, fillcolor=white](0.10069458,-1.5999999)
\psdots[linecolor=black, dotstyle=o, dotsize=0.2, fillcolor=white](0.9006946,-1.5999999)
\psdots[linecolor=black, dotstyle=o, dotsize=0.2, fillcolor=white](1.3006946,-1.1999999)
\psdots[linecolor=black, dotstyle=o, dotsize=0.2, fillcolor=white](1.7006946,-0.79999995)
\psdots[linecolor=black, dotstyle=o, dotsize=0.2, fillcolor=white](2.5006945,-0.79999995)
\psdots[linecolor=black, dotstyle=o, dotsize=0.2, fillcolor=white](2.1006947,-1.1999999)
\psdots[linecolor=black, dotstyle=o, dotsize=0.2, fillcolor=white](1.7006946,-1.5999999)
\psdots[linecolor=black, dotstyle=o, dotsize=0.2, fillcolor=white](2.5006945,-1.5999999)
\psdots[linecolor=black, dotstyle=o, dotsize=0.2, fillcolor=white](2.9006946,-1.1999999)
\psdots[linecolor=black, dotstyle=o, dotsize=0.2, fillcolor=white](3.3006945,-0.79999995)
\psdots[linecolor=black, dotstyle=o, dotsize=0.2, fillcolor=white](4.1006947,-0.79999995)
\psdots[linecolor=black, dotstyle=o, dotsize=0.2, fillcolor=white](3.7006946,-1.1999999)
\psdots[linecolor=black, dotstyle=o, dotsize=0.2, fillcolor=white](3.3006945,-1.5999999)
\psdots[linecolor=black, dotstyle=o, dotsize=0.2, fillcolor=white](4.1006947,-1.5999999)
\psdots[linecolor=black, dotstyle=o, dotsize=0.2, fillcolor=white](4.5006948,-1.1999999)
\psdots[linecolor=black, dotstyle=o, dotsize=0.2, fillcolor=white](4.9006944,-0.79999995)
\psdots[linecolor=black, dotstyle=o, dotsize=0.2, fillcolor=white](5.7006946,-0.79999995)
\psdots[linecolor=black, dotstyle=o, dotsize=0.2, fillcolor=white](5.3006945,-1.1999999)
\psdots[linecolor=black, dotstyle=o, dotsize=0.2, fillcolor=white](4.9006944,-1.5999999)
\psdots[linecolor=black, dotstyle=o, dotsize=0.2, fillcolor=white](6.1006947,-1.1999999)
\psdots[linecolor=black, dotstyle=o, dotsize=0.2, fillcolor=white](5.7006946,-1.5999999)
\psdots[linecolor=black, dotstyle=o, dotsize=0.2, fillcolor=white](6.5006948,-1.5999999)
\psdots[linecolor=black, dotsize=0.2](3.3006945,1.6)
\psdots[linecolor=black, dotsize=0.2](2.9006946,1.2)
\psdots[linecolor=black, dotsize=0.2](3.7006946,1.2)
\psdots[linecolor=black, dotsize=0.2](3.3006945,0.8000001)
\psdots[linecolor=black, dotsize=0.2](2.5006945,0.8000001)
\psdots[linecolor=black, dotsize=0.2](2.9006946,0.40000007)
\psdots[linecolor=black, dotsize=0.2](3.7006946,0.40000007)
\psdots[linecolor=black, dotsize=0.2](3.3006945,0.0)
\psdots[linecolor=black, dotsize=0.2](2.9006946,-0.39999992)
\psdots[linecolor=black, dotsize=0.2](3.7006946,-0.39999992)
\psdots[linecolor=black, dotsize=0.2](3.3006945,-0.79999995)
\psdots[linecolor=black, dotsize=0.2](4.1006947,-0.79999995)
\psdots[linecolor=black, dotsize=0.2](3.7006946,-1.1999999)
\psdots[linecolor=black, dotsize=0.2](4.5006948,-1.1999999)
\psdots[linecolor=black, dotsize=0.2](4.1006947,-1.5999999)
\psdots[linecolor=black, dotstyle=o, dotsize=0.2, fillcolor=white](13.300694,1.6)
\psdots[linecolor=black, dotstyle=o, dotsize=0.2, fillcolor=white](12.900695,1.2)
\psdots[linecolor=black, dotstyle=o, dotsize=0.2, fillcolor=white](12.500694,0.8000001)
\psdots[linecolor=black, dotstyle=o, dotsize=0.2, fillcolor=white](12.100695,0.40000007)
\psdots[linecolor=black, dotstyle=o, dotsize=0.2, fillcolor=white](11.700695,0.0)
\psdots[linecolor=black, dotstyle=o, dotsize=0.2, fillcolor=white](11.300694,-0.39999992)
\psdots[linecolor=black, dotstyle=o, dotsize=0.2, fillcolor=white](10.900695,-0.79999995)
\psdots[linecolor=black, dotstyle=o, dotsize=0.2, fillcolor=white](10.500694,-1.1999999)
\psdots[linecolor=black, dotstyle=o, dotsize=0.2, fillcolor=white](10.100695,-1.5999999)
\psdots[linecolor=black, dotstyle=o, dotsize=0.2, fillcolor=white](10.900695,-1.5999999)
\psdots[linecolor=black, dotstyle=o, dotsize=0.2, fillcolor=white](11.300694,-1.1999999)
\psdots[linecolor=black, dotstyle=o, dotsize=0.2, fillcolor=white](11.700695,-0.79999995)
\psdots[linecolor=black, dotstyle=o, dotsize=0.2, fillcolor=white](12.100695,-0.39999992)
\psdots[linecolor=black, dotstyle=o, dotsize=0.2, fillcolor=white](12.500694,0.0)
\psdots[linecolor=black, dotstyle=o, dotsize=0.2, fillcolor=white](12.900695,0.40000007)
\psdots[linecolor=black, dotstyle=o, dotsize=0.2, fillcolor=white](13.300694,0.8000001)
\psdots[linecolor=black, dotstyle=o, dotsize=0.2, fillcolor=white](13.700695,1.2)
\psdots[linecolor=black, dotstyle=o, dotsize=0.2, fillcolor=white](14.100695,0.8000001)
\psdots[linecolor=black, dotstyle=o, dotsize=0.2, fillcolor=white](14.500694,0.40000007)
\psdots[linecolor=black, dotstyle=o, dotsize=0.2, fillcolor=white](14.900695,0.0)
\psdots[linecolor=black, dotstyle=o, dotsize=0.2, fillcolor=white](15.300694,-0.39999992)
\psdots[linecolor=black, dotstyle=o, dotsize=0.2, fillcolor=white](15.700695,-0.79999995)
\psdots[linecolor=black, dotstyle=o, dotsize=0.2, fillcolor=white](16.100695,-1.1999999)
\psdots[linecolor=black, dotstyle=o, dotsize=0.2, fillcolor=white](16.500694,-1.5999999)
\psdots[linecolor=black, dotstyle=o, dotsize=0.2, fillcolor=white](15.700695,-1.5999999)
\psdots[linecolor=black, dotstyle=o, dotsize=0.2, fillcolor=white](14.900695,-1.5999999)
\psdots[linecolor=black, dotstyle=o, dotsize=0.2, fillcolor=white](14.100695,-1.5999999)
\psdots[linecolor=black, dotstyle=o, dotsize=0.2, fillcolor=white](13.300694,-1.5999999)
\psdots[linecolor=black, dotstyle=o, dotsize=0.2, fillcolor=white](12.500694,-1.5999999)
\psdots[linecolor=black, dotstyle=o, dotsize=0.2, fillcolor=white](11.700695,-1.5999999)
\psdots[linecolor=black, dotstyle=o, dotsize=0.2, fillcolor=white](12.100695,-1.1999999)
\psdots[linecolor=black, dotstyle=o, dotsize=0.2, fillcolor=white](12.900695,-1.1999999)
\psdots[linecolor=black, dotstyle=o, dotsize=0.2, fillcolor=white](13.700695,-1.1999999)
\psdots[linecolor=black, dotstyle=o, dotsize=0.2, fillcolor=white](14.500694,-1.1999999)
\psdots[linecolor=black, dotstyle=o, dotsize=0.2, fillcolor=white](15.300694,-1.1999999)
\psdots[linecolor=black, dotstyle=o, dotsize=0.2, fillcolor=white](14.900695,-0.79999995)
\psdots[linecolor=black, dotstyle=o, dotsize=0.2, fillcolor=white](14.500694,-0.39999992)
\psdots[linecolor=black, dotstyle=o, dotsize=0.2, fillcolor=white](14.100695,0.0)
\psdots[linecolor=black, dotstyle=o, dotsize=0.2, fillcolor=white](13.700695,0.40000007)
\psdots[linecolor=black, dotstyle=o, dotsize=0.2, fillcolor=white](13.300694,0.0)
\psdots[linecolor=black, dotstyle=o, dotsize=0.2, fillcolor=white](12.900695,-0.39999992)
\psdots[linecolor=black, dotstyle=o, dotsize=0.2, fillcolor=white](12.500694,-0.79999995)
\psdots[linecolor=black, dotstyle=o, dotsize=0.2, fillcolor=white](13.300694,-0.79999995)
\psdots[linecolor=black, dotstyle=o, dotsize=0.2, fillcolor=white](14.100695,-0.79999995)
\psdots[linecolor=black, dotstyle=o, dotsize=0.2, fillcolor=white](13.700695,-0.39999992)
\psdots[linecolor=black, dotstyle=o, dotsize=0.2, fillcolor=white](3.3006945,1.6)
\psdots[linecolor=black, dotstyle=o, dotsize=0.2, fillcolor=white](3.3006945,0.0)
\psdots[linecolor=black, dotstyle=o, dotsize=0.2, fillcolor=white](4.1006947,-1.5999999)
\psdots[linecolor=black, dotstyle=o, dotsize=0.2, fillcolor=white](2.9006946,-0.39999992)
\psdots[linecolor=black, dotstyle=o, dotsize=0.2, fillcolor=white](3.7006946,-0.39999992)
\psdots[linecolor=black, dotstyle=o, dotsize=0.2, fillcolor=white](3.3006945,-0.79999995)
\psdots[linecolor=black, dotstyle=o, dotsize=0.2, fillcolor=white](4.1006947,-0.79999995)
\psdots[linecolor=black, dotstyle=o, dotsize=0.2, fillcolor=white](3.7006946,-1.1999999)
\psdots[linecolor=black, dotstyle=o, dotsize=0.2, fillcolor=white](4.5006948,-1.1999999)
\psdots[linecolor=black, dotsize=0.2](12.900695,-0.39999992)
\psdots[linecolor=black, dotsize=0.2](13.700695,-0.39999992)
\psdots[linecolor=black, dotsize=0.2](13.300694,-0.79999995)
\psdots[linecolor=black, dotsize=0.2](14.100695,-0.79999995)
\psdots[linecolor=black, dotsize=0.2](13.700695,-1.1999999)
\psdots[linecolor=black, dotsize=0.2](14.500694,-1.1999999)
\end{pspicture}
}
}

\def\configurationD{
\psscalebox{.5 .5} 
{
\begin{pspicture}[shift=-1.7](0,-1.7006946)(6.601389,1.7006946)
\psdots[linecolor=black, dotstyle=o, dotsize=0.2, fillcolor=white](3.3006945,1.6)
\psdots[linecolor=black, dotstyle=o, dotsize=0.2, fillcolor=white](2.9006946,1.2)
\psdots[linecolor=black, dotstyle=o, dotsize=0.2, fillcolor=white](3.7006946,1.2)
\psdots[linecolor=black, dotstyle=o, dotsize=0.2, fillcolor=white](3.3006945,0.8000001)
\psdots[linecolor=black, dotstyle=o, dotsize=0.2, fillcolor=white](2.5006945,0.8000001)
\psdots[linecolor=black, dotstyle=o, dotsize=0.2, fillcolor=white](4.1006947,0.8000001)
\psdots[linecolor=black, dotstyle=o, dotsize=0.2, fillcolor=white](4.5006948,0.40000007)
\psdots[linecolor=black, dotstyle=o, dotsize=0.2, fillcolor=white](3.7006946,0.40000007)
\psdots[linecolor=black, dotstyle=o, dotsize=0.2, fillcolor=white](2.9006946,0.40000007)
\psdots[linecolor=black, dotstyle=o, dotsize=0.2, fillcolor=white](2.1006947,0.40000007)
\psdots[linecolor=black, dotstyle=o, dotsize=0.2, fillcolor=white](1.7006946,0.0)
\psdots[linecolor=black, dotstyle=o, dotsize=0.2, fillcolor=white](2.5006945,0.0)
\psdots[linecolor=black, dotstyle=o, dotsize=0.2, fillcolor=white](3.3006945,0.0)
\psdots[linecolor=black, dotstyle=o, dotsize=0.2, fillcolor=white](4.1006947,0.0)
\psdots[linecolor=black, dotstyle=o, dotsize=0.2, fillcolor=white](4.9006944,0.0)
\psdots[linecolor=black, dotstyle=o, dotsize=0.2, fillcolor=white](4.5006948,-0.39999992)
\psdots[linecolor=black, dotstyle=o, dotsize=0.2, fillcolor=white](5.3006945,-0.39999992)
\psdots[linecolor=black, dotstyle=o, dotsize=0.2, fillcolor=white](3.7006946,-0.39999992)
\psdots[linecolor=black, dotstyle=o, dotsize=0.2, fillcolor=white](2.9006946,-0.39999992)
\psdots[linecolor=black, dotstyle=o, dotsize=0.2, fillcolor=white](2.1006947,-0.39999992)
\psdots[linecolor=black, dotstyle=o, dotsize=0.2, fillcolor=white](1.3006946,-0.39999992)
\psdots[linecolor=black, dotstyle=o, dotsize=0.2, fillcolor=white](0.9006946,-0.79999995)
\psdots[linecolor=black, dotstyle=o, dotsize=0.2, fillcolor=white](0.5006946,-1.1999999)
\psdots[linecolor=black, dotstyle=o, dotsize=0.2, fillcolor=white](0.10069458,-1.5999999)
\psdots[linecolor=black, dotstyle=o, dotsize=0.2, fillcolor=white](0.9006946,-1.5999999)
\psdots[linecolor=black, dotstyle=o, dotsize=0.2, fillcolor=white](1.3006946,-1.1999999)
\psdots[linecolor=black, dotstyle=o, dotsize=0.2, fillcolor=white](1.7006946,-0.79999995)
\psdots[linecolor=black, dotstyle=o, dotsize=0.2, fillcolor=white](2.5006945,-0.79999995)
\psdots[linecolor=black, dotstyle=o, dotsize=0.2, fillcolor=white](2.1006947,-1.1999999)
\psdots[linecolor=black, dotstyle=o, dotsize=0.2, fillcolor=white](1.7006946,-1.5999999)
\psdots[linecolor=black, dotstyle=o, dotsize=0.2, fillcolor=white](2.5006945,-1.5999999)
\psdots[linecolor=black, dotstyle=o, dotsize=0.2, fillcolor=white](2.9006946,-1.1999999)
\psdots[linecolor=black, dotstyle=o, dotsize=0.2, fillcolor=white](3.3006945,-0.79999995)
\psdots[linecolor=black, dotstyle=o, dotsize=0.2, fillcolor=white](4.1006947,-0.79999995)
\psdots[linecolor=black, dotstyle=o, dotsize=0.2, fillcolor=white](3.7006946,-1.1999999)
\psdots[linecolor=black, dotstyle=o, dotsize=0.2, fillcolor=white](3.3006945,-1.5999999)
\psdots[linecolor=black, dotstyle=o, dotsize=0.2, fillcolor=white](4.1006947,-1.5999999)
\psdots[linecolor=black, dotstyle=o, dotsize=0.2, fillcolor=white](4.5006948,-1.1999999)
\psdots[linecolor=black, dotstyle=o, dotsize=0.2, fillcolor=white](4.9006944,-0.79999995)
\psdots[linecolor=black, dotstyle=o, dotsize=0.2, fillcolor=white](5.7006946,-0.79999995)
\psdots[linecolor=black, dotstyle=o, dotsize=0.2, fillcolor=white](5.3006945,-1.1999999)
\psdots[linecolor=black, dotstyle=o, dotsize=0.2, fillcolor=white](4.9006944,-1.5999999)
\psdots[linecolor=black, dotstyle=o, dotsize=0.2, fillcolor=white](6.1006947,-1.1999999)
\psdots[linecolor=black, dotstyle=o, dotsize=0.2, fillcolor=white](5.7006946,-1.5999999)
\psdots[linecolor=black, dotstyle=o, dotsize=0.2, fillcolor=white](6.5006948,-1.5999999)
\psdots[linecolor=black, dotsize=0.2](3.3006945,1.6)
\psdots[linecolor=black, dotsize=0.2](2.9006946,1.2)
\psdots[linecolor=black, dotsize=0.2](3.7006946,1.2)
\psdots[linecolor=black, dotsize=0.2](3.3006945,0.8000001)
\psdots[linecolor=black, dotsize=0.2](2.5006945,0.8000001)
\psdots[linecolor=black, dotsize=0.2](2.9006946,0.40000007)
\psdots[linecolor=black, dotsize=0.2](3.7006946,0.40000007)
\psdots[linecolor=black, dotsize=0.2](3.3006945,0.0)
\psdots[linecolor=black, dotsize=0.2](2.9006946,-0.39999992)
\psdots[linecolor=black, dotsize=0.2](3.7006946,-0.39999992)
\psdots[linecolor=black, dotsize=0.2](3.3006945,-0.79999995)
\psdots[linecolor=black, dotsize=0.2](4.1006947,-0.79999995)
\psdots[linecolor=black, dotsize=0.2](3.7006946,-1.1999999)
\psdots[linecolor=black, dotsize=0.2](4.5006948,-1.1999999)
\psdots[linecolor=black, dotsize=0.2](4.1006947,-1.5999999)
\end{pspicture}
}
}

\def\configurationC{
\psscalebox{.3 .3} 
{
\begin{pspicture}[shift=-5.6](0,-5.6)(32.8,5.6)
\rput(23.8,4.6){\scalebox{2}{$\emptyset$}}
\rput(6.6,4.6){\scalebox{2}{$\emptyset$}}
\rput(3,1.2){\scalebox{2}{$324324$}}
\rput(3,-2.4){\scalebox{2}{$324324$}}
\rput(10.2,1.2){\scalebox{2}{$324324$}}
\rput(7.8,-2.4){\scalebox{2}{$324324$}}
\rput(12.6,-2.4){\scalebox{2}{$324324$}}
\rput(7.8,-6){\scalebox{2}{$324324$}}
\rput(4.6,-6){\scalebox{2}{$324324$}}
\rput(1.4,-6){\scalebox{2}{$324324$}}
\rput(11,-6){\scalebox{2}{$324324$}}
\rput(14.2,-6){\scalebox{2}{$324324$}}
\rput(18.6,-6){\scalebox{2}{$676878$}}
\rput(21.8,-6){\scalebox{2}{$676878$}}
\rput(25,-6){\scalebox{2}{$676878$}}
\rput(28.2,-6){\scalebox{2}{$676878$}}
\rput(31.4,-6){\scalebox{2}{$676878$}}
\rput(25,-2.4){\scalebox{2}{$676878$}}
\rput(20.2,-2.4){\scalebox{2}{$676878$}}
\rput(29.8,-2.4){\scalebox{2}{$676878$}}
\rput(20.2,1.2){\scalebox{2}{$676878$}}
\rput(27.4,1.2){\scalebox{2}{$676878$}}
\psframe[linecolor=black, linewidth=0.04, linestyle=dotted, dotsep=0.10583334cm, dimen=outer](2.8,-3.6)(0.0,-5.6)
\psframe[linecolor=black, linewidth=0.04, linestyle=dotted, dotsep=0.10583334cm, dimen=outer](6.0,-3.6)(3.2,-5.6)
\psframe[linecolor=black, linewidth=0.04, linestyle=dotted, dotsep=0.10583334cm, dimen=outer](9.2,-3.6)(6.4,-5.6)
\psframe[linecolor=black, linewidth=0.04, linestyle=dotted, dotsep=0.10583334cm, dimen=outer](12.4,-3.6)(9.6,-5.6)
\psframe[linecolor=black, linewidth=0.04, linestyle=dotted, dotsep=0.10583334cm, dimen=outer](15.6,-3.6)(12.8,-5.6)
\psframe[linecolor=black, linewidth=0.04, linestyle=dotted, dotsep=0.10583334cm, dimen=outer](20.0,-3.6)(17.2,-5.6)
\psframe[linecolor=black, linewidth=0.04, linestyle=dotted, dotsep=0.10583334cm, dimen=outer](23.2,-3.6)(20.4,-5.6)
\psframe[linecolor=black, linewidth=0.04, linestyle=dotted, dotsep=0.10583334cm, dimen=outer](26.4,-3.6)(23.6,-5.6)
\psframe[linecolor=black, linewidth=0.04, linestyle=dotted, dotsep=0.10583334cm, dimen=outer](29.6,-3.6)(26.8,-5.6)
\psframe[linecolor=black, linewidth=0.04, linestyle=dotted, dotsep=0.10583334cm, dimen=outer](32.8,-3.6)(30.0,-5.6)
\psframe[linecolor=black, linewidth=0.04, linestyle=dotted, dotsep=0.10583334cm, dimen=outer](9.2,0.0)(6.4,-2.0)
\psframe[linecolor=black, linewidth=0.04, linestyle=dotted, dotsep=0.10583334cm, dimen=outer](4.4,0.0)(1.6,-2.0)
\psframe[linecolor=black, linewidth=0.04, linestyle=dotted, dotsep=0.10583334cm, dimen=outer](14.0,0.0)(11.2,-2.0)
\psframe[linecolor=black, linewidth=0.04, linestyle=dotted, dotsep=0.10583334cm, dimen=outer](11.6,3.6)(8.8,1.6)
\psframe[linecolor=black, linewidth=0.04, linestyle=dotted, dotsep=0.10583334cm, dimen=outer](4.4,3.6)(1.6,1.6)
\psframe[linecolor=black, linewidth=0.04, linestyle=dotted, dotsep=0.10583334cm, dimen=outer](8.0,5.6)(5.2,3.6)
\psframe[linecolor=black, linewidth=0.04, linestyle=dotted, dotsep=0.10583334cm, dimen=outer](21.6,0.0)(18.8,-2.0)
\psframe[linecolor=black, linewidth=0.04, linestyle=dotted, dotsep=0.10583334cm, dimen=outer](26.4,0.0)(23.6,-2.0)
\psframe[linecolor=black, linewidth=0.04, linestyle=dotted, dotsep=0.10583334cm, dimen=outer](31.2,0.0)(28.4,-2.0)
\psframe[linecolor=black, linewidth=0.04, linestyle=dotted, dotsep=0.10583334cm, dimen=outer](21.6,3.6)(18.8,1.6)
\psframe[linecolor=black, linewidth=0.04, linestyle=dotted, dotsep=0.10583334cm, dimen=outer](28.8,3.6)(26.0,1.6)
\psframe[linecolor=black, linewidth=0.04, linestyle=dotted, dotsep=0.10583334cm, dimen=outer](25.2,5.6)(22.4,3.6)
\psline[linecolor=black, linewidth=0.04, arrowsize=0.05291667cm 2.0,arrowlength=1.4,arrowinset=0.0]{->}(5.2,4.8)(3.2,4.8)(3.2,3.6)
\psline[linecolor=black, linewidth=0.04, arrowsize=0.05291667cm 2.0,arrowlength=1.4,arrowinset=0.0]{->}(8.0,4.8)(10.0,4.8)(10.0,3.6)
\psline[linecolor=black, linewidth=0.04, arrowsize=0.05291667cm 2.0,arrowlength=1.4,arrowinset=0.0]{->}(3.2,0.8)(3.2,0.0)
\psline[linecolor=black, linewidth=0.04, arrowsize=0.05291667cm 2.0,arrowlength=1.4,arrowinset=0.0]{->}(1.6,-1.2)(0.8,-1.2)(0.8,-3.6)
\psline[linecolor=black, linewidth=0.04, arrowsize=0.05291667cm 2.0,arrowlength=1.4,arrowinset=0.0]{->}(4.4,-1.2)(5.2,-1.2)(5.2,-3.6)
\psline[linecolor=black, linewidth=0.04, arrowsize=0.05291667cm 2.0,arrowlength=1.4,arrowinset=0.0]{->}(8.0,-2.8)(8.0,-3.6)
\psline[linecolor=black, linewidth=0.04, arrowsize=0.05291667cm 2.0,arrowlength=1.4,arrowinset=0.0]{->}(11.2,-1.2)(10.4,-1.2)(10.4,-3.6)
\psline[linecolor=black, linewidth=0.04, arrowsize=0.05291667cm 2.0,arrowlength=1.4,arrowinset=0.0]{->}(14.0,-1.2)(14.8,-1.2)(14.8,-3.6)
\psline[linecolor=black, linewidth=0.04, arrowsize=0.05291667cm 2.0,arrowlength=1.4,arrowinset=0.0]{->}(8.8,2.4)(8.0,2.4)(8.0,0.0)
\psline[linecolor=black, linewidth=0.04, arrowsize=0.05291667cm 2.0,arrowlength=1.4,arrowinset=0.0]{->}(11.6,2.4)(12.4,2.4)(12.4,0.0)
\psline[linecolor=black, linewidth=0.04, arrowsize=0.05291667cm 2.0,arrowlength=1.4,arrowinset=0.0]{->}(18.8,-1.2)(18.0,-1.2)(18.0,-3.6)
\psline[linecolor=black, linewidth=0.04, arrowsize=0.05291667cm 2.0,arrowlength=1.4,arrowinset=0.0]{->}(21.6,-1.2)(22.4,-1.2)(22.4,-3.6)
\psline[linecolor=black, linewidth=0.04, arrowsize=0.05291667cm 2.0,arrowlength=1.4,arrowinset=0.0]{->}(25.2,-2.8)(25.2,-3.6)
\psline[linecolor=black, linewidth=0.04, arrowsize=0.05291667cm 2.0,arrowlength=1.4,arrowinset=0.0]{->}(28.4,-1.2)(27.6,-1.2)(27.6,-3.6)
\psline[linecolor=black, linewidth=0.04, arrowsize=0.05291667cm 2.0,arrowlength=1.4,arrowinset=0.0]{->}(31.2,-1.2)(32.0,-1.2)(32.0,-3.6)
\psline[linecolor=black, linewidth=0.04, arrowsize=0.05291667cm 2.0,arrowlength=1.4,arrowinset=0.0]{->}(26.0,2.4)(25.2,2.4)(25.2,0.0)
\psline[linecolor=black, linewidth=0.04, arrowsize=0.05291667cm 2.0,arrowlength=1.4,arrowinset=0.0]{->}(28.8,2.4)(29.6,2.4)(29.6,0.0)
\psline[linecolor=black, linewidth=0.04, arrowsize=0.05291667cm 2.0,arrowlength=1.4,arrowinset=0.0]{->}(20.4,0.8)(20.4,0.0)
\psline[linecolor=black, linewidth=0.04, arrowsize=0.05291667cm 2.0,arrowlength=1.4,arrowinset=0.0]{->}(22.4,4.4)(20.4,4.4)(20.4,3.6)
\psline[linecolor=black, linewidth=0.04, arrowsize=0.05291667cm 2.0,arrowlength=1.4,arrowinset=0.0]{->}(25.2,4.4)(27.2,4.4)(27.2,3.6)
\psline[linecolor=black, linewidth=0.04, dotsize=0.07055555cm 2.0]{-*}(9.6,1.6)(9.6,2.4)
\psline[linecolor=black, linewidth=0.04, dotsize=0.07055555cm 2.0]{-*}(10.8,1.6)(10.8,2.4)
\psline[linecolor=black, linewidth=0.04, dotsize=0.07055555cm 2.0]{-*}(7.2,-2.0)(7.2,-1.2)
\psline[linecolor=black, linewidth=0.04, dotsize=0.07055555cm 2.0]{-*}(8.4,-2.0)(8.4,-1.2)
\psline[linecolor=black, linewidth=0.04, dotsize=0.07055555cm 2.0]{-*}(12.0,-2.0)(12.0,-1.2)
\psline[linecolor=black, linewidth=0.04, dotsize=0.07055555cm 2.0]{-*}(13.2,-2.0)(13.2,-1.2)
\psline[linecolor=black, linewidth=0.04, dotsize=0.07055555cm 2.0]{-*}(7.2,-5.6)(7.2,-4.8)
\psline[linecolor=black, linewidth=0.04, dotsize=0.07055555cm 2.0]{-*}(8.4,-5.6)(8.4,-4.8)
\psline[linecolor=black, linewidth=0.04, dotsize=0.07055555cm 2.0]{-*}(10.4,-5.6)(10.4,-4.8)
\psline[linecolor=black, linewidth=0.04, dotsize=0.07055555cm 2.0]{-*}(11.6,-5.6)(11.6,-4.8)
\psline[linecolor=black, linewidth=0.04, dotsize=0.07055555cm 2.0]{-*}(13.6,-5.6)(13.6,-4.8)
\psline[linecolor=black, linewidth=0.04, dotsize=0.07055555cm 2.0]{-*}(14.8,-5.6)(14.8,-4.8)
\psline[linecolor=black, linewidth=0.04, dotsize=0.07055555cm 2.0]{-*}(14.0,-5.6)(14.0,-4.8)
\psline[linecolor=black, linewidth=0.04, dotsize=0.07055555cm 2.0]{-*}(15.2,-5.6)(15.2,-4.8)
\psline[linecolor=black, linewidth=0.04, dotsize=0.07055555cm 2.0]{-*}(14.4,-5.6)(14.4,-4.8)
\psline[linecolor=black, linewidth=0.04, dotsize=0.07055555cm 2.0]{-*}(13.2,-5.6)(13.2,-4.8)
\psline[linecolor=black, linewidth=0.04, dotsize=0.07055555cm 2.0]{-*}(11.6,-2.0)(11.6,-1.2)
\psline[linecolor=black, linewidth=0.04, dotsize=0.07055555cm 2.0]{-*}(12.8,-2.0)(12.8,-1.2)
\psline[linecolor=black, linewidth=0.04, dotsize=0.07055555cm 2.0]{-*}(10.0,-5.6)(10.0,-4.8)
\psline[linecolor=black, linewidth=0.04, dotsize=0.07055555cm 2.0]{-*}(11.2,-5.6)(11.2,-4.8)
\psline[linecolor=black, linewidth=0.04, dotsize=0.07055555cm 2.0]{-*}(7.6,-5.6)(7.6,-4.8)
\psline[linecolor=black, linewidth=0.04, dotsize=0.07055555cm 2.0]{-*}(8.8,-5.6)(8.8,-4.8)
\psline[linecolor=black, linewidth=0.04, dotsize=0.07055555cm 2.0]{-*}(2.0,-2.0)(2.0,-1.2)
\psline[linecolor=black, linewidth=0.04, dotsize=0.07055555cm 2.0]{-*}(3.2,-2.0)(3.2,-1.2)
\psline[linecolor=black, linewidth=0.04, dotsize=0.07055555cm 2.0]{-*}(3.6,-5.6)(3.6,-4.8)
\psline[linecolor=black, linewidth=0.04, dotsize=0.07055555cm 2.0]{-*}(4.4,-5.6)(4.4,-4.8)
\psline[linecolor=black, linewidth=0.04, dotsize=0.07055555cm 2.0]{-*}(4.8,-5.6)(4.8,-4.8)
\psline[linecolor=black, linewidth=0.04, dotsize=0.07055555cm 2.0]{-*}(5.6,-5.6)(5.6,-4.8)
\psline[linecolor=black, linewidth=0.04, dotsize=0.07055555cm 2.0]{-*}(0.4,-5.6)(0.4,-4.8)
\psline[linecolor=black, linewidth=0.04, dotsize=0.07055555cm 2.0]{-*}(1.6,-5.6)(1.6,-4.8)
\psbezier[linecolor=black, linewidth=0.04](2.4,1.6)(2.4,2.0)(2.4,2.0)(2.4,2.4)(2.4,2.8)(3.6,2.8)(3.6,2.4)(3.6,2.0)(3.6,2.4)(3.6,1.6)
\psbezier[linecolor=black, linewidth=0.04](2.4,-2.0)(2.4,-1.6)(2.4,-1.6)(2.4,-1.2)(2.4,-0.8)(3.6,-0.8)(3.6,-1.2)(3.6,-1.6)(3.6,-1.6)(3.6,-2.0)
\psbezier[linecolor=black, linewidth=0.04](6.8,-2.0)(6.8,-1.6)(6.8,-1.6)(6.8,-1.2)(6.8,-0.8)(8.0,-0.8)(8.0,-1.2)(8.0,-1.6)(8.0,-1.6)(8.0,-2.0)
\psbezier[linecolor=black, linewidth=0.04](0.8,-5.6)(0.8,-4.8)(0.8,-5.2)(0.8,-4.8)(0.8,-4.4)(2.0,-4.4)(2.0,-4.8)(2.0,-5.2)(2.0,-4.8)(2.0,-5.6)
\psbezier[linecolor=black, linewidth=0.04](1.2,-5.6)(1.2,-5.2)(1.2,-5.2)(1.2,-4.8)(1.2,-4.4)(2.4,-4.4)(2.4,-4.8)(2.4,-5.2)(2.4,-4.8)(2.4,-5.6)
\psbezier[linecolor=black, linewidth=0.04](4.0,-5.6)(4.0,-5.2)(4.0,-5.2)(4.0,-4.8)(4.0,-4.4)(5.2,-4.4)(5.2,-4.8)(5.2,-5.2)(5.2,-4.8)(5.2,-5.6)
\psbezier[linecolor=black, linewidth=0.04](6.8,-5.6)(6.8,-5.2)(6.8,-5.2)(6.8,-4.8)(6.8,-4.4)(8.0,-4.4)(8.0,-4.8)(8.0,-5.2)(8.0,-4.8)(8.0,-5.6)
\psbezier[linecolor=black, linewidth=0.04](10.8,-5.6)(10.8,-5.2)(10.8,-5.2)(10.8,-4.8)(10.8,-4.4)(12.0,-4.4)(12.0,-4.8)(12.0,-5.2)(12.0,-4.8)(12.0,-5.6)
\psline[linecolor=black, linewidth=0.04, dotsize=0.07055555cm 2.0]{-*}(26.4,1.6)(26.4,2.4)
\psline[linecolor=black, linewidth=0.04, dotsize=0.07055555cm 2.0]{-*}(27.2,1.6)(27.2,2.4)
\psline[linecolor=black, linewidth=0.04, dotsize=0.07055555cm 2.0]{-*}(24.0,-2.0)(24.0,-1.2)
\psline[linecolor=black, linewidth=0.04, dotsize=0.07055555cm 2.0]{-*}(24.8,-2.0)(24.8,-1.2)
\psline[linecolor=black, linewidth=0.04, dotsize=0.07055555cm 2.0]{-*}(28.8,-2.0)(28.8,-1.2)
\psline[linecolor=black, linewidth=0.04, dotsize=0.07055555cm 2.0]{-*}(29.6,-2.0)(29.6,-1.2)
\psline[linecolor=black, linewidth=0.04, dotsize=0.07055555cm 2.0]{-*}(24.0,-5.6)(24.0,-4.8)
\psline[linecolor=black, linewidth=0.04, dotsize=0.07055555cm 2.0]{-*}(24.8,-5.6)(24.8,-4.8)
\psline[linecolor=black, linewidth=0.04, dotsize=0.07055555cm 2.0]{-*}(27.2,-5.6)(27.2,-4.8)
\psline[linecolor=black, linewidth=0.04, dotsize=0.07055555cm 2.0]{-*}(28.0,-5.6)(28.0,-4.8)
\psline[linecolor=black, linewidth=0.04, dotsize=0.07055555cm 2.0]{-*}(30.4,-5.6)(30.4,-4.8)
\psline[linecolor=black, linewidth=0.04, dotsize=0.07055555cm 2.0]{-*}(31.2,-5.6)(31.2,-4.8)
\psline[linecolor=black, linewidth=0.04, dotsize=0.07055555cm 2.0]{-*}(29.2,-2.0)(29.2,-1.2)
\psline[linecolor=black, linewidth=0.04, dotsize=0.07055555cm 2.0]{-*}(30.4,-2.0)(30.4,-1.2)
\psline[linecolor=black, linewidth=0.04, dotsize=0.07055555cm 2.0]{-*}(30.8,-5.6)(30.8,-4.8)
\psline[linecolor=black, linewidth=0.04, dotsize=0.07055555cm 2.0]{-*}(31.6,-5.6)(31.6,-4.8)
\psline[linecolor=black, linewidth=0.04, dotsize=0.07055555cm 2.0]{-*}(32.0,-5.6)(32.0,-4.8)
\psline[linecolor=black, linewidth=0.04, dotsize=0.07055555cm 2.0]{-*}(32.4,-5.6)(32.4,-4.8)
\psline[linecolor=black, linewidth=0.04, dotsize=0.07055555cm 2.0]{-*}(27.6,-5.6)(27.6,-4.8)
\psline[linecolor=black, linewidth=0.04, dotsize=0.07055555cm 2.0]{-*}(28.8,-5.6)(28.8,-4.8)
\psline[linecolor=black, linewidth=0.04, dotsize=0.07055555cm 2.0]{-*}(25.2,-5.6)(25.2,-4.8)
\psline[linecolor=black, linewidth=0.04, dotsize=0.07055555cm 2.0]{-*}(26.0,-5.6)(26.0,-4.8)
\psline[linecolor=black, linewidth=0.04, dotsize=0.07055555cm 2.0]{-*}(19.6,-2.0)(19.6,-1.2)
\psline[linecolor=black, linewidth=0.04, dotsize=0.07055555cm 2.0]{-*}(20.8,-2.0)(20.8,-1.2)
\psline[linecolor=black, linewidth=0.04, dotsize=0.07055555cm 2.0]{-*}(18.0,-5.6)(18.0,-4.8)
\psline[linecolor=black, linewidth=0.04, dotsize=0.07055555cm 2.0]{-*}(19.2,-5.6)(19.2,-4.8)
\psline[linecolor=black, linewidth=0.04, dotsize=0.07055555cm 2.0]{-*}(21.2,-5.6)(21.2,-4.8)
\psline[linecolor=black, linewidth=0.04, dotsize=0.07055555cm 2.0]{-*}(22.4,-5.6)(22.4,-4.8)
\psline[linecolor=black, linewidth=0.04, dotsize=0.07055555cm 2.0]{-*}(22.0,-5.6)(22.0,-4.8)
\psline[linecolor=black, linewidth=0.04, dotsize=0.07055555cm 2.0]{-*}(22.8,-5.6)(22.8,-4.8)
\psbezier[linecolor=black, linewidth=0.04](19.2,1.6)(19.2,2.4)(19.2,2.0)(19.2,2.4)(19.2,2.8)(20.0,2.8)(20.0,2.4)(20.0,2.0)(20.0,2.4)(20.0,1.6)
\psbezier[linecolor=black, linewidth=0.04](19.2,-2.0)(19.2,-1.2)(19.2,-1.6)(19.2,-1.2)(19.2,-0.8)(20.0,-0.8)(20.0,-1.2)(20.0,-1.6)(20.0,-1.2)(20.0,-2.0)
\psbezier[linecolor=black, linewidth=0.04](24.4,-2.0)(24.4,-1.6)(24.4,-1.6)(24.4,-1.2)(24.4,-0.8)(25.6,-0.8)(25.6,-1.2)(25.6,-1.6)(25.6,-1.2)(25.6,-2.0)
\psbezier[linecolor=black, linewidth=0.04](17.6,-5.6)(17.6,-4.8)(17.6,-5.2)(17.6,-4.8)(17.6,-4.4)(18.4,-4.4)(18.4,-4.8)(18.4,-5.2)(18.4,-4.8)(18.4,-5.6)
\psbezier[linecolor=black, linewidth=0.04](20.8,-5.6)(20.8,-4.8)(20.8,-5.2)(20.8,-4.8)(20.8,-4.4)(21.6,-4.4)(21.6,-4.8)(21.6,-5.2)(21.6,-4.8)(21.6,-5.6)
\psbezier[linecolor=black, linewidth=0.04](18.8,-5.6)(18.8,-4.8)(18.8,-5.2)(18.8,-4.8)(18.8,-4.4)(19.6,-4.4)(19.6,-4.8)(19.6,-5.2)(19.6,-4.8)(19.6,-5.6)
\psbezier[linecolor=black, linewidth=0.04](24.4,-5.6)(24.4,-5.2)(24.4,-5.2)(24.4,-4.8)(24.4,-4.4)(25.6,-4.4)(25.6,-4.8)(25.6,-5.2)(25.6,-5.2)(25.6,-5.6)
\psbezier[linecolor=black, linewidth=0.04](28.4,-5.6)(28.4,-4.8)(28.4,-5.2)(28.4,-4.8)(28.4,-4.4)(29.2,-4.4)(29.2,-4.8)(29.2,-5.2)(29.2,-4.8)(29.2,-5.6)
\end{pspicture}
}
}

\def\configurationB{
\psscalebox{.5 .5} 
{
\begin{pspicture}[shift=-1.9](0,-1.8985577)(2.5971155,1.8985577)
\psdots[linecolor=black, dotsize=0.2](0.8985577,1.8000001)
\psdots[linecolor=black, dotsize=0.2](1.6985577,1.8000001)
\psdots[linecolor=black, dotsize=0.2](1.2985578,1.4000001)
\psdots[linecolor=black, dotsize=0.2](2.0985577,1.4000001)
\psdots[linecolor=black, dotsize=0.2](1.6985577,1.0)
\psdots[linecolor=black, dotsize=0.2](2.4985578,1.0)
\psdots[linecolor=black, dotsize=0.2](1.2985578,0.6)
\psdots[linecolor=black, dotsize=0.2](2.0985577,0.6)
\psdots[linecolor=black, dotsize=0.2](1.6985577,0.20000003)
\psdots[linecolor=black, dotsize=0.2](0.8985577,0.20000003)
\psdots[linecolor=black, dotsize=0.2](1.2985578,-0.19999996)
\psdots[linecolor=black, dotsize=0.2](2.0985577,-0.19999996)
\psdots[linecolor=black, dotsize=0.2](1.6985577,-0.59999996)
\psdots[linecolor=black, dotsize=0.2](0.8985577,-0.59999996)
\psdots[linecolor=black, dotsize=0.2](1.2985578,-0.99999994)
\psdots[linecolor=black, dotsize=0.2](0.49855775,-0.99999994)
\psdots[linecolor=black, dotsize=0.2](0.8985577,-1.4)
\psdots[linecolor=black, dotsize=0.2](0.09855774,-1.4)
\psdots[linecolor=black, dotsize=0.2](0.49855775,-1.8)
\psdots[linecolor=black, dotsize=0.2](1.2985578,-1.8)
\end{pspicture}
}
}

\def\configurationA{
\psscalebox{.5 .5} 
{
\begin{pspicture}[shift=-1.5](0,-1.5006945)(5.801389,1.5006945)
\psdots[linecolor=black, dotstyle=o, dotsize=0.2, fillcolor=white](2.9006946,1.4000001)
\psdots[linecolor=black, dotstyle=o, dotsize=0.2, fillcolor=white](2.5006945,1.0000001)
\psdots[linecolor=black, dotstyle=o, dotsize=0.2, fillcolor=white](3.3006945,1.0000001)
\psdots[linecolor=black, dotstyle=o, dotsize=0.2, fillcolor=white](2.9006946,0.6000001)
\psdots[linecolor=black, dotstyle=o, dotsize=0.2, fillcolor=white](2.1006947,0.6000001)
\psdots[linecolor=black, dotstyle=o, dotsize=0.2, fillcolor=white](3.7006946,0.6000001)
\psdots[linecolor=black, dotstyle=o, dotsize=0.2, fillcolor=white](3.3006945,0.20000008)
\psdots[linecolor=black, dotstyle=o, dotsize=0.2, fillcolor=white](2.5006945,0.20000008)
\psdots[linecolor=black, dotstyle=o, dotsize=0.2, fillcolor=white](1.7006946,0.20000008)
\psdots[linecolor=black, dotstyle=o, dotsize=0.2, fillcolor=white](4.1006947,0.20000008)
\psdots[linecolor=black, dotstyle=o, dotsize=0.2, fillcolor=white](3.7006946,-0.19999993)
\psdots[linecolor=black, dotstyle=o, dotsize=0.2, fillcolor=white](2.9006946,-0.19999993)
\psdots[linecolor=black, dotstyle=o, dotsize=0.2, fillcolor=white](2.1006947,-0.19999993)
\psdots[linecolor=black, dotstyle=o, dotsize=0.2, fillcolor=white](1.3006946,-0.19999993)
\psdots[linecolor=black, dotstyle=o, dotsize=0.2, fillcolor=white](4.5006948,-0.19999993)
\psdots[linecolor=black, dotstyle=o, dotsize=0.2, fillcolor=white](4.9006944,-0.5999999)
\psdots[linecolor=black, dotstyle=o, dotsize=0.2, fillcolor=white](4.1006947,-0.5999999)
\psdots[linecolor=black, dotstyle=o, dotsize=0.2, fillcolor=white](3.3006945,-0.5999999)
\psdots[linecolor=black, dotstyle=o, dotsize=0.2, fillcolor=white](2.5006945,-0.5999999)
\psdots[linecolor=black, dotstyle=o, dotsize=0.2, fillcolor=white](1.7006946,-0.5999999)
\psdots[linecolor=black, dotstyle=o, dotsize=0.2, fillcolor=white](0.9006946,-0.5999999)
\psdots[linecolor=black, dotstyle=o, dotsize=0.2, fillcolor=white](0.5006946,-0.99999994)
\psdots[linecolor=black, dotstyle=o, dotsize=0.2, fillcolor=white](1.3006946,-0.99999994)
\psdots[linecolor=black, dotstyle=o, dotsize=0.2, fillcolor=white](2.1006947,-0.99999994)
\psdots[linecolor=black, dotstyle=o, dotsize=0.2, fillcolor=white](2.9006946,-0.99999994)
\psdots[linecolor=black, dotstyle=o, dotsize=0.2, fillcolor=white](3.7006946,-0.99999994)
\psdots[linecolor=black, dotstyle=o, dotsize=0.2, fillcolor=white](4.5006948,-0.99999994)
\psdots[linecolor=black, dotstyle=o, dotsize=0.2, fillcolor=white](5.3006945,-0.99999994)
\psdots[linecolor=black, dotstyle=o, dotsize=0.2, fillcolor=white](5.7006946,-1.4)
\psdots[linecolor=black, dotstyle=o, dotsize=0.2, fillcolor=white](4.9006944,-1.4)
\psdots[linecolor=black, dotstyle=o, dotsize=0.2, fillcolor=white](4.1006947,-1.4)
\psdots[linecolor=black, dotstyle=o, dotsize=0.2, fillcolor=white](3.3006945,-1.4)
\psdots[linecolor=black, dotstyle=o, dotsize=0.2, fillcolor=white](2.5006945,-1.4)
\psdots[linecolor=black, dotstyle=o, dotsize=0.2, fillcolor=white](1.7006946,-1.4)
\psdots[linecolor=black, dotstyle=o, dotsize=0.2, fillcolor=white](0.9006946,-1.4)
\psdots[linecolor=black, dotstyle=o, dotsize=0.2, fillcolor=white](0.10069458,-1.4)
\psdots[linecolor=black, dotsize=0.2](2.9006946,1.4000001)
\psdots[linecolor=black, dotsize=0.2](3.3006945,1.0000001)
\psdots[linecolor=black, dotsize=0.2](2.5006945,1.0000001)
\psdots[linecolor=black, dotsize=0.2](2.9006946,0.6000001)
\psdots[linecolor=black, dotsize=0.2](3.7006946,0.6000001)
\psdots[linecolor=black, dotsize=0.2](3.3006945,0.20000008)
\psdots[linecolor=black, dotsize=0.2](4.1006947,0.20000008)
\psdots[linecolor=black, dotsize=0.2](3.7006946,-0.19999993)
\psdots[linecolor=black, dotsize=0.2](4.5006948,-0.19999993)
\psdots[linecolor=black, dotsize=0.2](4.1006947,-0.5999999)
\psdots[linecolor=black, dotsize=0.2](4.9006944,-0.5999999)
\psdots[linecolor=black, dotsize=0.2](4.5006948,-0.99999994)
\psdots[linecolor=black, dotsize=0.2](5.3006945,-0.99999994)
\psdots[linecolor=black, dotsize=0.2](4.9006944,-1.4)
\end{pspicture}
}
}

\def\configuration{
\psscalebox{.5 .5} 
{
\begin{pspicture}[shift=-3.7](0,-3.7006946)(14.601389,3.7006946)
\psdots[linecolor=black, dotstyle=o, dotsize=0.2, fillcolor=white](7.3006945,3.6000001)
\psdots[linecolor=black, dotstyle=o, dotsize=0.2, fillcolor=white](6.9006944,3.2)
\psdots[linecolor=black, dotstyle=o, dotsize=0.2, fillcolor=white](7.7006946,3.2)
\psdots[linecolor=black, dotstyle=o, dotsize=0.2, fillcolor=white](7.3006945,2.8000002)
\psdots[linecolor=black, dotstyle=o, dotsize=0.2, fillcolor=white](6.5006943,2.8000002)
\psdots[linecolor=black, dotstyle=o, dotsize=0.2, fillcolor=white](8.100695,2.8000002)
\psdots[linecolor=black, dotstyle=o, dotsize=0.2, fillcolor=white](7.7006946,2.4)
\psdots[linecolor=black, dotstyle=o, dotsize=0.2, fillcolor=white](6.9006944,2.4)
\psdots[linecolor=black, dotstyle=o, dotsize=0.2, fillcolor=white](6.1006947,2.4)
\psdots[linecolor=black, dotstyle=o, dotsize=0.2, fillcolor=white](8.500694,2.4)
\psdots[linecolor=black, dotstyle=o, dotsize=0.2, fillcolor=white](8.100695,2.0)
\psdots[linecolor=black, dotstyle=o, dotsize=0.2, fillcolor=white](7.3006945,2.0)
\psdots[linecolor=black, dotstyle=o, dotsize=0.2, fillcolor=white](6.5006943,2.0)
\psdots[linecolor=black, dotstyle=o, dotsize=0.2, fillcolor=white](5.7006946,2.0)
\psdots[linecolor=black, dotstyle=o, dotsize=0.2, fillcolor=white](8.900695,2.0)
\psdots[linecolor=black, dotstyle=o, dotsize=0.2, fillcolor=white](9.300694,1.6)
\psdots[linecolor=black, dotstyle=o, dotsize=0.2, fillcolor=white](8.500694,1.6)
\psdots[linecolor=black, dotstyle=o, dotsize=0.2, fillcolor=white](7.7006946,1.6)
\psdots[linecolor=black, dotstyle=o, dotsize=0.2, fillcolor=white](6.9006944,1.6)
\psdots[linecolor=black, dotstyle=o, dotsize=0.2, fillcolor=white](6.1006947,1.6)
\psdots[linecolor=black, dotstyle=o, dotsize=0.2, fillcolor=white](5.3006945,1.6)
\psdots[linecolor=black, dotstyle=o, dotsize=0.2, fillcolor=white](4.9006944,1.2)
\psdots[linecolor=black, dotstyle=o, dotsize=0.2, fillcolor=white](5.7006946,1.2)
\psdots[linecolor=black, dotstyle=o, dotsize=0.2, fillcolor=white](6.5006943,1.2)
\psdots[linecolor=black, dotstyle=o, dotsize=0.2, fillcolor=white](7.3006945,1.2)
\psdots[linecolor=black, dotstyle=o, dotsize=0.2, fillcolor=white](8.100695,1.2)
\psdots[linecolor=black, dotstyle=o, dotsize=0.2, fillcolor=white](8.900695,1.2)
\psdots[linecolor=black, dotstyle=o, dotsize=0.2, fillcolor=white](9.700694,1.2)
\psdots[linecolor=black, dotstyle=o, dotsize=0.2, fillcolor=white](10.100695,0.8000001)
\psdots[linecolor=black, dotstyle=o, dotsize=0.2, fillcolor=white](9.300694,0.8000001)
\psdots[linecolor=black, dotstyle=o, dotsize=0.2, fillcolor=white](8.500694,0.8000001)
\psdots[linecolor=black, dotstyle=o, dotsize=0.2, fillcolor=white](7.7006946,0.8000001)
\psdots[linecolor=black, dotstyle=o, dotsize=0.2, fillcolor=white](6.9006944,0.8000001)
\psdots[linecolor=black, dotstyle=o, dotsize=0.2, fillcolor=white](6.1006947,0.8000001)
\psdots[linecolor=black, dotstyle=o, dotsize=0.2, fillcolor=white](5.3006945,0.8000001)
\psdots[linecolor=black, dotstyle=o, dotsize=0.2, fillcolor=white](4.5006943,0.8000001)
\psdots[linecolor=black, dotstyle=o, dotsize=0.2, fillcolor=white](4.1006947,0.40000007)
\psdots[linecolor=black, dotstyle=o, dotsize=0.2, fillcolor=white](4.9006944,0.40000007)
\psdots[linecolor=black, dotstyle=o, dotsize=0.2, fillcolor=white](3.7006943,0.0)
\psdots[linecolor=black, dotstyle=o, dotsize=0.2, fillcolor=white](3.3006945,-0.39999992)
\psdots[linecolor=black, dotstyle=o, dotsize=0.2, fillcolor=white](2.9006944,-0.79999995)
\psdots[linecolor=black, dotstyle=o, dotsize=0.2, fillcolor=white](3.7006943,-0.79999995)
\psdots[linecolor=black, dotstyle=o, dotsize=0.2, fillcolor=white](4.1006947,-0.39999992)
\psdots[linecolor=black, dotstyle=o, dotsize=0.2, fillcolor=white](4.5006943,0.0)
\psdots[linecolor=black, dotstyle=o, dotsize=0.2, fillcolor=white](4.5006943,-0.79999995)
\psdots[linecolor=black, dotstyle=o, dotsize=0.2, fillcolor=white](4.9006944,-0.39999992)
\psdots[linecolor=black, dotstyle=o, dotsize=0.2, fillcolor=white](5.3006945,0.0)
\psdots[linecolor=black, dotstyle=o, dotsize=0.2, fillcolor=white](5.7006946,0.40000007)
\psdots[linecolor=black, dotstyle=o, dotsize=0.2, fillcolor=white](5.3006945,-0.79999995)
\psdots[linecolor=black, dotstyle=o, dotsize=0.2, fillcolor=white](5.7006946,-0.39999992)
\psdots[linecolor=black, dotstyle=o, dotsize=0.2, fillcolor=white](6.1006947,0.0)
\psdots[linecolor=black, dotstyle=o, dotsize=0.2, fillcolor=white](6.5006943,0.40000007)
\psdots[linecolor=black, dotstyle=o, dotsize=0.2, fillcolor=white](6.9006944,0.0)
\psdots[linecolor=black, dotstyle=o, dotsize=0.2, fillcolor=white](7.3006945,0.40000007)
\psdots[linecolor=black, dotstyle=o, dotsize=0.2, fillcolor=white](6.5006943,-0.39999992)
\psdots[linecolor=black, dotstyle=o, dotsize=0.2, fillcolor=white](6.1006947,-0.79999995)
\psdots[linecolor=black, dotstyle=o, dotsize=0.2, fillcolor=white](8.100695,0.40000007)
\psdots[linecolor=black, dotstyle=o, dotsize=0.2, fillcolor=white](7.7006946,0.0)
\psdots[linecolor=black, dotstyle=o, dotsize=0.2, fillcolor=white](7.3006945,-0.39999992)
\psdots[linecolor=black, dotstyle=o, dotsize=0.2, fillcolor=white](6.9006944,-0.79999995)
\psdots[linecolor=black, dotstyle=o, dotsize=0.2, fillcolor=white](8.900695,0.40000007)
\psdots[linecolor=black, dotstyle=o, dotsize=0.2, fillcolor=white](8.500694,0.0)
\psdots[linecolor=black, dotstyle=o, dotsize=0.2, fillcolor=white](8.100695,-0.39999992)
\psdots[linecolor=black, dotstyle=o, dotsize=0.2, fillcolor=white](7.7006946,-0.79999995)
\psdots[linecolor=black, dotstyle=o, dotsize=0.2, fillcolor=white](8.500694,-0.79999995)
\psdots[linecolor=black, dotstyle=o, dotsize=0.2, fillcolor=white](8.900695,-0.39999992)
\psdots[linecolor=black, dotstyle=o, dotsize=0.2, fillcolor=white](9.300694,0.0)
\psdots[linecolor=black, dotstyle=o, dotsize=0.2, fillcolor=white](9.700694,0.40000007)
\psdots[linecolor=black, dotstyle=o, dotsize=0.2, fillcolor=white](10.500694,0.40000007)
\psdots[linecolor=black, dotstyle=o, dotsize=0.2, fillcolor=white](10.100695,0.0)
\psdots[linecolor=black, dotstyle=o, dotsize=0.2, fillcolor=white](9.700694,-0.39999992)
\psdots[linecolor=black, dotstyle=o, dotsize=0.2, fillcolor=white](9.300694,-0.79999995)
\psdots[linecolor=black, dotstyle=o, dotsize=0.2, fillcolor=white](10.100695,-0.79999995)
\psdots[linecolor=black, dotstyle=o, dotsize=0.2, fillcolor=white](10.500694,-0.39999992)
\psdots[linecolor=black, dotstyle=o, dotsize=0.2, fillcolor=white](10.900695,0.0)
\psdots[linecolor=black, dotstyle=o, dotsize=0.2, fillcolor=white](11.300694,-0.39999992)
\psdots[linecolor=black, dotstyle=o, dotsize=0.2, fillcolor=white](11.700694,-0.79999995)
\psdots[linecolor=black, dotstyle=o, dotsize=0.2, fillcolor=white](10.900695,-0.79999995)
\psdots[linecolor=black, dotstyle=o, dotsize=0.2, fillcolor=white](12.100695,-1.1999999)
\psdots[linecolor=black, dotstyle=o, dotsize=0.2, fillcolor=white](11.300694,-1.1999999)
\psdots[linecolor=black, dotstyle=o, dotsize=0.2, fillcolor=white](10.500694,-1.1999999)
\psdots[linecolor=black, dotstyle=o, dotsize=0.2, fillcolor=white](9.700694,-1.1999999)
\psdots[linecolor=black, dotstyle=o, dotsize=0.2, fillcolor=white](8.900695,-1.1999999)
\psdots[linecolor=black, dotstyle=o, dotsize=0.2, fillcolor=white](8.100695,-1.1999999)
\psdots[linecolor=black, dotstyle=o, dotsize=0.2, fillcolor=white](7.3006945,-1.1999999)
\psdots[linecolor=black, dotstyle=o, dotsize=0.2, fillcolor=white](6.5006943,-1.1999999)
\psdots[linecolor=black, dotstyle=o, dotsize=0.2, fillcolor=white](5.7006946,-1.1999999)
\psdots[linecolor=black, dotstyle=o, dotsize=0.2, fillcolor=white](4.9006944,-1.1999999)
\psdots[linecolor=black, dotstyle=o, dotsize=0.2, fillcolor=white](4.1006947,-1.1999999)
\psdots[linecolor=black, dotstyle=o, dotsize=0.2, fillcolor=white](3.3006945,-1.1999999)
\psdots[linecolor=black, dotstyle=o, dotsize=0.2, fillcolor=white](2.5006945,-1.1999999)
\psdots[linecolor=black, dotsize=0.2](7.3006945,3.6000001)
\psdots[linecolor=black, dotsize=0.2](6.9006944,3.2)
\psdots[linecolor=black, dotsize=0.2](7.7006946,3.2)
\psdots[linecolor=black, dotsize=0.2](7.3006945,2.8000002)
\psdots[linecolor=black, dotsize=0.2](8.100695,2.8000002)
\psdots[linecolor=black, dotsize=0.2](7.7006946,2.4)
\psdots[linecolor=black, dotsize=0.2](6.9006944,2.4)
\psdots[linecolor=black, dotsize=0.2](7.3006945,2.0)
\psdots[linecolor=black, dotsize=0.2](6.1006947,1.6)
\psdots[linecolor=black, dotsize=0.2](5.7006946,1.2)
\psdots[linecolor=black, dotsize=0.2](6.1006947,0.8000001)
\psdots[linecolor=black, dotsize=0.2](5.3006945,0.8000001)
\psdots[linecolor=black, dotsize=0.2](5.7006946,0.40000007)
\psdots[linecolor=black, dotsize=0.2](6.5006943,0.40000007)
\psdots[linecolor=black, dotsize=0.2](6.1006947,0.0)
\psdots[linecolor=black, dotsize=0.2](5.3006945,0.0)
\psdots[linecolor=black, dotsize=0.2](5.7006946,-0.39999992)
\psdots[linecolor=black, dotsize=0.2](5.3006945,-0.79999995)
\psdots[linecolor=black, dotsize=0.2](5.7006946,-1.1999999)
\psdots[linecolor=black, dotstyle=o, dotsize=0.2, fillcolor=white](2.1006944,-1.5999999)
\psdots[linecolor=black, dotstyle=o, dotsize=0.2, fillcolor=white](2.9006944,-1.5999999)
\psdots[linecolor=black, dotstyle=o, dotsize=0.2, fillcolor=white](3.7006943,-1.5999999)
\psdots[linecolor=black, dotstyle=o, dotsize=0.2, fillcolor=white](4.5006943,-1.5999999)
\psdots[linecolor=black, dotstyle=o, dotsize=0.2, fillcolor=white](5.3006945,-1.5999999)
\psdots[linecolor=black, dotstyle=o, dotsize=0.2, fillcolor=white](6.1006947,-1.5999999)
\psdots[linecolor=black, dotstyle=o, dotsize=0.2, fillcolor=white](6.9006944,-1.5999999)
\psdots[linecolor=black, dotstyle=o, dotsize=0.2, fillcolor=white](7.7006946,-1.5999999)
\psdots[linecolor=black, dotstyle=o, dotsize=0.2, fillcolor=white](8.500694,-1.5999999)
\psdots[linecolor=black, dotstyle=o, dotsize=0.2, fillcolor=white](9.300694,-1.5999999)
\psdots[linecolor=black, dotstyle=o, dotsize=0.2, fillcolor=white](10.100695,-1.5999999)
\psdots[linecolor=black, dotstyle=o, dotsize=0.2, fillcolor=white](10.900695,-1.5999999)
\psdots[linecolor=black, dotstyle=o, dotsize=0.2, fillcolor=white](11.700694,-1.5999999)
\psdots[linecolor=black, dotstyle=o, dotsize=0.2, fillcolor=white](12.500694,-1.5999999)
\psdots[linecolor=black, dotstyle=o, dotsize=0.2, fillcolor=white](12.900695,-1.9999999)
\psdots[linecolor=black, dotstyle=o, dotsize=0.2, fillcolor=white](12.100695,-1.9999999)
\psdots[linecolor=black, dotstyle=o, dotsize=0.2, fillcolor=white](11.300694,-1.9999999)
\psdots[linecolor=black, dotstyle=o, dotsize=0.2, fillcolor=white](10.500694,-1.9999999)
\psdots[linecolor=black, dotstyle=o, dotsize=0.2, fillcolor=white](9.700694,-1.9999999)
\psdots[linecolor=black, dotstyle=o, dotsize=0.2, fillcolor=white](8.900695,-1.9999999)
\psdots[linecolor=black, dotstyle=o, dotsize=0.2, fillcolor=white](8.100695,-1.9999999)
\psdots[linecolor=black, dotstyle=o, dotsize=0.2, fillcolor=white](7.3006945,-1.9999999)
\psdots[linecolor=black, dotstyle=o, dotsize=0.2, fillcolor=white](6.5006943,-1.9999999)
\psdots[linecolor=black, dotstyle=o, dotsize=0.2, fillcolor=white](5.7006946,-1.9999999)
\psdots[linecolor=black, dotstyle=o, dotsize=0.2, fillcolor=white](4.9006944,-1.9999999)
\psdots[linecolor=black, dotstyle=o, dotsize=0.2, fillcolor=white](4.1006947,-1.9999999)
\psdots[linecolor=black, dotstyle=o, dotsize=0.2, fillcolor=white](3.3006945,-1.9999999)
\psdots[linecolor=black, dotstyle=o, dotsize=0.2, fillcolor=white](2.5006945,-1.9999999)
\psdots[linecolor=black, dotstyle=o, dotsize=0.2, fillcolor=white](1.7006944,-1.9999999)
\psdots[linecolor=black, dotstyle=o, dotsize=0.2, fillcolor=white](1.3006945,-2.3999999)
\psdots[linecolor=black, dotstyle=o, dotsize=0.2, fillcolor=white](2.1006944,-2.3999999)
\psdots[linecolor=black, dotstyle=o, dotsize=0.2, fillcolor=white](0.90069443,-2.8)
\psdots[linecolor=black, dotstyle=o, dotsize=0.2, fillcolor=white](1.7006944,-2.8)
\psdots[linecolor=black, dotstyle=o, dotsize=0.2, fillcolor=white](2.5006945,-2.8)
\psdots[linecolor=black, dotstyle=o, dotsize=0.2, fillcolor=white](2.9006944,-2.3999999)
\psdots[linecolor=black, dotstyle=o, dotsize=0.2, fillcolor=white](3.3006945,-2.8)
\psdots[linecolor=black, dotstyle=o, dotsize=0.2, fillcolor=white](3.7006943,-2.3999999)
\psdots[linecolor=black, dotstyle=o, dotsize=0.2, fillcolor=white](4.1006947,-2.8)
\psdots[linecolor=black, dotstyle=o, dotsize=0.2, fillcolor=white](4.5006943,-2.3999999)
\psdots[linecolor=black, dotstyle=o, dotsize=0.2, fillcolor=white](4.9006944,-2.8)
\psdots[linecolor=black, dotstyle=o, dotsize=0.2, fillcolor=white](5.3006945,-2.3999999)
\psdots[linecolor=black, dotstyle=o, dotsize=0.2, fillcolor=white](5.7006946,-2.8)
\psdots[linecolor=black, dotstyle=o, dotsize=0.2, fillcolor=white](6.1006947,-2.3999999)
\psdots[linecolor=black, dotstyle=o, dotsize=0.2, fillcolor=white](6.5006943,-2.8)
\psdots[linecolor=black, dotstyle=o, dotsize=0.2, fillcolor=white](6.9006944,-2.3999999)
\psdots[linecolor=black, dotstyle=o, dotsize=0.2, fillcolor=white](7.3006945,-2.8)
\psdots[linecolor=black, dotstyle=o, dotsize=0.2, fillcolor=white](7.7006946,-2.3999999)
\psdots[linecolor=black, dotstyle=o, dotsize=0.2, fillcolor=white](8.100695,-2.8)
\psdots[linecolor=black, dotstyle=o, dotsize=0.2, fillcolor=white](8.500694,-2.3999999)
\psdots[linecolor=black, dotstyle=o, dotsize=0.2, fillcolor=white](8.900695,-2.8)
\psdots[linecolor=black, dotstyle=o, dotsize=0.2, fillcolor=white](9.300694,-2.3999999)
\psdots[linecolor=black, dotstyle=o, dotsize=0.2, fillcolor=white](9.700694,-2.8)
\psdots[linecolor=black, dotstyle=o, dotsize=0.2, fillcolor=white](10.100695,-2.3999999)
\psdots[linecolor=black, dotstyle=o, dotsize=0.2, fillcolor=white](10.500694,-2.8)
\psdots[linecolor=black, dotstyle=o, dotsize=0.2, fillcolor=white](10.900695,-2.3999999)
\psdots[linecolor=black, dotstyle=o, dotsize=0.2, fillcolor=white](11.300694,-2.8)
\psdots[linecolor=black, dotstyle=o, dotsize=0.2, fillcolor=white](11.700694,-2.3999999)
\psdots[linecolor=black, dotstyle=o, dotsize=0.2, fillcolor=white](12.100695,-2.8)
\psdots[linecolor=black, dotstyle=o, dotsize=0.2, fillcolor=white](12.500694,-2.3999999)
\psdots[linecolor=black, dotstyle=o, dotsize=0.2, fillcolor=white](12.900695,-2.8)
\psdots[linecolor=black, dotstyle=o, dotsize=0.2, fillcolor=white](13.300694,-2.3999999)
\psdots[linecolor=black, dotstyle=o, dotsize=0.2, fillcolor=white](13.700694,-2.8)
\psdots[linecolor=black, dotstyle=o, dotsize=0.2, fillcolor=white](14.100695,-3.1999998)
\psdots[linecolor=black, dotstyle=o, dotsize=0.2, fillcolor=white](14.500694,-3.6)
\psdots[linecolor=black, dotstyle=o, dotsize=0.2, fillcolor=white](13.700694,-3.6)
\psdots[linecolor=black, dotstyle=o, dotsize=0.2, fillcolor=white](13.300694,-3.1999998)
\psdots[linecolor=black, dotstyle=o, dotsize=0.2, fillcolor=white](12.900695,-3.6)
\psdots[linecolor=black, dotstyle=o, dotsize=0.2, fillcolor=white](12.500694,-3.1999998)
\psdots[linecolor=black, dotstyle=o, dotsize=0.2, fillcolor=white](12.100695,-3.6)
\psdots[linecolor=black, dotstyle=o, dotsize=0.2, fillcolor=white](11.700694,-3.1999998)
\psdots[linecolor=black, dotstyle=o, dotsize=0.2, fillcolor=white](11.300694,-3.6)
\psdots[linecolor=black, dotstyle=o, dotsize=0.2, fillcolor=white](10.900695,-3.1999998)
\psdots[linecolor=black, dotstyle=o, dotsize=0.2, fillcolor=white](10.500694,-3.6)
\psdots[linecolor=black, dotstyle=o, dotsize=0.2, fillcolor=white](10.100695,-3.1999998)
\psdots[linecolor=black, dotstyle=o, dotsize=0.2, fillcolor=white](9.700694,-3.6)
\psdots[linecolor=black, dotstyle=o, dotsize=0.2, fillcolor=white](9.300694,-3.1999998)
\psdots[linecolor=black, dotstyle=o, dotsize=0.2, fillcolor=white](8.900695,-3.6)
\psdots[linecolor=black, dotstyle=o, dotsize=0.2, fillcolor=white](8.500694,-3.1999998)
\psdots[linecolor=black, dotstyle=o, dotsize=0.2, fillcolor=white](8.100695,-3.6)
\psdots[linecolor=black, dotstyle=o, dotsize=0.2, fillcolor=white](7.7006946,-3.1999998)
\psdots[linecolor=black, dotstyle=o, dotsize=0.2, fillcolor=white](7.3006945,-3.6)
\psdots[linecolor=black, dotstyle=o, dotsize=0.2, fillcolor=white](6.9006944,-3.1999998)
\psdots[linecolor=black, dotstyle=o, dotsize=0.2, fillcolor=white](6.5006943,-3.6)
\psdots[linecolor=black, dotstyle=o, dotsize=0.2, fillcolor=white](6.1006947,-3.1999998)
\psdots[linecolor=black, dotstyle=o, dotsize=0.2, fillcolor=white](5.7006946,-3.6)
\psdots[linecolor=black, dotstyle=o, dotsize=0.2, fillcolor=white](5.3006945,-3.1999998)
\psdots[linecolor=black, dotstyle=o, dotsize=0.2, fillcolor=white](4.9006944,-3.6)
\psdots[linecolor=black, dotstyle=o, dotsize=0.2, fillcolor=white](4.5006943,-3.1999998)
\psdots[linecolor=black, dotstyle=o, dotsize=0.2, fillcolor=white](4.1006947,-3.6)
\psdots[linecolor=black, dotstyle=o, dotsize=0.2, fillcolor=white](3.7006943,-3.1999998)
\psdots[linecolor=black, dotstyle=o, dotsize=0.2, fillcolor=white](3.3006945,-3.6)
\psdots[linecolor=black, dotstyle=o, dotsize=0.2, fillcolor=white](2.9006944,-3.1999998)
\psdots[linecolor=black, dotstyle=o, dotsize=0.2, fillcolor=white](2.5006945,-3.6)
\psdots[linecolor=black, dotstyle=o, dotsize=0.2, fillcolor=white](2.1006944,-3.1999998)
\psdots[linecolor=black, dotstyle=o, dotsize=0.2, fillcolor=white](1.7006944,-3.6)
\psdots[linecolor=black, dotstyle=o, dotsize=0.2, fillcolor=white](1.3006945,-3.1999998)
\psdots[linecolor=black, dotstyle=o, dotsize=0.2, fillcolor=white](0.90069443,-3.6)
\psdots[linecolor=black, dotstyle=o, dotsize=0.2, fillcolor=white](0.50069445,-3.1999998)
\psdots[linecolor=black, dotstyle=o, dotsize=0.2, fillcolor=white](0.100694425,-3.6)
\psdots[linecolor=black, dotsize=0.2](12.500694,-2.3999999)
\psdots[linecolor=black, dotsize=0.2](12.100695,-2.8)
\psdots[linecolor=black, dotsize=0.2](12.500694,-3.1999998)
\psdots[linecolor=black, dotsize=0.2](11.700694,-3.1999998)
\psdots[linecolor=black, dotsize=0.2](12.100695,-3.6)
\end{pspicture}
}}

\def\treePagliacci{
\psscalebox{.33 .33} 
{
\begin{pspicture}[shift=-7.4](0,-7.4)(38.0,7.4)
\rput(15.4,4.7){\scalebox{1.5}{$2\, 1\, 3\, 2\, 4\, 3\, 5\, 4\, 6\, 5$}}
\rput(16.6,-7.8){\scalebox{1.5}{$2\, 1\, 3\, 2\, 4\, 3\, 5\, 4\, 6\, 5$}}
\rput(5.8,1.95){\scalebox{1.5}{$2\, 1\, 3\, 2\, 4\, 3\, 5\, 4\, 6\, 5$}}
\rput(24.25,1.95){\scalebox{1.5}{$2\, 1\, 3\, 2\, 4\, 3\, 5\, 4\, 6\, 5$}}
\rput(5.8,-1.25){\scalebox{1.5}{$2\, 1\, 3\, 2\, 4\, 3\, 5\, 4\, 6\, 5$}}
\rput(19,-1.25){\scalebox{1.5}{$2\, 1\, 3\, 2\, 4\, 3\, 5\, 4\, 6\, 5$}}
\rput(29.8,-1.25){\scalebox{1.5}{$2\, 1\, 3\, 2\, 4\, 3\, 5\, 4\, 6\, 5$}}
\rput(2.2,-4.5){\scalebox{1.5}{$2\, 1\, 3\, 2\, 4\, 3\, 5\, 4\, 6\, 5$}}
\rput(9.35,-4.5){\scalebox{1.5}{$2\, 1\, 3\, 2\, 4\, 3\, 5\, 4\, 6\, 5$}}
\rput(19,-4.5){\scalebox{1.5}{$2\, 1\, 3\, 2\, 4\, 3\, 5\, 4\, 6\, 5$}}
\rput(26.2,-4.5){\scalebox{1.5}{$2\, 1\, 3\, 2\, 4\, 3\, 5\, 4\, 6\, 5$}}
\rput(33.4,-4.5){\scalebox{1.5}{$2\, 1\, 3\, 2\, 4\, 3\, 5\, 4\, 6\, 5$}}
\rput(21.4,-7.8){\scalebox{1.5}{$2\, 1\, 3\, 2\, 4\, 3\, 5\, 4\, 6\, 5$}}
\rput(26.2,-7.8){\scalebox{1.5}{$2\, 1\, 3\, 2\, 4\, 3\, 5\, 4\, 6\, 5$}}
\rput(31,-7.8){\scalebox{1.5}{$2\, 1\, 3\, 2\, 4\, 3\, 5\, 4\, 6\, 5$}}
\rput(35.8,-7.8){\scalebox{1.5}{$2\, 1\, 3\, 2\, 4\, 3\, 5\, 4\, 6\, 5$}}
\rput(11.8,-7.8){\scalebox{1.5}{$2\, 1\, 3\, 2\, 4\, 3\, 5\, 4\, 6\, 5$}}
\rput(7,-7.8){\scalebox{1.5}{$2\, 1\, 3\, 2\, 4\, 3\, 5\, 4\, 6\, 5$}}
\rput(2.2,-7.8){\scalebox{1.5}{$2\, 1\, 3\, 2\, 4\, 3\, 5\, 4\, 6\, 5$}}

\psframe[linecolor=black, linewidth=0.04, linestyle=dotted, dotsep=0.10583334cm, dimen=outer](14.0,-5.4)(9.6,-7.4)
\psframe[linecolor=black, linewidth=0.04, linestyle=dotted, dotsep=0.10583334cm, dimen=outer](18.8,-5.4)(14.4,-7.4)
\psframe[linecolor=black, linewidth=0.04, linestyle=dotted, dotsep=0.10583334cm, dimen=outer](23.6,-5.4)(19.2,-7.4)
\psframe[linecolor=black, linewidth=0.04, linestyle=dotted, dotsep=0.10583334cm, dimen=outer](28.4,-5.4)(24.0,-7.4)
\psframe[linecolor=black, linewidth=0.04, linestyle=dotted, dotsep=0.10583334cm, dimen=outer](33.2,-5.4)(28.8,-7.4)
\psframe[linecolor=black, linewidth=0.04, linestyle=dotted, dotsep=0.10583334cm, dimen=outer](9.2,-5.4)(4.8,-7.4)
\psframe[linecolor=black, linewidth=0.04, linestyle=dotted, dotsep=0.10583334cm, dimen=outer](4.4,-5.4)(0.0,-7.4)
\psframe[linecolor=black, linewidth=0.04, linestyle=dotted, dotsep=0.10583334cm, dimen=outer](38.0,-5.4)(33.6,-7.4)
\psline[linecolor=black, linewidth=0.04, dotsize=0.07055555cm 2.0]{-*}(0.8,-7.4)(0.8,-6.6)
\psline[linecolor=black, linewidth=0.04, dotsize=0.07055555cm 2.0]{-*}(5.6,-7.4)(5.6,-6.6)
\psline[linecolor=black, linewidth=0.04, dotsize=0.07055555cm 2.0]{-*}(10.4,-7.4)(10.4,-6.6)
\psline[linecolor=black, linewidth=0.04, dotsize=0.07055555cm 2.0]{-*}(20.0,-7.4)(20.0,-6.6)
\psline[linecolor=black, linewidth=0.04, dotsize=0.07055555cm 2.0]{-*}(24.8,-7.4)(24.8,-6.6)
\psline[linecolor=black, linewidth=0.04, dotsize=0.07055555cm 2.0]{-*}(29.6,-7.4)(29.6,-6.6)
\psframe[linecolor=black, linewidth=0.04, linestyle=dotted, dotsep=0.10583334cm, dimen=outer](4.4,-2.2)(0.0,-4.2)
\psframe[linecolor=black, linewidth=0.04, linestyle=dotted, dotsep=0.10583334cm, dimen=outer](11.6,-2.2)(7.2,-4.2)
\psframe[linecolor=black, linewidth=0.04, linestyle=dotted, dotsep=0.10583334cm, dimen=outer](21.2,-2.2)(16.8,-4.2)
\psframe[linecolor=black, linewidth=0.04, linestyle=dotted, dotsep=0.10583334cm, dimen=outer](28.4,-2.2)(24.0,-4.2)
\psframe[linecolor=black, linewidth=0.04, linestyle=dotted, dotsep=0.10583334cm, dimen=outer](35.6,-2.2)(31.2,-4.2)
\psframe[linecolor=black, linewidth=0.04, linestyle=dotted, dotsep=0.10583334cm, dimen=outer](8.0,1.0)(3.6,-1.0)
\psframe[linecolor=black, linewidth=0.04, linestyle=dotted, dotsep=0.10583334cm, dimen=outer](21.2,1.0)(16.8,-1.0)
\psframe[linecolor=black, linewidth=0.04, linestyle=dotted, dotsep=0.10583334cm, dimen=outer](32.0,1.0)(27.6,-1.0)
\psframe[linecolor=black, linewidth=0.04, linestyle=dotted, dotsep=0.10583334cm, dimen=outer](8.0,4.2)(3.6,2.2)
\psframe[linecolor=black, linewidth=0.04, linestyle=dotted, dotsep=0.10583334cm, dimen=outer](26.4,4.2)(22.0,2.2)
\psframe[linecolor=black, linewidth=0.04, linestyle=dotted, dotsep=0.10583334cm, dimen=outer](17.6,7.4)(13.2,5.0)
\psline[linecolor=black, linewidth=0.04, arrowsize=0.05291667cm 2.0,arrowlength=1.4,arrowinset=0.0]{->}(13.2,6.2)(6.0,6.2)(6.0,4.2)
\psline[linecolor=black, linewidth=0.04, arrowsize=0.05291667cm 2.0,arrowlength=1.4,arrowinset=0.0]{->}(17.6,6.2)(24.0,6.2)(24.0,4.2)
\psline[linecolor=black, linewidth=0.04, arrowsize=0.05291667cm 2.0,arrowlength=1.4,arrowinset=0.0]{->}(22.0,3.0)(19.2,3.0)(19.2,1.0)
\psline[linecolor=black, linewidth=0.04, arrowsize=0.05291667cm 2.0,arrowlength=1.4,arrowinset=0.0]{->}(26.4,3.0)(29.6,3.0)(29.6,1.0)
\psline[linecolor=black, linewidth=0.04, arrowsize=0.05291667cm 2.0,arrowlength=1.4,arrowinset=0.0]{->}(19.2,-1.4)(19.2,-2.2)
\psline[linecolor=black, linewidth=0.04, arrowsize=0.05291667cm 2.0,arrowlength=1.4,arrowinset=0.0]{->}(27.6,-0.2)(26.4,-0.2)(26.4,-2.2)
\psline[linecolor=black, linewidth=0.04, arrowsize=0.05291667cm 2.0,arrowlength=1.4,arrowinset=0.0]{->}(32.0,-0.2)(33.2,-0.2)(33.2,-2.2)
\psline[linecolor=black, linewidth=0.04, arrowsize=0.05291667cm 2.0,arrowlength=1.4,arrowinset=0.0]{->}(8.0,-0.2)(9.2,-0.2)(9.2,-2.2)
\psline[linecolor=black, linewidth=0.04, arrowsize=0.05291667cm 2.0,arrowlength=1.4,arrowinset=0.0]{->}(3.6,-0.2)(2.4,-0.2)(2.4,-2.2)
\psline[linecolor=black, linewidth=0.04, arrowsize=0.05291667cm 2.0,arrowlength=1.4,arrowinset=0.0]{->}(2.4,-4.6)(2.4,-5.4)
\psline[linecolor=black, linewidth=0.04, arrowsize=0.05291667cm 2.0,arrowlength=1.4,arrowinset=0.0]{->}(7.2,-3.4)(6.4,-3.4)(6.4,-5.4)
\psline[linecolor=black, linewidth=0.04, arrowsize=0.05291667cm 2.0,arrowlength=1.4,arrowinset=0.0]{->}(11.6,-3.0)(12.4,-3.0)(12.4,-5.4)
\psline[linecolor=black, linewidth=0.04, arrowsize=0.05291667cm 2.0,arrowlength=1.4,arrowinset=0.0]{->}(16.8,-3.0)(16.0,-3.0)(16.0,-5.4)
\psline[linecolor=black, linewidth=0.04, arrowsize=0.05291667cm 2.0,arrowlength=1.4,arrowinset=0.0]{->}(21.2,-3.0)(22.0,-3.0)(22.0,-5.4)
\psline[linecolor=black, linewidth=0.04, arrowsize=0.05291667cm 2.0,arrowlength=1.4,arrowinset=0.0]{->}(26.4,-4.6)(26.4,-5.4)
\psline[linecolor=black, linewidth=0.04, arrowsize=0.05291667cm 2.0,arrowlength=1.4,arrowinset=0.0]{->}(31.2,-3.0)(30.4,-3.0)(30.4,-5.4)
\psline[linecolor=black, linewidth=0.04, arrowsize=0.05291667cm 2.0,arrowlength=1.4,arrowinset=0.0]{->}(35.6,-3.0)(36.4,-3.0)(36.4,-5.4)
\psbezier[linecolor=black, linewidth=0.04](4.0,2.2)(4.0,2.7333333)(4.0,2.4666667)(4.0,3.0)(4.0,3.5333333)(5.2,3.5333333)(5.2,3.0)(5.2,2.4666667)(5.2,2.7333333)(5.2,2.2)
\psline[linecolor=black, linewidth=0.04, dotsize=0.07055555cm 2.0]{-*}(17.6,-1.0)(17.6,-0.2)
\psline[linecolor=black, linewidth=0.04, dotsize=0.07055555cm 2.0]{-*}(17.2,-1.0)(17.2,-0.2)
\psline[linecolor=black, linewidth=0.04, dotsize=0.07055555cm 2.0]{-*}(18.4,-1.0)(18.4,-0.2)
\psbezier[linecolor=black, linewidth=0.04](18.0,-1.0)(18.0,-0.46666667)(18.0,-0.73333335)(18.0,-0.2)(18.0,0.33333334)(19.2,0.33333334)(19.2,-0.2)(19.2,-0.73333335)(19.2,-0.46666667)(19.2,-1.0)
\psline[linecolor=black, linewidth=0.04, dotsize=0.07055555cm 2.0]{-*}(1.2,-7.4)(1.2,-6.6)
\psbezier[linecolor=black, linewidth=0.04](0.4,-7.4)(0.4,-7.0444446)(0.4,-6.9555554)(0.4,-6.6)(0.4,-6.2444444)(1.6,-6.2444444)(1.6,-6.6)(1.6,-6.9555554)(1.6,-7.0444446)(1.6,-7.4)
\psbezier[linecolor=black, linewidth=0.04](2.0,-7.4)(2.0,-7.0444446)(2.0,-6.9555554)(2.0,-6.6)(2.0,-6.2444444)(3.2,-6.2444444)(3.2,-6.6)(3.2,-6.9555554)(3.2,-7.0444446)(3.2,-7.4)
\psbezier[linecolor=black, linewidth=0.04](5.2,-7.4)(5.2,-6.866667)(5.2,-7.133333)(5.2,-6.6)(5.2,-6.0666666)(6.4,-6.0666666)(6.4,-6.6)(6.4,-7.133333)(6.4,-6.866667)(6.4,-7.4)
\psbezier[linecolor=black, linewidth=0.04](7.6,-7.4)(7.6,-6.866667)(7.6,-7.133333)(7.6,-6.6)(7.6,-6.0666666)(8.8,-6.0666666)(8.8,-6.6)(8.8,-7.133333)(8.8,-6.866667)(8.8,-7.4)
\psbezier[linecolor=black, linewidth=0.04](10.0,-7.4)(10.0,-6.866667)(10.0,-7.133333)(10.0,-6.6)(10.0,-6.0666666)(11.2,-6.0666666)(11.2,-6.6)(11.2,-7.133333)(11.2,-6.866667)(11.2,-7.4)
\psbezier[linecolor=black, linewidth=0.04](19.6,-7.4)(19.6,-6.866667)(19.6,-7.133333)(19.6,-6.6)(19.6,-6.0666666)(20.8,-6.0666666)(20.8,-6.6)(20.8,-7.133333)(20.8,-6.866667)(20.8,-7.4)
\psbezier[linecolor=black, linewidth=0.04](26.0,-7.4)(26.0,-6.866667)(26.0,-7.133333)(26.0,-6.6)(26.0,-6.0666666)(27.2,-6.0666666)(27.2,-6.6)(27.2,-7.133333)(27.2,-6.866667)(27.2,-7.4)
\psbezier[linecolor=black, linewidth=0.04](31.6,-7.4)(31.6,-6.866667)(31.6,-7.133333)(31.6,-6.6)(31.6,-6.0666666)(32.8,-6.0666666)(32.8,-6.6)(32.8,-7.133333)(32.8,-6.866667)(32.8,-7.4)
\psline[linecolor=black, linewidth=0.04, dotsize=0.07055555cm 2.0]{-*}(2.4,-7.4)(2.4,-6.6)
\psline[linecolor=black, linewidth=0.04, dotsize=0.07055555cm 2.0]{-*}(2.8,-7.4)(2.8,-6.6)
\psline[linecolor=black, linewidth=0.04, dotsize=0.07055555cm 2.0]{-*}(4.0,-7.4)(4.0,-6.6)
\psline[linecolor=black, linewidth=0.04, dotsize=0.07055555cm 2.0]{-*}(6.0,-7.4)(6.0,-6.6)
\psline[linecolor=black, linewidth=0.04, dotsize=0.07055555cm 2.0]{-*}(10.8,-7.4)(10.8,-6.6)
\psline[linecolor=black, linewidth=0.04, dotsize=0.07055555cm 2.0]{-*}(11.6,-7.4)(11.6,-6.6)
\psline[linecolor=black, linewidth=0.04, dotsize=0.07055555cm 2.0]{-*}(12.0,-7.4)(12.0,-6.6)
\psline[linecolor=black, linewidth=0.04, dotsize=0.07055555cm 2.0]{-*}(12.4,-7.4)(12.4,-6.6)
\psline[linecolor=black, linewidth=0.04, dotsize=0.07055555cm 2.0]{-*}(12.8,-7.4)(12.8,-6.6)
\psline[linecolor=black, linewidth=0.04, dotsize=0.07055555cm 2.0]{-*}(13.6,-7.4)(13.6,-6.6)
\psline[linecolor=black, linewidth=0.04, dotsize=0.07055555cm 2.0]{-*}(20.4,-7.4)(20.4,-6.6)
\psline[linecolor=black, linewidth=0.04, dotsize=0.07055555cm 2.0]{-*}(21.2,-7.4)(21.2,-6.6)
\psline[linecolor=black, linewidth=0.04, dotsize=0.07055555cm 2.0]{-*}(21.6,-7.4)(21.6,-6.6)
\psline[linecolor=black, linewidth=0.04, dotsize=0.07055555cm 2.0]{-*}(22.0,-7.4)(22.0,-6.6)
\psline[linecolor=black, linewidth=0.04, dotsize=0.07055555cm 2.0]{-*}(22.4,-7.4)(22.4,-6.6)
\psline[linecolor=black, linewidth=0.04, dotsize=0.07055555cm 2.0]{-*}(23.2,-7.4)(23.2,-6.6)
\psline[linecolor=black, linewidth=0.04, dotsize=0.07055555cm 2.0]{-*}(24.4,-7.4)(24.4,-6.6)
\psline[linecolor=black, linewidth=0.04, dotsize=0.07055555cm 2.0]{-*}(25.2,-7.4)(25.2,-6.6)
\psline[linecolor=black, linewidth=0.04, dotsize=0.07055555cm 2.0]{-*}(25.6,-7.4)(25.6,-6.6)
\psline[linecolor=black, linewidth=0.04, dotsize=0.07055555cm 2.0]{-*}(26.4,-7.4)(26.4,-6.6)
\psline[linecolor=black, linewidth=0.04, dotsize=0.07055555cm 2.0]{-*}(26.8,-7.4)(26.8,-6.6)
\psline[linecolor=black, linewidth=0.04, dotsize=0.07055555cm 2.0]{-*}(28.0,-7.4)(28.0,-6.4)
\psline[linecolor=black, linewidth=0.04, dotsize=0.07055555cm 2.0]{-*}(29.2,-7.4)(29.2,-6.6)
\psline[linecolor=black, linewidth=0.04, dotsize=0.07055555cm 2.0]{-*}(30.0,-7.4)(30.0,-6.6)
\psline[linecolor=black, linewidth=0.04, dotsize=0.07055555cm 2.0]{-*}(30.4,-7.4)(30.4,-6.6)
\psline[linecolor=black, linewidth=0.04, dotsize=0.07055555cm 2.0]{-*}(30.8,-7.4)(30.8,-6.6)
\psline[linecolor=black, linewidth=0.04, dotsize=0.07055555cm 2.0]{-*}(31.2,-7.4)(31.2,-6.6)
\psline[linecolor=black, linewidth=0.04, dotsize=0.07055555cm 2.0]{-*}(32.0,-7.4)(32.0,-6.6)
\psline[linecolor=black, linewidth=0.04, dotsize=0.07055555cm 2.0]{-*}(0.8,-4.2)(0.8,-3.4)
\psline[linecolor=black, linewidth=0.04, dotsize=0.07055555cm 2.0]{-*}(1.2,-4.2)(1.2,-3.4)
\psbezier[linecolor=black, linewidth=0.04](0.4,-4.2)(0.4,-3.6666667)(0.4,-3.9333334)(0.4,-3.4)(0.4,-2.8666666)(1.6,-2.8666666)(1.6,-3.4)(1.6,-3.9333334)(1.6,-3.6666667)(1.6,-4.2)
\psbezier[linecolor=black, linewidth=0.04](2.0,-4.2)(2.0,-3.6666667)(2.0,-3.9333334)(2.0,-3.4)(2.0,-2.8666666)(3.2,-2.8666666)(3.2,-3.4)(3.2,-3.9333334)(3.2,-3.6666667)(3.2,-4.2)
\psline[linecolor=black, linewidth=0.04, dotsize=0.07055555cm 2.0]{-*}(2.4,-4.2)(2.4,-3.4)
\psline[linecolor=black, linewidth=0.04, dotsize=0.07055555cm 2.0]{-*}(3.6,-7.4)(3.6,-6.6)
\psline[linecolor=black, linewidth=0.04, dotsize=0.07055555cm 2.0]{-*}(13.2,-7.4)(13.2,-6.6)
\psline[linecolor=black, linewidth=0.04, dotsize=0.07055555cm 2.0]{-*}(22.8,-7.4)(22.8,-6.6)
\psline[linecolor=black, linewidth=0.04, dotsize=0.07055555cm 2.0]{-*}(27.6,-7.4)(27.6,-6.6)
\psline[linecolor=black, linewidth=0.04, dotsize=0.07055555cm 2.0]{-*}(32.4,-7.4)(32.4,-6.6)
\psline[linecolor=black, linewidth=0.04, dotsize=0.07055555cm 2.0]{-*}(3.6,-4.2)(3.6,-3.4)
\psline[linecolor=black, linewidth=0.04, dotsize=0.07055555cm 2.0]{-*}(8.4,-4.2)(8.4,-3.4)
\psbezier[linecolor=black, linewidth=0.04](7.6,-4.2)(7.6,-3.6666667)(7.6,-3.9333334)(7.6,-3.4)(7.6,-2.8666666)(8.8,-2.8666666)(8.8,-3.4)(8.8,-3.9333334)(8.8,-3.6666667)(8.8,-4.2)
\psline[linecolor=black, linewidth=0.04, dotsize=0.07055555cm 2.0]{-*}(9.6,-4.2)(9.6,-3.4)
\psline[linecolor=black, linewidth=0.04, dotsize=0.07055555cm 2.0]{-*}(10.8,-4.2)(10.8,-3.4)
\psline[linecolor=black, linewidth=0.04, dotsize=0.07055555cm 2.0]{-*}(8.0,-4.2)(8.0,-3.4)
\psline[linecolor=black, linewidth=0.04, dotsize=0.07055555cm 2.0]{-*}(14.0,5.0)(14.0,5.8)
\psline[linecolor=black, linewidth=0.04, dotsize=0.07055555cm 2.0]{-*}(16.8,5.0)(16.8,5.8)
\psline[linecolor=black, linewidth=0.04, dotsize=0.07055555cm 2.0]{-*}(4.4,2.2)(4.4,3.0)
\psline[linecolor=black, linewidth=0.04, dotsize=0.07055555cm 2.0]{-*}(7.2,2.2)(7.2,3.0)
\psbezier[linecolor=black, linewidth=0.04](4.0,-1.0)(4.0,-0.46666667)(4.0,-0.73333335)(4.0,-0.2)(4.0,0.33333334)(5.2,0.33333334)(5.2,-0.2)(5.2,-0.73333335)(5.2,-0.46666667)(5.2,-1.0)
\psline[linecolor=black, linewidth=0.04, dotsize=0.07055555cm 2.0]{-*}(4.4,-1.0)(4.4,-0.2)
\psline[linecolor=black, linewidth=0.04, dotsize=0.07055555cm 2.0]{-*}(7.2,-1.0)(7.2,-0.2)
\psline[linecolor=black, linewidth=0.04, dotsize=0.07055555cm 2.0]{-*}(4.8,-1.0)(4.8,-0.2)
\psline[linecolor=black, linewidth=0.04, dotsize=0.07055555cm 2.0]{-*}(6.0,-1.0)(6.0,-0.2)
\psline[linecolor=black, linewidth=0.04, dotsize=0.07055555cm 2.0]{-*}(9.2,-4.2)(9.2,-3.4)
\psline[linecolor=black, linewidth=0.04, dotsize=0.07055555cm 2.0]{-*}(10.4,-4.2)(10.4,-3.4)
\psline[linecolor=black, linewidth=0.04, dotsize=0.07055555cm 2.0]{-*}(6.8,-7.4)(6.8,-6.6)
\psline[linecolor=black, linewidth=0.04, dotsize=0.07055555cm 2.0]{-*}(7.2,-7.4)(7.2,-6.6)
\psline[linecolor=black, linewidth=0.04, dotsize=0.07055555cm 2.0]{-*}(8.0,-7.4)(8.0,-6.6)
\psline[linecolor=black, linewidth=0.04, dotsize=0.07055555cm 2.0]{-*}(8.4,-7.4)(8.4,-6.6)
\psline[linecolor=black, linewidth=0.04, dotsize=0.07055555cm 2.0]{-*}(22.8,2.2)(22.8,3.0)
\psline[linecolor=black, linewidth=0.04, dotsize=0.07055555cm 2.0]{-*}(25.6,2.2)(25.6,3.0)
\psline[linecolor=black, linewidth=0.04, dotsize=0.07055555cm 2.0]{-*}(22.4,2.2)(22.4,3.0)
\psline[linecolor=black, linewidth=0.04, dotsize=0.07055555cm 2.0]{-*}(23.6,2.2)(23.6,3.0)
\psline[linecolor=black, linewidth=0.04, arrowsize=0.05291667cm 2.0,arrowlength=1.4,arrowinset=0.0]{->}(5.6,1.8)(5.6,1.0)
\psline[linecolor=black, linewidth=0.04, dotsize=0.07055555cm 2.0]{-*}(20.4,-1.0)(20.4,-0.2)
\psline[linecolor=black, linewidth=0.04, dotsize=0.07055555cm 2.0]{-*}(28.4,-1.0)(28.4,-0.2)
\psline[linecolor=black, linewidth=0.04, dotsize=0.07055555cm 2.0]{-*}(31.2,-1.0)(31.2,-0.2)
\psline[linecolor=black, linewidth=0.04, dotsize=0.07055555cm 2.0]{-*}(28.0,-1.0)(28.0,-0.2)
\psline[linecolor=black, linewidth=0.04, dotsize=0.07055555cm 2.0]{-*}(29.2,-1.0)(29.2,-0.2)
\psline[linecolor=black, linewidth=0.04, dotsize=0.07055555cm 2.0]{-*}(28.8,-1.0)(28.8,-0.2)
\psline[linecolor=black, linewidth=0.04, dotsize=0.07055555cm 2.0]{-*}(30.0,-1.0)(30.0,-0.2)
\psline[linecolor=black, linewidth=0.04, dotsize=0.07055555cm 2.0]{-*}(17.6,-4.2)(17.6,-3.4)
\psline[linecolor=black, linewidth=0.04, dotsize=0.07055555cm 2.0]{-*}(17.2,-4.2)(17.2,-3.4)
\psline[linecolor=black, linewidth=0.04, dotsize=0.07055555cm 2.0]{-*}(18.4,-4.2)(18.4,-3.4)
\psbezier[linecolor=black, linewidth=0.04](18.0,-4.2)(18.0,-3.6666667)(18.0,-3.9333334)(18.0,-3.4)(18.0,-2.8666666)(19.2,-2.8666666)(19.2,-3.4)(19.2,-3.9333334)(19.2,-3.6666667)(19.2,-4.2)
\psline[linecolor=black, linewidth=0.04, dotsize=0.07055555cm 2.0]{-*}(20.4,-4.2)(20.4,-3.4)
\psline[linecolor=black, linewidth=0.04, dotsize=0.07055555cm 2.0]{-*}(18.8,-4.2)(18.8,-3.4)
\psline[linecolor=black, linewidth=0.04, dotsize=0.07055555cm 2.0]{-*}(20.0,-4.2)(20.0,-3.4)
\psline[linecolor=black, linewidth=0.04, dotsize=0.07055555cm 2.0]{-*}(24.8,-4.2)(24.8,-3.4)
\psline[linecolor=black, linewidth=0.04, dotsize=0.07055555cm 2.0]{-*}(27.6,-4.2)(27.6,-3.4)
\psline[linecolor=black, linewidth=0.04, dotsize=0.07055555cm 2.0]{-*}(24.4,-4.2)(24.4,-3.4)
\psline[linecolor=black, linewidth=0.04, dotsize=0.07055555cm 2.0]{-*}(25.6,-4.2)(25.6,-3.4)
\psline[linecolor=black, linewidth=0.04, dotsize=0.07055555cm 2.0]{-*}(25.2,-4.2)(25.2,-3.4)
\psline[linecolor=black, linewidth=0.04, dotsize=0.07055555cm 2.0]{-*}(26.4,-4.2)(26.4,-3.4)
\psline[linecolor=black, linewidth=0.04, dotsize=0.07055555cm 2.0]{-*}(32.0,-4.2)(32.0,-3.4)
\psline[linecolor=black, linewidth=0.04, dotsize=0.07055555cm 2.0]{-*}(34.8,-4.2)(34.8,-3.4)
\psline[linecolor=black, linewidth=0.04, dotsize=0.07055555cm 2.0]{-*}(31.6,-4.2)(31.6,-3.4)
\psline[linecolor=black, linewidth=0.04, dotsize=0.07055555cm 2.0]{-*}(32.8,-4.2)(32.8,-3.4)
\psline[linecolor=black, linewidth=0.04, dotsize=0.07055555cm 2.0]{-*}(32.4,-4.2)(32.4,-3.4)
\psline[linecolor=black, linewidth=0.04, dotsize=0.07055555cm 2.0]{-*}(33.6,-4.2)(33.6,-3.4)
\psline[linecolor=black, linewidth=0.04, dotsize=0.07055555cm 2.0]{-*}(33.2,-4.2)(33.2,-3.4)
\psline[linecolor=black, linewidth=0.04, dotsize=0.07055555cm 2.0]{-*}(34.4,-4.2)(34.4,-3.4)
\psbezier[linecolor=black, linewidth=0.04](26.0,-4.2)(26.0,-3.6666667)(26.0,-3.9333334)(26.0,-3.4)(26.0,-2.8666666)(27.2,-2.8666666)(27.2,-3.4)(27.2,-3.9333334)(27.2,-3.6666667)(27.2,-4.2)
\psline[linecolor=black, linewidth=0.04, dotsize=0.07055555cm 2.0]{-*}(15.2,-7.4)(15.2,-6.6)
\psbezier[linecolor=black, linewidth=0.04](15.6,-7.4)(15.6,-6.866667)(15.6,-7.133333)(15.6,-6.6)(15.6,-6.0666666)(16.8,-6.0666666)(16.8,-6.6)(16.8,-7.133333)(16.8,-6.866667)(16.8,-7.4)
\psbezier[linecolor=black, linewidth=0.04](17.2,-7.4)(17.2,-6.866667)(17.2,-7.133333)(17.2,-6.6)(17.2,-6.0666666)(18.4,-6.0666666)(18.4,-6.6)(18.4,-7.133333)(18.4,-6.866667)(18.4,-7.4)
\psline[linecolor=black, linewidth=0.04, dotsize=0.07055555cm 2.0]{-*}(16.0,-7.4)(16.0,-6.6)
\psline[linecolor=black, linewidth=0.04, dotsize=0.07055555cm 2.0]{-*}(16.4,-7.4)(16.4,-6.6)
\psline[linecolor=black, linewidth=0.04, dotsize=0.07055555cm 2.0]{-*}(17.6,-7.4)(17.6,-6.6)
\psline[linecolor=black, linewidth=0.04, dotsize=0.07055555cm 2.0]{-*}(18.0,-7.4)(18.0,-6.6)
\psline[linecolor=black, linewidth=0.04, dotsize=0.07055555cm 2.0]{-*}(14.8,-7.4)(14.8,-6.6)
\psline[linecolor=black, linewidth=0.04, dotsize=0.07055555cm 2.0]{-*}(34.0,-7.4)(34.0,-6.6)
\psline[linecolor=black, linewidth=0.04, dotsize=0.07055555cm 2.0]{-*}(34.4,-7.4)(34.4,-6.6)
\psline[linecolor=black, linewidth=0.04, dotsize=0.07055555cm 2.0]{-*}(34.8,-7.4)(34.8,-6.6)
\psline[linecolor=black, linewidth=0.04, dotsize=0.07055555cm 2.0]{-*}(35.2,-7.4)(35.2,-6.6)
\psline[linecolor=black, linewidth=0.04, dotsize=0.07055555cm 2.0]{-*}(35.6,-7.4)(35.6,-6.6)
\psline[linecolor=black, linewidth=0.04, dotsize=0.07055555cm 2.0]{-*}(36.0,-7.4)(36.0,-6.6)
\psline[linecolor=black, linewidth=0.04, dotsize=0.07055555cm 2.0]{-*}(36.4,-7.4)(36.4,-6.6)
\psline[linecolor=black, linewidth=0.04, dotsize=0.07055555cm 2.0]{-*}(36.8,-7.4)(36.8,-6.6)
\psline[linecolor=black, linewidth=0.04, dotsize=0.07055555cm 2.0]{-*}(37.2,-7.4)(37.2,-6.6)
\psline[linecolor=black, linewidth=0.04, dotsize=0.07055555cm 2.0]{-*}(37.6,-7.4)(37.6,-6.6)
\end{pspicture}
}
}

\def\sixvalent{ 
\psscalebox{1 1} 
{
\begin{pspicture}[shift=-1.6](0,-1.6019659)(7.2,1.6019659)
\definecolor{colour0}{rgb}{0.6,0.6,0.6}
\definecolor{colour1}{rgb}{0.0,0.2,0.8}
\definecolor{colour2}{rgb}{1.0,0.2,0.2}
\psframe[linecolor=black, linewidth=0.04, linestyle=dotted, dotsep=0.10583334cm, dimen=outer](7.2,1.598034)(0.0,-1.6019659)
\psbezier[linecolor=black, linewidth=0.04](2.8,-1.6019659)(2.8,-0.80196595)(2.8,-0.0019659423)(3.6,-0.0019659423828125)(4.4,-0.0019659423)(5.2,-0.40196595)(5.2,-1.6019659)
\psbezier[linecolor=colour0, linewidth=0.04](3.6,-0.0019659423)(4.5142856,0.39803407)(4.4362936,0.89417297)(5.428571,0.3980340576171875)(6.4208493,-0.09810488)(6.8,-0.0019659423)(6.8,-1.6019659)
\psbezier[linecolor=colour0, linewidth=0.04](4.4,-1.6019659)(4.4,-0.80196595)(4.0,-0.40196595)(3.6,-0.0019659423828125)(3.2,0.39803407)(2.4,-0.40196595)(2.8,1.598034)
\psbezier[linecolor=black, linewidth=0.04](3.6,-0.0019659423)(3.6,0.7980341)(4.0,0.7980341)(4.0,1.5980340576171874)
\psline[linecolor=colour1, linewidth=0.04, dotsize=0.07055555cm 2.0]{-*}(3.6,-1.6019659)(3.6,-0.80196595)
\psline[linecolor=colour2, linewidth=0.04, dotsize=0.07055555cm 2.0]{-*}(6.0,-1.6019659)(6.0,-0.80196595)
\rput(1.5,0){\scalebox{2}{$\cdots $} }
\rput(1.5,-2){\scalebox{2}{$\cdots $} }
\rput(2.7,-2){\scalebox{1}{$s_i$}}
\rput(3.6,-2){\scalebox{1}{$s_{i-1}$}}
\rput(4.5,-2){\scalebox{1}{$s_{i+1}$}}
\rput(5.3,-2){\scalebox{1}{$s_{i}$}}
\rput(6.0,-2){\scalebox{1}{$s_{i+2}$}}
\rput(6.9,-2){\scalebox{1}{$s_{i+1}$}}
\end{pspicture}
}
 }

\def\simplest{
\psscalebox{.5 .5} 
{
\begin{pspicture}[shift=-4.2](0,-4.207071)(17.2,4.207071)
\definecolor{colour0}{rgb}{1.0,0.2,0.2}
\rput(7.2,3.5){\scalebox{2}{$\emptyset$}}
\psframe[linecolor=black, linewidth=0.04, linestyle=dotted, dotsep=0.10583334cm, dimen=outer](1.6,-2.5929291)(0.0,-4.1929293)
\psframe[linecolor=black, linewidth=0.04, linestyle=dotted, dotsep=0.10583334cm, dimen=outer](5.2,-2.5929291)(3.6,-4.1929293)
\psframe[linecolor=black, linewidth=0.04, linestyle=dotted, dotsep=0.10583334cm, dimen=outer](10.0,-2.5929291)(8.4,-4.1929293)
\psframe[linecolor=black, linewidth=0.04, linestyle=dotted, dotsep=0.10583334cm, dimen=outer](13.6,-2.5929291)(12.0,-4.1929293)
\psframe[linecolor=black, linewidth=0.04, linestyle=dotted, dotsep=0.10583334cm, dimen=outer](17.2,-2.5929291)(15.6,-4.1929293)
\psframe[linecolor=black, linewidth=0.04, linestyle=dotted, dotsep=0.10583334cm, dimen=outer](10.0,-0.19292907)(8.4,-1.792929)
\psframe[linecolor=black, linewidth=0.04, linestyle=dotted, dotsep=0.10583334cm, dimen=outer](15.2,-0.19292907)(13.6,-1.792929)
\psframe[linecolor=black, linewidth=0.04, linestyle=dotted, dotsep=0.10583334cm, dimen=outer](3.6,-0.19292907)(2.0,-1.792929)
\psframe[linecolor=black, linewidth=0.04, linestyle=dotted, dotsep=0.10583334cm, dimen=outer](3.6,2.2070708)(2.0,0.6070709)
\psframe[linecolor=black, linewidth=0.04, linestyle=dotted, dotsep=0.10583334cm, dimen=outer](12.8,2.2070708)(11.2,0.6070709)
\psframe[linecolor=black, linewidth=0.04, linestyle=dotted, dotsep=0.10583334cm, dimen=outer](8.0,4.207071)(6.4,2.607071)
\psline[linecolor=black, linewidth=0.04, arrowsize=0.05291667cm 2.0,arrowlength=1.4,arrowinset=0.0]{->}(6.4,3.4070709)(2.8,3.4070709)(2.8,2.2070708)
\psline[linecolor=black, linewidth=0.04, arrowsize=0.05291667cm 2.0,arrowlength=1.4,arrowinset=0.0]{->}(8.0,3.4070709)(12.0,3.4070709)(12.0,2.2070708)
\psline[linecolor=black, linewidth=0.04, arrowsize=0.05291667cm 2.0,arrowlength=1.4,arrowinset=0.0]{->}(2.8,0.6070709)(2.8,-0.19292907)
\psline[linecolor=black, linewidth=0.04, arrowsize=0.05291667cm 2.0,arrowlength=1.4,arrowinset=0.0]{->}(11.2,1.4070709)(9.2,1.4070709)(9.2,-0.19292907)
\psline[linecolor=black, linewidth=0.04, arrowsize=0.05291667cm 2.0,arrowlength=1.4,arrowinset=0.0]{->}(12.8,1.4070709)(14.4,1.4070709)(14.4,-0.19292907)
\psline[linecolor=black, linewidth=0.04, arrowsize=0.05291667cm 2.0,arrowlength=1.4,arrowinset=0.0]{->}(2.0,-0.9929291)(0.8,-0.9929291)(0.8,-2.5929291)
\psline[linecolor=black, linewidth=0.04, arrowsize=0.05291667cm 2.0,arrowlength=1.4,arrowinset=0.0]{->}(3.6,-0.9929291)(4.4,-0.9929291)(4.4,-2.5929291)
\psline[linecolor=black, linewidth=0.04, arrowsize=0.05291667cm 2.0,arrowlength=1.4,arrowinset=0.0]{->}(9.2,-1.792929)(9.2,-2.5929291)
\psline[linecolor=black, linewidth=0.04, arrowsize=0.05291667cm 2.0,arrowlength=1.4,arrowinset=0.0]{->}(13.6,-0.9929291)(12.8,-0.9929291)(12.8,-2.5929291)
\psline[linecolor=black, linewidth=0.04, arrowsize=0.05291667cm 2.0,arrowlength=1.4,arrowinset=0.0]{->}(15.2,-0.9929291)(16.4,-0.9929291)(16.4,-2.5929291)
\psline[linecolor=colour0, linewidth=0.04](2.4,2.2070708)(2.4,0.6070709)
\psline[linecolor=colour0, linewidth=0.04](9.2,-0.19292907)(9.2,-1.792929)
\psline[linecolor=colour0, linewidth=0.04](1.2,-2.5929291)(1.2,-4.1929293)
\psline[linecolor=colour0, linewidth=0.04, dotsize=0.07055555cm 2.0]{-*}(11.6,0.6070709)(11.6,1.4070709)
\psline[linecolor=colour0, linewidth=0.04, dotsize=0.07055555cm 2.0]{-*}(14.0,-1.792929)(14.0,-0.9929291)
\psline[linecolor=colour0, linewidth=0.04, dotsize=0.07055555cm 2.0]{-*}(14.4,-1.792929)(14.4,-0.9929291)
\psline[linecolor=colour0, linewidth=0.04, dotsize=0.07055555cm 2.0]{-*}(8.8,-1.792929)(8.8,-0.9929291)
\psline[linecolor=colour0, linewidth=0.04, dotsize=0.07055555cm 2.0]{-*}(4.8,-4.1929293)(4.8,-3.392929)
\psline[linecolor=colour0, linewidth=0.04, dotsize=0.07055555cm 2.0]{-*}(8.8,-4.1929293)(8.8,-3.392929)
\psline[linecolor=colour0, linewidth=0.04, dotsize=0.07055555cm 2.0]{-*}(12.4,-4.1929293)(12.4,-3.392929)
\psline[linecolor=colour0, linewidth=0.04, dotsize=0.07055555cm 2.0]{-*}(12.8,-4.1929293)(12.8,-3.392929)
\psline[linecolor=colour0, linewidth=0.04, dotsize=0.07055555cm 2.0]{-*}(16.0,-4.1929293)(16.0,-3.392929)(16.0,-3.392929)(16.0,-3.392929)
\psline[linecolor=colour0, linewidth=0.04, dotsize=0.07055555cm 2.0]{-*}(16.4,-4.1929293)(16.4,-3.392929)
\psline[linecolor=colour0, linewidth=0.04, dotsize=0.07055555cm 2.0]{-*}(16.8,-4.1929293)(16.8,-3.392929)
\psline[linecolor=colour0, linewidth=0.04](13.2,-2.5929291)(13.2,-4.1929293)
\psbezier[linecolor=colour0, linewidth=0.04](2.351515,-1.7808079)(2.4,-0.9929291)(2.8,-0.59292907)(2.8363636,-1.780807865027226)
\psbezier[linecolor=colour0, linewidth=0.04](0.4,-4.1929293)(0.4,-2.992929)(0.8,-3.392929)(0.8,-4.192929077148437)
\psbezier[linecolor=colour0, linewidth=0.04](4.0,-4.1929293)(3.6,-3.7929292)(4.4,-2.5929291)(4.4,-4.192929077148437)
\psbezier[linecolor=colour0, linewidth=0.04](9.2,-4.1929293)(8.8,-2.992929)(10.0,-3.392929)(9.6,-4.192929077148437)
\end{pspicture}
}}

\def\fibotreeone{
\psscalebox{.5 .5} 
{
\begin{pspicture}[shift=-1.09](0,-1.0985577)(4.1971154,1.0985577)
\psdots[linecolor=black, dotsize=0.2](2.0985577,1.0)
\psdots[linecolor=black, dotsize=0.2](0.09855774,-1.0)
\psdots[linecolor=black, dotsize=0.2](4.098558,-1.0)
\psline[linecolor=black, linewidth=0.04](2.0985577,1.0)(0.09855774,-1.0)
\psline[linecolor=black, linewidth=0.04](2.0985577,1.0)(4.098558,-1.0)
\end{pspicture}
}
}

\def\CaseDtwo{
\psscalebox{.7 .7} 
{
\begin{pspicture}[shift=-1.2](0,-1.2)(10.4,1.2)
\psframe[linecolor=black, linewidth=0.04, linestyle=dotted, dotsep=0.10583334cm, dimen=outer](4.0,1.2)(0.0,-1.2)
\psline[linecolor=black, linewidth=0.04, dotsize=0.07055555cm 2.0]{-*}(0.8,-1.2)(0.8,0.0)
\psline[linecolor=black, linewidth=0.04, dotsize=0.07055555cm 2.0]{-*}(3.2,-1.2)(3.2,0.0)
\psframe[linecolor=black, linewidth=0.04, linestyle=dotted, dotsep=0.10583334cm, dimen=outer](10.4,1.2)(6.4,-1.2)
\psbezier[linecolor=black, linewidth=0.04](7.2,-1.2)(7.2,-0.8)(6.8,-0.4)(7.6,0.0)(8.4,0.4)(8.342507,0.51449573)(9.2,0.0)(10.057493,-0.51449573)(9.6,-0.8)(9.6,-1.2)
\rput(.8,-1.6){\scalebox{1.5}{$s_i$}}
\rput(3.3,-1.6){\scalebox{1.5}{$s_i$}}
\rput(7.2,-1.6){\scalebox{1.5}{$s_i$}}
\rput(9.7,-1.6){\scalebox{1.5}{$s_i$}}
\end{pspicture}
}
}

\def\CaseDone{ 
\psscalebox{.7 .7} 
{
\begin{pspicture}[shift=-1.2](0,-1.2)(4.0,1.2)
\definecolor{colour3}{rgb}{0.6,0.6,0.6}
\psframe[linecolor=black, linewidth=0.04, linestyle=dotted, dotsep=0.10583334cm, dimen=outer](4.0,1.2)(0.0,-1.2)
\psline[linecolor=black, linewidth=0.04, dotsize=0.07055555cm 2.0]{-*}(0.8,-1.2)(0.8,0.0)
\psline[linecolor=colour3, linewidth=0.04](2.0,-1.2)(2.0,1.2)
\psline[linecolor=black, linewidth=0.04, dotsize=0.07055555cm 2.0]{-*}(3.2,-1.2)(3.2,0.0)
\rput(.8,-1.6){\scalebox{1.5}{$s_i$}}
\rput(2,-1.6){\scalebox{1.5}{$s_{i+1}$}}
\rput(3.3,-1.6){\scalebox{1.5}{$s_i$}}
\end{pspicture}
}
}

\def\CaseCthree{
\psscalebox{.7 .7} 
{
\begin{pspicture}[shift=-1.2](0,-1.2)(4.0,1.2)
\psframe[linecolor=black, linewidth=0.04, linestyle=dotted, dotsep=0.10583334cm, dimen=outer](4.0,1.2)(0.0,-1.2)
\psline[linecolor=black, linewidth=0.04, dotsize=0.07055555cm 2.0]{-*}(0.8,-1.2)(0.8,0.0)
\psline[linecolor=black, linewidth=0.04](3.2,-1.2)(3.2,1.2)
\rput(.8,-1.6){\scalebox{1.5}{$s_i$}}
\rput(3.3,-1.6){\scalebox{1.5}{$s_i$}}
\end{pspicture}
}
}

\def\CaseCtwo{ 
\psscalebox{.7 .7} 
{
\begin{pspicture}[shift=-1.2](0,-1.2)(4.0,1.2)
\definecolor{colour3}{rgb}{0.6,0.6,0.6}
\psframe[linecolor=black, linewidth=0.04, linestyle=dotted, dotsep=0.10583334cm, dimen=outer](4.0,1.2)(0.0,-1.2)
\psline[linecolor=colour3, linewidth=0.04](2.0,1.2)(2.0,-1.2)
\psline[linecolor=black, linewidth=0.04, dotsize=0.07055555cm 2.0]{-*}(0.8,-1.2)(0.8,0.0)
\psline[linecolor=black, linewidth=0.04](3.2,-1.2)(3.2,1.2)
\rput(.8,-1.6){\scalebox{1.5}{$s_i$}}
\rput(2,-1.6){\scalebox{1.5}{$s_{i+1}$}}
\rput(3.3,-1.6){\scalebox{1.5}{$s_i$}}
\end{pspicture}
}
 }

\def\CaseCone{
\psscalebox{.7 .7} 
{
\begin{pspicture}[shift=-1.2](0,-1.2)(4.0,1.2)
\definecolor{colour3}{rgb}{0.6,0.6,0.6}
\psframe[linecolor=black, linewidth=0.04, linestyle=dotted, dotsep=0.10583334cm, dimen=outer](4.0,1.2)(0.0,-1.2)
\psline[linecolor=black, linewidth=0.04](0.8,1.2)(0.8,-1.2)
\psline[linecolor=colour3, linewidth=0.04](2.0,1.2)(2.0,-1.2)
\psline[linecolor=black, linewidth=0.04, dotsize=0.07055555cm 2.0]{-*}(3.2,-1.2)(3.2,0.0)
\rput(.8,-1.6){\scalebox{1.5}{$s_i$}}
\rput(2,-1.6){\scalebox{1.5}{$s_{i+1}$}}
\rput(3.3,-1.6){\scalebox{1.5}{$s_i$}}
\end{pspicture}
}
}

\def\CaseBone{
\psscalebox{.7 .7} 
{
\begin{pspicture}[shift=-1.2](0,-1.2)(4.0,1.2)
\definecolor{colour3}{rgb}{0.6,0.6,0.6}
\psframe[linecolor=black, linewidth=0.04, linestyle=dotted, dotsep=0.10583334cm, dimen=outer](4.0,1.2)(0.0,-1.2)
\psline[linecolor=black, linewidth=0.04](0.8,1.2)(0.8,-1.2)
\psline[linecolor=black, linewidth=0.04](3.2,1.2)(3.2,-1.2)
\psline[linecolor=colour3, linewidth=0.04](2.0,1.2)(2.0,-1.2)
\rput(.8,-1.6){\scalebox{1.5}{$s_i$}}
\rput(2,-1.6){\scalebox{1.5}{$s_{i+1}$}}
\rput(3.3,-1.6){\scalebox{1.5}{$s_i$}}
\end{pspicture}
}
}

\def\CaseAtwo{
\psscalebox{.5 .5} 
{
\begin{pspicture}[shift=-1.2](0,-1.2)(1.6,1.2)
\psframe[linecolor=black, linewidth=0.04, linestyle=dotted, dotsep=0.10583334cm, dimen=outer](1.6,1.2)(0.0,-1.2)
\psline[linecolor=black, linewidth=0.04](0.8,1.2)(0.8,-1.2)
\rput(.8,-1.6){\scalebox{2}{$s_i$}}
\end{pspicture}
}
}

\def\CaseAone{
\psscalebox{.5 .5} 
{
\begin{pspicture}[shift=-1.2](0,-1.2)(1.6,1.2)
\psframe[linecolor=black, linewidth=0.04, linestyle=dotted, dotsep=0.10583334cm, dimen=outer](1.6,1.2)(0.0,-1.2)
\psline[linecolor=black, linewidth=0.04, dotsize=0.07055555cm 2.0]{-*}(0.8,-1.2)(0.8,0.0)
\rput(.8,-1.6){\scalebox{2}{$s_i$}}
\end{pspicture}
}}

\def\exatreeUDB{   
\psscalebox{.5 .5} 
{
\begin{pspicture}[shift=-1.2](0,-1.2144687)(2.4,1.2144687)
\definecolor{colour2}{rgb}{1.0,0.0,0.2}
\definecolor{colour0}{rgb}{0.2,0.2,1.0}
\psframe[linecolor=black, linewidth=0.04, linestyle=dotted, dotsep=0.10583334cm, dimen=outer](2.4,1.1944687)(0.0,-1.2055314)
\psbezier[linecolor=colour2, linewidth=0.04](0.4,-1.2055314)(0.4,-0.40553132)(0.4,0.3944687)(0.8,0.39446868896484377)(1.2,0.3944687)(1.2,0.3944687)(1.6,0.3944687)(2.0,0.3944687)(2.0,-0.40553132)(2.0,-1.2055314)
\psline[linecolor=colour0, linewidth=0.04](1.2,-0.40553132)(0.8,-1.2055314)
\psline[linecolor=colour0, linewidth=0.04](1.2,-0.40553132)(1.6,-1.2055314)
\psline[linecolor=colour0, linewidth=0.04](1.2,-0.40553132)(1.2,1.1944687)(1.2,1.1944687)
\psline[linecolor=green, linewidth=0.04](1.2,-1.2055314)(1.2,-0.40553132)
\psline[linecolor=green, linewidth=0.04](1.2,-0.40553132)(1.6,0.3944687)(1.6,1.1944687)(1.6,1.1944687)(1.6,1.1944687)
\psline[linecolor=green, linewidth=0.04](1.2,-0.40553132)(0.8,0.3944687)(0.8,1.1944687)
\end{pspicture}
}
 }

\def\exatreeUDA{
\psscalebox{.5 .5} 
{
\begin{pspicture}[shift=-1.2](0,-1.2)(2.4,1.2)
\definecolor{colour0}{rgb}{0.2,0.2,1.0}
\definecolor{colour1}{rgb}{0.2,1.0,0.2}
\definecolor{colour2}{rgb}{1.0,0.0,0.2}
\psframe[linecolor=black, linewidth=0.04, linestyle=dotted, dotsep=0.10583334cm, dimen=outer](2.4,1.2)(0.0,-1.2)
\psline[linecolor=colour0, linewidth=0.04](0.8,1.2)(0.8,-1.2)
\psline[linecolor=colour0, linewidth=0.04](1.6,1.2)(1.6,-1.2)
\psline[linecolor=colour1, linewidth=0.04](1.2,1.2)(1.2,-1.2)
\psbezier[linecolor=colour2, linewidth=0.04](0.4,-1.2)(0.4,-0.4)(0.4,0.0)(0.8,0.0)(1.2,0.0)(1.2,0.0)(1.6,0.0)(2.0,0.0)(2.0,-0.4)(2.0,-1.2)
\end{pspicture}
}
}

\def\exatreeUD{
\psscalebox{.5 .5} 
{
\begin{pspicture}[shift=-1.2](0,-1.2)(2.0,1.2)
\definecolor{colour0}{rgb}{0.2,0.2,1.0}
\definecolor{colour1}{rgb}{0.2,1.0,0.2}
\psframe[linecolor=black, linewidth=0.04, linestyle=dotted, dotsep=0.10583334cm, dimen=outer](2.0,1.2)(0.0,-1.2)
\psline[linecolor=red, linewidth=0.04](0.4,1.2)(0.4,-1.2)
\psline[linecolor=colour0, linewidth=0.04](0.8,1.2)(0.8,-1.2)
\psline[linecolor=colour0, linewidth=0.04](1.6,1.2)(1.6,-1.2)
\psline[linecolor=colour1, linewidth=0.04](1.2,1.2)(1.2,-1.2)
\end{pspicture}
}
}

\def\Fibo{
\psscalebox{.5 .5} 
{
\begin{pspicture}[shift=-2.4](0,-2.4985573)(16.197115,2.4985573)
\psdots[linecolor=black, dotsize=0.2](0.09855774,-2.4000003)
\psline[linecolor=black, linewidth=0.04](0.09855774,-2.4000003)(0.09855774,-1.6000003)
\psline[linecolor=black, linewidth=0.04](0.8985577,-2.4000003)(1.2985578,-1.6000003)(1.6985577,-2.4000003)
\psline[linecolor=black, linewidth=0.04](2.4985578,-2.4000003)(2.8985577,-1.6000003)(3.2985578,-2.4000003)
\psline[linecolor=black, linewidth=0.04](4.098558,-1.6000003)(4.098558,-2.4000003)
\psline[linecolor=black, linewidth=0.04](4.8985577,-2.4000003)(5.2985578,-1.6000003)(5.698558,-2.4000003)
\psline[linecolor=black, linewidth=0.04](6.4985576,-2.4000003)(6.8985577,-1.6000003)(7.2985578,-2.4000003)(7.2985578,-2.4000003)
\psline[linecolor=black, linewidth=0.04](8.098557,-1.6000003)(8.098557,-2.4000003)
\psline[linecolor=black, linewidth=0.04](8.898558,-2.4000003)(9.298557,-1.6000003)(9.698558,-2.4000003)
\psline[linecolor=black, linewidth=0.04](10.498558,-1.6000003)(10.498558,-2.4000003)
\psline[linecolor=black, linewidth=0.04](11.298557,-2.4000003)(11.698558,-1.6000003)(12.098557,-2.4000003)
\psline[linecolor=black, linewidth=0.04](0.09855774,-1.6000003)(0.49855775,-0.8000003)(1.2985578,-1.6000003)(1.2985578,-1.6000003)
\psline[linecolor=black, linewidth=0.04](2.8985577,-1.6000003)(2.8985577,-0.8000003)
\psline[linecolor=black, linewidth=0.04](4.098558,-1.6000003)(4.4985576,-0.8000003)(5.2985578,-1.6000003)
\psline[linecolor=black, linewidth=0.04](6.8985577,-1.6000003)(6.8985577,-0.8000003)(6.8985577,-0.8000003)(6.8985577,-0.8000003)
\psline[linecolor=black, linewidth=0.04](8.098557,-1.6000003)(8.498558,-0.8000003)(9.298557,-1.6000003)
\psline[linecolor=black, linewidth=0.04](10.498558,-1.6000003)(10.898558,-0.8000003)(11.698558,-1.6000003)
\psline[linecolor=black, linewidth=0.04](0.49855775,-0.8000003)(0.49855775,0.0)
\psline[linecolor=black, linewidth=0.04](2.8985577,-0.8000003)(3.6985579,0.0)(4.4985576,-0.8000003)
\psline[linecolor=black, linewidth=0.04](6.8985577,-0.8000003)(7.698558,0.0)(8.498558,-0.8000003)
\psline[linecolor=black, linewidth=0.04](10.898558,-0.8000003)(10.898558,0.0)
\psline[linecolor=black, linewidth=0.04](12.898558,-2.4000003)(13.298557,-1.6000003)(13.698558,-2.4000003)
\psline[linecolor=black, linewidth=0.04](14.498558,-2.4000003)(14.498558,-1.6000003)
\psline[linecolor=black, linewidth=0.04](15.298557,-2.4000003)(15.698558,-1.6000003)(16.098558,-2.4000003)
\psline[linecolor=black, linewidth=0.04](13.298557,-1.6000003)(13.298557,-0.8000003)
\psline[linecolor=black, linewidth=0.04](14.498558,-1.6000003)(14.898558,-0.8000003)(15.698558,-1.6000003)
\psline[linecolor=black, linewidth=0.04](13.298557,-0.8000003)(14.098557,0.0)(14.898558,-0.8000003)
\psline[linecolor=black, linewidth=0.04](0.49855775,0.0)(2.0985577,0.7999997)(3.6985579,0.0)
\psline[linecolor=black, linewidth=0.04](7.698558,0.0)(7.698558,0.7999997)
\psline[linecolor=black, linewidth=0.04](10.898558,0.0)(12.498558,0.7999997)(14.098557,0.0)
\psline[linecolor=black, linewidth=0.04](2.0985577,0.7999997)(2.0985577,1.5999997)
\psline[linecolor=black, linewidth=0.04](7.698558,0.7999997)(10.098557,1.5999997)(12.498558,0.7999997)(12.498558,0.7999997)
\psline[linecolor=black, linewidth=0.04](2.0985577,1.5999997)(6.098558,2.3999996)(10.098557,1.5999997)
\psdots[linecolor=black, dotsize=0.2](0.8985577,-2.4000003)
\psdots[linecolor=black, dotsize=0.2](1.6985577,-2.4000003)
\psdots[linecolor=black, dotsize=0.2](2.4985578,-2.4000003)
\psdots[linecolor=black, dotsize=0.2](3.2985578,-2.4000003)
\psdots[linecolor=black, dotsize=0.2](4.098558,-2.4000003)
\psdots[linecolor=black, dotsize=0.2](4.8985577,-2.4000003)
\psdots[linecolor=black, dotsize=0.2](5.698558,-2.4000003)
\psdots[linecolor=black, dotsize=0.2](6.4985576,-2.4000003)
\psdots[linecolor=black, dotsize=0.2](7.2985578,-2.4000003)
\psdots[linecolor=black, dotsize=0.2](8.098557,-2.4000003)
\psdots[linecolor=black, dotsize=0.2](8.898558,-2.4000003)
\psdots[linecolor=black, dotsize=0.2](9.698558,-2.4000003)
\psdots[linecolor=black, dotsize=0.2](10.498558,-2.4000003)
\psdots[linecolor=black, dotsize=0.2](11.298557,-2.4000003)
\psdots[linecolor=black, dotsize=0.2](12.098557,-2.4000003)
\psdots[linecolor=black, dotsize=0.2](12.898558,-2.4000003)
\psdots[linecolor=black, dotsize=0.2](13.698558,-2.4000003)
\psdots[linecolor=black, dotsize=0.2](14.498558,-2.4000003)
\psdots[linecolor=black, dotsize=0.2](15.298557,-2.4000003)
\psdots[linecolor=black, dotsize=0.2](16.098558,-2.4000003)
\psdots[linecolor=black, dotsize=0.2](15.698558,-1.6000003)
\psdots[linecolor=black, dotsize=0.2](14.498558,-1.6000003)
\psdots[linecolor=black, dotsize=0.2](14.898558,-0.8000003)
\psdots[linecolor=black, dotsize=0.2](14.098557,0.0)
\psdots[linecolor=black, dotsize=0.2](13.298557,-0.8000003)
\psdots[linecolor=black, dotsize=0.2](13.298557,-1.6000003)
\psdots[linecolor=black, dotsize=0.2](11.698558,-1.6000003)
\psdots[linecolor=black, dotsize=0.2](10.898558,-0.8000003)
\psdots[linecolor=black, dotsize=0.2](10.498558,-1.6000003)
\psdots[linecolor=black, dotsize=0.2](9.298557,-1.6000003)
\psdots[linecolor=black, dotsize=0.2](8.098557,-1.6000003)
\psdots[linecolor=black, dotsize=0.2](8.498558,-0.8000003)
\psdots[linecolor=black, dotsize=0.2](7.698558,0.0)
\psdots[linecolor=black, dotsize=0.2](6.8985577,-0.8000003)
\psdots[linecolor=black, dotsize=0.2](6.8985577,-1.6000003)
\psdots[linecolor=black, dotsize=0.2](5.2985578,-1.6000003)
\psdots[linecolor=black, dotsize=0.2](4.098558,-1.6000003)
\psdots[linecolor=black, dotsize=0.2](4.4985576,-0.8000003)
\psdots[linecolor=black, dotsize=0.2](3.6985579,0.0)
\psdots[linecolor=black, dotsize=0.2](2.8985577,-0.8000003)
\psdots[linecolor=black, dotsize=0.2](2.8985577,-1.6000003)
\psdots[linecolor=black, dotsize=0.2](1.2985578,-1.6000003)
\psdots[linecolor=black, dotsize=0.2](0.49855775,-0.8000003)
\psdots[linecolor=black, dotsize=0.2](0.09855774,-1.6000003)
\psdots[linecolor=black, dotsize=0.2](0.49855775,0.0)
\psdots[linecolor=black, dotsize=0.2](2.0985577,0.7999997)
\psdots[linecolor=black, dotsize=0.2](2.0985577,1.5999997)
\psdots[linecolor=black, dotsize=0.2](6.098558,2.3999996)
\psdots[linecolor=black, dotsize=0.2](7.698558,0.7999997)
\psdots[linecolor=black, dotsize=0.2](10.098557,1.5999997)
\psdots[linecolor=black, dotsize=0.2](12.498558,0.7999997)
\psdots[linecolor=black, dotsize=0.2](10.898558,0.0)
\end{pspicture}
}
 }

\def\treeproofA{\psscalebox{.5 .5} 
{
\begin{pspicture}[shift=-2.2](0,-2.2)(4.8,2.2)
\definecolor{colour0}{rgb}{1.0,0.2,0.2}
\psframe[linecolor=black, linewidth=0.04, linestyle=dotted, dotsep=0.10583334cm, dimen=outer](3.2,2.2)(1.6,0.6)
\psframe[linecolor=black, linewidth=0.04, linestyle=dotted, dotsep=0.10583334cm, dimen=outer](1.6,-0.6)(0.0,-2.2)
\psframe[linecolor=black, linewidth=0.04, linestyle=dotted, dotsep=0.10583334cm, dimen=outer](4.8,-0.6)(3.2,-2.2)
\psline[linecolor=black, linewidth=0.04, arrowsize=0.05291667cm 2.0,arrowlength=1.4,arrowinset=0.0]{->}(1.6,1.4)(0.8,1.4)(0.8,-0.2)
\psline[linecolor=black, linewidth=0.04, arrowsize=0.05291667cm 2.0,arrowlength=1.4,arrowinset=0.0]{->}(3.2,1.4)(4.0,1.4)(4.0,-0.2)
\psline[linecolor=colour0, linewidth=0.04](1.2,-0.6)(1.2,-2.2)
\psline[linecolor=colour0, linewidth=0.04, dotsize=0.07055555cm 2.0]{-*}(4.4,-2.2)(4.4,-1.4)
\rput(2.4,1.5){\scalebox{3}{${l}$}}
\rput(2.4,0.2){\scalebox{1.5}{${\underline{v}'}$}}
\rput(2.4,2.4){\scalebox{1.5}{${u}$}}
\rput(0.6,-1.3){\scalebox{3}{${l}$}}
\rput(3.8,-1.3){\scalebox{3}{${l}$}}
\rput(3.9,-2.5){\scalebox{1.5}{${\underline{v}}$}}
\rput(3.9,-.4){\scalebox{1.5}{${u}$}}
\rput(.6,-2.5){\scalebox{1.5}{${\underline{v}}$}}
\rput(.6,-.4){\scalebox{1.5}{${us}$}}
\end{pspicture}
}
}

\def\treeproof{\psscalebox{.5 .5} 
{
\begin{pspicture}[shift=-2.2](0,-2.2)(4.8,2.2)
\definecolor{colour0}{rgb}{1.0,0.2,0.2}
\psframe[linecolor=black, linewidth=0.04, linestyle=dotted, dotsep=0.10583334cm, dimen=outer](3.2,2.2)(1.6,0.6)
\psframe[linecolor=black, linewidth=0.04, linestyle=dotted, dotsep=0.10583334cm, dimen=outer](1.6,-0.6)(0.0,-2.2)
\psframe[linecolor=black, linewidth=0.04, linestyle=dotted, dotsep=0.10583334cm, dimen=outer](4.8,-0.6)(3.2,-2.2)
\psline[linecolor=black, linewidth=0.04, arrowsize=0.05291667cm 2.0,arrowlength=1.4,arrowinset=0.0]{->}(1.6,1.4)(0.8,1.4)(0.8,-0.2)
\psline[linecolor=black, linewidth=0.04, arrowsize=0.05291667cm 2.0,arrowlength=1.4,arrowinset=0.0]{->}(3.2,1.4)(4.0,1.4)(4.0,-0.2)
\psline[linecolor=colour0, linewidth=0.04](1.2,-0.6)(1.2,-2.2)
\psline[linecolor=colour0, linewidth=0.04, dotsize=0.07055555cm 2.0]{-*}(4.4,-2.2)(4.4,-1.4)
\rput(2.4,1.5){\scalebox{3}{${l}$}}
\rput(2.4,0.2){\scalebox{1.5}{${\underline{v}'}$}}
\rput(2.4,2.4){\scalebox{1.5}{${us}$}}
\rput(0.6,-1.3){\scalebox{3}{${l}$}}
\rput(3.8,-1.3){\scalebox{3}{${l}$}}
\rput(3.9,-2.5){\scalebox{1.5}{${\underline{v}}$}}
\rput(3.9,-.4){\scalebox{1.5}{${us}$}}
\rput(.6,-2.5){\scalebox{1.5}{${\underline{v}}$}}
\rput(.6,-.4){\scalebox{1.5}{${u}$}}
\end{pspicture}
}
}

\def\exampletreeB{ \psscalebox{.5 .5} 
{
\begin{pspicture}[shift=-1.2](0,-1.2)(15.6,1.2)
\definecolor{colour2}{rgb}{0.2,0.2,1.0}
\definecolor{colour0}{rgb}{1.0,0.2,0.2}
\definecolor{colour3}{rgb}{0.2,1.0,0.2}
\psframe[linecolor=black, linewidth=0.04, linestyle=dotted, dotsep=0.10583334cm, dimen=outer](2.8,1.2)(0.0,-1.2)
\psframe[linecolor=black, linewidth=0.04, linestyle=dotted, dotsep=0.10583334cm, dimen=outer](6.0,1.2)(3.2,-1.2)
\psframe[linecolor=black, linewidth=0.04, linestyle=dotted, dotsep=0.10583334cm, dimen=outer](9.2,1.2)(6.4,-1.2)
\psframe[linecolor=black, linewidth=0.04, linestyle=dotted, dotsep=0.10583334cm, dimen=outer](12.4,1.2)(9.6,-1.2)
\psframe[linecolor=black, linewidth=0.04, linestyle=dotted, dotsep=0.10583334cm, dimen=outer](15.6,1.2)(12.8,-1.2)
\psline[linecolor=colour2, linewidth=0.04, dotsize=0.07055555cm 2.0]{-*}(0.4,-1.2)(0.4,-0.4)(0.4,-0.4)
\psline[linecolor=colour2, linewidth=0.04, dotsize=0.07055555cm 2.0]{-*}(1.6,-1.2)(1.6,-0.4)
\psline[linecolor=colour2, linewidth=0.04, dotsize=0.07055555cm 2.0]{-*}(6.8,-1.2)(6.8,-0.4)
\psline[linecolor=colour2, linewidth=0.04, dotsize=0.07055555cm 2.0]{-*}(8.0,-1.2)(8.0,-0.4)(8.0,-0.4)
\psline[linecolor=colour2, linewidth=0.04, dotsize=0.07055555cm 2.0]{-*}(10.0,-1.2)(10.0,-0.4)
\psline[linecolor=colour2, linewidth=0.04, dotsize=0.07055555cm 2.0]{-*}(11.2,-1.2)(11.2,-0.4)
\psline[linecolor=colour2, linewidth=0.04, dotsize=0.07055555cm 2.0]{-*}(13.2,-1.2)(13.2,-0.4)
\psline[linecolor=colour2, linewidth=0.04, dotsize=0.07055555cm 2.0]{-*}(14.4,-1.2)(14.4,-0.4)(14.4,-0.4)
\psline[linecolor=colour0, linewidth=0.04, dotsize=0.07055555cm 2.0]{-*}(1.2,-1.2)(1.2,-0.4)
\psline[linecolor=colour0, linewidth=0.04, dotsize=0.07055555cm 2.0]{-*}(2.4,-1.2)(2.4,-0.4)
\psline[linecolor=colour0, linewidth=0.04, dotsize=0.07055555cm 2.0]{-*}(4.4,-1.2)(4.4,-0.4)
\psline[linecolor=colour0, linewidth=0.04, dotsize=0.07055555cm 2.0]{-*}(5.6,-1.2)(5.6,-0.4)
\psline[linecolor=colour0, linewidth=0.04, dotsize=0.07055555cm 2.0]{-*}(7.6,-1.2)(7.6,-0.4)
\psline[linecolor=colour0, linewidth=0.04, dotsize=0.07055555cm 2.0]{-*}(8.8,-1.2)(8.8,-0.4)
\psline[linecolor=colour3, linewidth=0.04, dotsize=0.07055555cm 2.0]{-*}(0.8,-1.2)(0.8,-0.4)
\psline[linecolor=colour3, linewidth=0.04, dotsize=0.07055555cm 2.0]{-*}(2.0,-1.2)(2.0,-0.4)
\psline[linecolor=colour3, linewidth=0.04, dotsize=0.07055555cm 2.0]{-*}(4.0,-1.2)(4.0,-0.4)
\psline[linecolor=colour3, linewidth=0.04, dotsize=0.07055555cm 2.0]{-*}(5.2,-1.2)(5.2,-0.4)
\psline[linecolor=colour3, linewidth=0.04, dotsize=0.07055555cm 2.0]{-*}(10.4,-1.2)(10.4,-0.4)
\psline[linecolor=colour3, linewidth=0.04, dotsize=0.07055555cm 2.0]{-*}(11.6,-1.2)(11.6,-0.4)
\psbezier[linecolor=colour2, linewidth=0.04](3.6,-1.2)(3.6,-0.4)(3.6,-0.4)(3.6,0.0)(3.6,0.4)(4.8,0.4)(4.8,0.0)(4.8,-0.4)(4.8,-0.4)(4.8,-1.2)
\psbezier[linecolor=colour0, linewidth=0.04](10.8,-1.2)(10.8,-0.4)(10.8,-0.4)(10.8,0.0)(10.8,0.4)(12.0,0.4)(12.0,0.0)(12.0,-0.4)(12.0,-0.4)(12.0,-1.2)
\psbezier[linecolor=colour0, linewidth=0.04](14.0,-1.2)(14.0,-0.4)(14.0,-0.4)(14.0,0.0)(14.0,0.4)(15.2,0.4)(15.2,0.0)(15.2,-0.4)(15.2,-0.4)(15.2,-1.2)
\psbezier[linecolor=colour3, linewidth=0.04](7.2,-1.2)(7.2,-0.4)(7.2,-0.4)(7.2,0.0)(7.2,0.4)(8.4,0.4)(8.4,0.0)(8.4,-0.4)(8.4,-0.4)(8.4,-1.2)
\psbezier[linecolor=colour3, linewidth=0.04](13.6,-1.2)(13.6,-0.4)(13.6,-0.4)(13.6,0.0)(13.6,0.4)(14.8,0.4)(14.8,0.0)(14.8,-0.4)(14.8,-0.4)(14.8,-1.2)
\end{pspicture}
}
}

\def\exampletreeA{ \psscalebox{.5 .5} 
{
\begin{pspicture}[shift=-1.2](0,-1.2025417)(15.6,1.2025417)
\definecolor{colour1}{rgb}{1.0,0.0,0.2}
\definecolor{colour2}{rgb}{0.2,0.2,1.0}
\definecolor{colour3}{rgb}{0.2,1.0,0.2}
\definecolor{colour0}{rgb}{1.0,0.2,0.2}
\psframe[linecolor=black, linewidth=0.04, linestyle=dotted, dotsep=0.10583334cm, dimen=outer](2.8,1.2025417)(0.0,-1.1974583)
\psframe[linecolor=black, linewidth=0.04, linestyle=dotted, dotsep=0.10583334cm, dimen=outer](6.0,1.2025417)(3.2,-1.1974583)
\psframe[linecolor=black, linewidth=0.04, linestyle=dotted, dotsep=0.10583334cm, dimen=outer](9.2,1.2025417)(6.4,-1.1974583)
\psframe[linecolor=black, linewidth=0.04, linestyle=dotted, dotsep=0.10583334cm, dimen=outer](12.4,1.2025417)(9.6,-1.1974583)
\psframe[linecolor=black, linewidth=0.04, linestyle=dotted, dotsep=0.10583334cm, dimen=outer](15.6,1.2025417)(12.8,-1.1974583)
\psline[linecolor=colour1, linewidth=0.04, dotsize=0.07055555cm 2.0]{-*}(0.4,-1.1974583)(0.4,-0.3974583)
\psline[linecolor=colour1, linewidth=0.04, dotsize=0.07055555cm 2.0]{-*}(1.2,-1.1974583)(1.2,-0.3974583)
\psline[linecolor=colour1, linewidth=0.04, dotsize=0.07055555cm 2.0]{-*}(2.4,-1.1974583)(2.4,-0.3974583)
\psline[linecolor=colour1, linewidth=0.04, dotsize=0.07055555cm 2.0]{-*}(5.6,-1.1974583)(5.6,-0.3974583)
\psline[linecolor=colour1, linewidth=0.04, dotsize=0.07055555cm 2.0]{-*}(6.8,-1.1974583)(6.8,-0.3974583)
\psline[linecolor=colour1, linewidth=0.04, dotsize=0.07055555cm 2.0]{-*}(8.8,-1.1974583)(8.8,-0.3974583)
\psline[linecolor=colour1, linewidth=0.04, dotsize=0.07055555cm 2.0]{-*}(10.0,-1.1974583)(10.0,-0.3974583)
\psline[linecolor=colour1, linewidth=0.04, dotsize=0.07055555cm 2.0]{-*}(14.0,-1.1974583)(14.0,-0.3974583)
\psline[linecolor=colour2, linewidth=0.04, dotsize=0.07055555cm 2.0]{-*}(0.8,-1.1974583)(0.8,-0.3974583)
\psline[linecolor=colour2, linewidth=0.04, dotsize=0.07055555cm 2.0]{-*}(2.0,-1.1974583)(2.0,-0.3974583)
\psline[linecolor=colour2, linewidth=0.04, dotsize=0.07055555cm 2.0]{-*}(4.0,-1.1974583)(4.0,-0.3974583)
\psline[linecolor=colour2, linewidth=0.04, dotsize=0.07055555cm 2.0]{-*}(5.2,-1.1974583)(5.2,-0.3974583)
\psline[linecolor=colour2, linewidth=0.04, dotsize=0.07055555cm 2.0]{-*}(10.4,-1.1974583)(10.4,-0.3974583)
\psline[linecolor=colour2, linewidth=0.04, dotsize=0.07055555cm 2.0]{-*}(11.6,-1.1974583)(11.6,-0.3974583)
\psline[linecolor=colour3, linewidth=0.04, dotsize=0.07055555cm 2.0]{-*}(1.6,-1.1974583)(1.6,-0.3974583)
\psline[linecolor=colour3, linewidth=0.04, dotsize=0.07055555cm 2.0]{-*}(4.8,-1.1974583)(4.8,-0.3974583)
\psline[linecolor=colour3, linewidth=0.04, dotsize=0.07055555cm 2.0]{-*}(8.0,-1.1974583)(8.0,-0.3974583)
\psline[linecolor=colour3, linewidth=0.04, dotsize=0.07055555cm 2.0]{-*}(11.2,-1.1974583)(11.2,-0.3974583)
\psline[linecolor=colour3, linewidth=0.04, dotsize=0.07055555cm 2.0]{-*}(14.4,-1.1974583)(14.4,-0.3974583)
\psbezier[linecolor=colour0, linewidth=0.04](3.6053598,-1.1975927)(3.357177,-0.22887951)(3.6,0.4025417)(3.9960327,0.40471016928805154)(4.3920655,0.40687865)(4.6422267,-0.22139378)(4.4053464,-1.1929326)
\psbezier[linecolor=colour0, linewidth=0.04](10.8,-1.1974583)(10.8,-0.3974583)(10.8,-0.3974583)(10.8,0.002541694641113281)(10.8,0.4025417)(12.0,0.4025417)(12.0,0.0025416946)(12.0,-0.3974583)(12.0,-0.3974583)(12.0,-1.1974583)
\psbezier[linecolor=colour0, linewidth=0.04](13.2,-1.1974583)(13.2,-0.5974583)(13.2,-0.2974583)(13.2,0.002541694641113281)(13.2,0.3025417)(15.2,0.4025417)(15.2,0.0025416946)(15.2,-0.3974583)(15.2,-0.5974583)(15.2,-1.1974583)
\psline[linecolor=colour0, linewidth=0.04, dotsize=0.07055555cm 2.0]{-*}(7.6,-1.1974583)(7.6,-0.3974583)
\psbezier[linecolor=colour2, linewidth=0.04](7.2,-1.1974583)(7.2,-0.5974583)(7.2,-0.5974583)(7.2,0.002541694641113281)(7.2,0.6025417)(8.4,0.6025417)(8.4,0.0025416946)(8.4,-0.5974583)(8.4,-0.5974583)(8.4,-1.1974583)
\psbezier[linecolor=colour2, linewidth=0.04](13.6,-1.1974583)(13.6,-0.3974583)(13.6,-0.7974583)(13.6,-0.3974583053588867)(13.6,0.0025416946)(14.8,0.0025416946)(14.8,-0.3974583)(14.8,-0.7974583)(14.8,-0.7974583)(14.8,-1.1974583)
\end{pspicture}
}}

\def\treesteptwo{\psscalebox{.5 .5} 
{
\begin{pspicture}[shift=-3](0,-3.0)(4.42,3.0)
\definecolor{colour1}{rgb}{0.2,0.2,1.0}
\definecolor{colour2}{rgb}{0.8,0.8,0.8}
\psframe[linecolor=black, linewidth=0.04, linestyle=dotted, dotsep=0.10583334cm, dimen=outer](3.22,3.0)(1.62,1.4)
\psframe[linecolor=black, linewidth=0.04, linestyle=dotted, dotsep=0.10583334cm, dimen=outer](4.42,0.2)(0.02,-3.0)
\psframe[linecolor=black, linewidth=0.04, linestyle=dotted, dotsep=0.10583334cm, dimen=outer](3.22,-1.4)(0.02,-3.0)
\psframe[linecolor=black, linewidth=0.04, linestyle=dotted, dotsep=0.10583334cm, dimen=outer](3.24,-1.4)(3.22,-1.42)
\psline[linecolor=black, linewidth=0.04, linestyle=dotted, dotsep=0.10583334cm](0.02,-0.6)(3.22,-0.6)(3.22,-1.4)
\psline[linecolor=colour1, linewidth=0.04](1.22,-1.4)(1.22,-1.8)
\psline[linecolor=black, linewidth=0.04](0.82,-1.4)(0.82,-1.8)
\psline[linecolor=black, linewidth=0.04](0.42,-1.4)(0.42,-1.8)
\psline[linecolor=black, linewidth=0.04](1.62,-1.4)(1.62,-1.8)
\psline[linecolor=black, linewidth=0.04](2.02,-1.4)(2.02,-1.8)
\psline[linecolor=black, linewidth=0.04](2.42,-1.4)(2.42,-1.8)
\psline[linecolor=black, linewidth=0.04](2.82,-1.4)(2.82,-1.8)
\psline[linecolor=black, linewidth=0.04](0.42,0.2)(0.42,-0.6)
\psline[linecolor=black, linewidth=0.04](0.82,0.2)(0.82,-0.6)
\psline[linecolor=black, linewidth=0.04](1.22,0.2)(1.22,-0.6)
\psline[linecolor=black, linewidth=0.04](1.62,0.2)(1.62,-0.6)
\psline[linecolor=black, linewidth=0.04](2.02,0.2)(2.02,-0.6)
\psline[linecolor=black, linewidth=0.04](2.42,0.2)(2.42,-0.6)
\psbezier[linecolor=colour1, linewidth=0.04](2.82,-0.6)(2.82,-0.2)(4.02,0.6)(4.02,-3.0)
\psline[linecolor=black, linewidth=0.04, arrowsize=0.05291667cm 2.0,arrowlength=1.4,arrowinset=0.0]{->}(2.42,1.4)(2.42,0.2)
\psframe[linecolor=black, linewidth=0.04, fillstyle=solid,fillcolor=colour2, dimen=outer](3.22,-0.6)(0.02,-1.4)
\rput(2.5,2.3){\scalebox{3}{$D$}}
\rput(1.6,-2.3){\scalebox{3}{$D$}}
\end{pspicture}
}
}

\def\treestepone{\psscalebox{.5 .5} 
{
\begin{pspicture}[shift=-1.8](0,-1.8)(4.8,1.8)
\definecolor{colour1}{rgb}{0.2,0.2,1.0}
\psframe[linecolor=black, linewidth=0.04, linestyle=dotted, dotsep=0.10583334cm, dimen=outer](2.0,-0.2)(0.0,-1.8)
\psframe[linecolor=black, linewidth=0.04, linestyle=dotted, dotsep=0.10583334cm, dimen=outer](5.2,-0.2)(3.2,-1.8)
\psframe[linecolor=black, linewidth=0.04, linestyle=dotted, dotsep=0.10583334cm, dimen=outer](3.6,1.8)(2.0,0.2)
\psline[linecolor=black, linewidth=0.04, arrowsize=0.05291667cm 2.0,arrowlength=1.4,arrowinset=0.0]{->}(2.0,1.0)(1.2,1.0)(1.2,-0.2)
\psline[linecolor=black, linewidth=0.04, arrowsize=0.05291667cm 2.0,arrowlength=1.4,arrowinset=0.0]{->}(3.6,1.0)(4.4,1.0)(4.4,-0.2)
\psline[linecolor=colour1, linewidth=0.04](1.6,-0.2)(1.6,-1.8)
\psline[linecolor=colour1, linewidth=0.04](4.8,-1.0)(4.8,-1.8)
\psdots[linecolor=colour1, dotsize=0.2](4.8,-1.0)
\rput(2.8,1){\scalebox{3}{$D$}}
\rput(0.8,-1){\scalebox{3}{$D$}}
\rput(4,-1){\scalebox{3}{$D$}}
\end{pspicture}
}
}

\def\univalent{\psscalebox{.5 .5} 
{
\begin{pspicture}[shift=-1.2](0,-1.2)(2.4,1.2)
\definecolor{colour0}{rgb}{1.0,0.2,0.2}
\pscircle[linecolor=black, linewidth=0.04, linestyle=dotted, dotsep=0.10583334cm, dimen=outer](1.2,0.0){1.2}
\psline[linecolor=colour0, linewidth=0.04](1.2,0.0)(1.2,-1.2)
\psdots[linecolor=colour0, dotsize=0.2](1.2,0.0)
\end{pspicture}
}}

\def\multivalent{\psscalebox{.5 .5} 
{
\begin{pspicture}[shift=-1.2](0,-1.2)(2.4,1.2)
\definecolor{colour0}{rgb}{1.0,0.2,0.2}
\definecolor{colour1}{rgb}{0.2,0.2,1.0}
\pscircle[linecolor=black, linewidth=0.04, linestyle=dotted, dotsep=0.10583334cm, dimen=outer](1.2,0.0){1.2}
\psline[linecolor=colour0, linewidth=0.04](1.2,0.0)(1.2,-1.2)
\psline[linecolor=colour0, linewidth=0.04](1.2,0.0)(0.0,0.0)
\psline[linecolor=colour0, linewidth=0.04](1.2,0.0)(1.2,1.2)
\psline[linecolor=colour0, linewidth=0.04](1.2,0.0)(2.4,0.0)
\psline[linecolor=colour1, linewidth=0.04](1.2,0.0)(0.4,0.8)
\psline[linecolor=colour1, linewidth=0.04](1.2,0.0)(2.0,0.8)
\psline[linecolor=colour1, linewidth=0.04](1.2,0.0)(2.0,-0.8)
\psline[linecolor=colour1, linewidth=0.04](1.2,0.0)(0.4,-0.8)
\end{pspicture}
}}

\def\exampleSgraph{\psscalebox{.5 .5} 
{
\begin{pspicture}[shift=-2.7](0,-2.607071)(6.434142,2.607071)
\definecolor{colour0}{rgb}{1.0,0.2,0.2}
\definecolor{colour2}{rgb}{0.2,1.0,0.2}
\definecolor{colour1}{rgb}{0.2,0.2,1.0}
\psline[linecolor=colour0, linewidth=0.04](0.014142135,-2.5929291)(0.8141421,-1.792929)(0.8141421,-0.9929291)
\psline[linecolor=colour0, linewidth=0.04](0.8141421,-1.792929)(1.6141422,-2.5929291)
\psbezier[linecolor=colour0, linewidth=0.04](0.8141421,-0.9929291)(0.8141421,-0.79292905)(0.8141421,-0.19292907)(1.6141422,-0.5929290771484375)(2.4141421,-0.9929291)(4.014142,-1.3929291)(4.814142,-0.19292907)
\psline[linecolor=colour2, linewidth=0.04](3.214142,-2.5929291)(3.214142,-0.19292907)
\psdots[linecolor=colour2, dotsize=0.2](3.214142,-0.19292907)
\psline[linecolor=colour1, linewidth=0.04](0.8141421,-1.792929)(0.8141421,-2.5929291)
\psbezier[linecolor=colour1, linewidth=0.04](0.79152906,-1.8271714)(1.1691972,-1.4060221)(1.6141422,-1.792929)(1.6141422,-1.7929290771484374)(1.6141422,-1.792929)(2.014142,-2.192929)(2.0332365,-2.5607672)
\psbezier[linecolor=colour1, linewidth=0.04](0.8141421,-1.792929)(0.41414213,-1.792929)(0.014142135,-0.59292907)(0.41414213,-0.1929290771484375)(0.8141421,0.20707092)(1.6141422,-0.19292907)(2.014142,0.20707092)(2.4141421,0.6070709)(2.4141421,0.6070709)(2.8141422,0.6070709)(3.214142,0.6070709)(2.8141422,1.0070709)(3.214142,1.0070709)
\psline[linecolor=colour1, linewidth=0.04](4.814142,-2.5929291)(4.814142,-0.19292907)
\psbezier[linecolor=colour1, linewidth=0.04](4.814142,-0.19292907)(5.2141423,0.6070709)(5.614142,0.20707092)(5.614142,0.6070709228515625)(5.614142,1.0070709)(5.614142,1.807071)(5.614142,2.607071)
\psbezier[linecolor=colour1, linewidth=0.04](3.214142,1.0070709)(3.214142,1.0070709)(3.6141422,0.6070709)(4.014142,1.0070709228515624)(4.414142,1.4070709)(4.814142,-0.19292907)(4.814142,-0.19292907)
\psbezier[linecolor=colour0, linewidth=0.04](4.814142,-0.19292907)(5.2141423,-0.19292907)(5.2141423,-0.19292907)(5.614142,-0.5929290771484375)(6.014142,-0.9929291)(6.414142,-2.192929)(6.414142,-2.5929291)
\psline[linecolor=colour0, linewidth=0.04](5.614142,-1.792929)(5.614142,-2.5929291)
\psdots[linecolor=colour0, dotsize=0.2](5.614142,-1.792929)
\psline[linecolor=colour0, linewidth=0.04](4.814142,-0.19292907)(4.814142,2.607071)
\psbezier[linecolor=colour2, linewidth=0.04](4.014142,-2.5929291)(3.6141422,-1.792929)(4.2566776,-1.1630716)(4.014142,-0.1929290771484375)(3.7716064,0.7772134)(3.6141422,-0.19292907)(3.214142,1.0070709)
\psbezier[linecolor=colour2, linewidth=0.04](3.214142,1.0070709)(3.6141422,1.4070709)(4.014142,1.807071)(4.014142,1.8070709228515625)(4.014142,1.807071)(4.014142,2.607071)(4.014142,2.607071)
\psbezier[linecolor=colour2, linewidth=0.04](3.214142,1.0070709)(2.4141421,1.4070709)(2.8141422,1.4070709)(2.4141421,1.8070709228515625)(2.014142,2.2070708)(1.6141422,1.807071)(1.6141422,1.4070709)
\psdots[linecolor=colour2, dotsize=0.2](1.6141422,1.4070709)
\psline[linecolor=colour1, linewidth=0.04](3.214142,1.0070709)(3.214142,1.807071)
\psdots[linecolor=colour1, dotsize=0.2](3.214142,1.807071)
\end{pspicture}
}
}

\def\levelone{\psscalebox{.5 .5} 
{
\begin{pspicture}[shift=-1.8](0,-1.8)(4.8,1.8)
\definecolor{colour0}{rgb}{1.0,0.2,0.2}
\psframe[linecolor=black, linewidth=0.04, linestyle=dotted, dotsep=0.10583334cm, dimen=outer](1.6,-0.2)(0.0,-1.8)
\psframe[linecolor=black, linewidth=0.04, linestyle=dotted, dotsep=0.10583334cm, dimen=outer](4.8,-0.2)(3.2,-1.8)
\psframe[linecolor=black, linewidth=0.04, linestyle=dotted, dotsep=0.10583334cm, dimen=outer](3.2,1.8)(1.6,0.2)
\psline[linecolor=black, linewidth=0.04, arrowsize=0.05291667cm 2.0,arrowlength=1.4,arrowinset=0.0]{->}(1.6,1.0)(0.8,1.0)(0.8,-0.2)
\psline[linecolor=black, linewidth=0.04, arrowsize=0.05291667cm 2.0,arrowlength=1.4,arrowinset=0.0]{->}(3.2,1.0)(4.0,1.0)(4.0,-0.2)
\psline[linecolor=colour0, linewidth=0.04](0.8,-0.2)(0.8,-1.8)
\psline[linecolor=colour0, linewidth=0.04](4.0,-1.0)(4.0,-1.8)
\psdots[linecolor=colour0, dotsize=0.2](4.0,-1.0)
\rput(2.4,1){\scalebox{3}{$\emptyset$}}
\end{pspicture}
}
}